\RequirePackage{fix-cm}
\documentclass[smallextended]{svjour3}       
\smartqed  

\usepackage{graphicx}          
\usepackage{amsmath,amssymb,mathtools,amsfonts}

\usepackage{hyperref}
\usepackage{multirow}
\usepackage{color}
\usepackage{psfrag}
\usepackage{subfigure}
\usepackage{empheq}
\usepackage{tabularx}
\usepackage{bm} 
\usepackage{psfrag}
\usepackage{float}
\usepackage{comment}
\usepackage{alphalph}

\usepackage[short]{optidef}

\usepackage[shortlabels]{enumitem}

\usepackage{tikz}
\usepackage{pgfplots}
\usetikzlibrary{arrows.meta}
\newlength\fwidth
\newlength\fheight
\usetikzlibrary{arrows}

\newcommand{\R}{{\mathbb{R}}}
\newcommand{\dd}{{\mathrm{d}}}
\newcommand{\partialder}[2]{\frac{\partial{#1}}{\partial{#2}}}
\newcommand{\sign}{{\mathrm{sign}}}

\newcommand{\Nsys}{{n_{f}}} 
\newcommand{\Nstg}{{n_\mathrm{s}}} 
\newcommand{\Nsubsys}{{n_\mathrm{sys}}} 
\newcommand{\NFE}{{N_\mathrm{FE}}} 
\newcommand{\Nswitch}{{N_\mathrm{sw}}} 

\newcommand{\I}{{\mathcal{I}}}

\newcommand{\LP}{{\mathrm{LP}}}
\newcommand{\irk}{{\mathrm{rk}}}

\newcommand{\fesd}{{\mathrm{fesd}}}
\newcommand{\Tctrl}{{T_{\mathrm{ctrl}}}}
\newcommand{\Nctrl}{N_{\mathrm{ctrl}}}

\newcommand{\ts}{{t_{\mathrm{s}}}}
\newcommand{\tsn}{t_{\mathrm{s},n}}
\newcommand{\tsnn}{{t_{\mathrm{s},n+1}}}
\newcommand{\tsnhat}{{\hat{t}_{\mathrm{s},n}}}
\newcommand{\tsnnhat}{{\hat{t}_{\mathrm{s},n+1}}}
\newcommand{\hatts}{{\hat{t}_{\mathrm{s}}}}

\newtheorem{theorem1}{Theorem}
\newtheorem{example1}{Example}
\newtheorem{remark1}[theorem1]{Remark}
\newtheorem{definition1}[theorem1]{Definition}
\newtheorem{assumption}[theorem1]{Assumption}

\newtheorem{proposition1}[theorem1]{Proposition}

\newtheorem{lemma1}[theorem1]{Lemma}
\allowdisplaybreaks
\begin{document}
\title{Finite Elements with Switch Detection for Direct Optimal Control of Nonsmooth Systems}

\titlerunning{FESD for Direct Optimal Control of Nonsmooth Systems}        

\author{Armin Nurkanovi\'c         
		\and  Mario Sperl
       \and  Sebastian Albrecht 
       \and  Moritz Diehl 
}


\institute{A. Nurkanovi\'c, \at
			Department of Microsystems Engineering (IMTEK), University of Freiburg, Germany,
              \email{armin.nurkanovic@imtek.uni-freiburg.de}           
              \and
                M. Sperl,\at
                Department of Mathematics, Chair of Applied Mathematics, University of Bayreuth, Germany,     
                \email{Mario.Sperl@uni-bayreuth.de}           
               \and
              S. Albrecht,\at
				Siemens Technology, Munich, Germany,	\email{sebastian.albrecht@siemens.com} 
               \and
              M. Diehl, \at
			  Department of Microsystems Engineering (IMTEK) and Department of Mathematics, University of Freiburg, Germany,  \email{moritz.diehl@imtek.uni-freiburg.de}
              \and 
}          
\date{Received: xxxx/ Accepted: xxxx}

\maketitle
\begin{abstract}
\sloppy{This paper introduces Finite Elements with Switch Detection (FESD)}, a numerical discretization method for nonsmooth differential equations.
We consider the Filippov convexification of these systems and a transformation into dynamic complementarity systems introduced by Stewart \cite{Stewart1990b}.
FESD is based on solving nonlinear complementarity problems and can automatically detect nonsmooth events in time. 
If standard \textcolor{black}{time-stepping} Runge-Kutta (RK) methods are naively applied to a nonsmooth ODE, the accuracy is at best of order one.
In FESD, we let the integrator step size be a degree of freedom. Additional complementarity conditions, which we call \textit{cross complementarities}, enable \textit{exact} switch detection, hence FESD can recover the high order accuracy that the RK methods enjoy for smooth ODE. 
Additional conditions called \textit{step equilibration} allow the step size to change only when switches occur and thus avoid spurious degrees of freedom.
Convergence results for the FESD method are derived, local uniqueness of the solution and convergence of numerical sensitivities are proven. 
The efficacy of FESD is demonstrated in several simulation and optimal control examples. 
In an optimal control problem benchmark with FESD, we achieve up to five orders of magnitude more accurate solutions than a standard \textcolor{black}{time-stepping} approach for the same computational time.

 \end{abstract}

\section{Introduction}\label{sec:introduction}
The goal of this paper is to develop high-accuracy numerical simulation and optimal control methods for Ordinary Differential Equations (ODE) with a discontinuous vector field. We assume the following structure
	\begin{align} \label{eq:pws1}
	\begin{split}
		\dot{x}(t) & =f_i(x(t),u(t)),\ \mathrm{if} \  x(t) \in R_i \subset \R^{n_x},\ i \in \mathcal{J} \coloneqq \{ 1,\dots,\Nsys  \},
	\end{split}
\end{align}
where $R_i$ are disjoint open sets and $f_i(\cdot)$ are smooth functions on an open neighborhood of $\overline{R}_i$, $\Nsys$ is a positive integer and $u(t)$ is an externally chosen control function.

This formulation of \textit{piecewise smooth} ODE falls into the class of \textit{hybrid systems} \cite{Brogliato2020,Stewart2011}. 
Many practical problems give rise to such ODE with a discontinuous right hand side (r.h.s.), e.g., in sliding mode control \cite{Acary2010b}, mechanics problems with Coulomb friction \cite{Stewart2000}, state constrained ODE derived from Pontryagin's maximum principle \cite{Pontryagin1962a}, electronic circuits \cite{Acary2010}, biological systems \cite{Acary2014}, vaccination strategies \cite{Cojocaru2008}, transportation systems and traffic flow networks \cite{Ban2012}, constrained optimization algorithms viewed as dynamic systems \cite{Hauswirth2021} and many more.
Systems with state jumps, including impact mechanics, robotics and hybrid systems with hysteresis can be transformed into systems matching the form of \eqref{eq:pws1} via the \textit{time-freezing reformulation} \cite{Nurkanovic2023a,Nurkanovic2022a,Nurkanovic2021}.
Consequently, efficient and accurate numerical optimal control algorithms for this class of systems are of great interest. 

\paragraph{Related work:}
High-accuracy simulation of ODE with a discontinuous r.h.s. is numerically difficult since an accurate location of the nonsmooth events in time is needed. 
Standard \textcolor{black}{time-stepping methods} with a fixed step size applied to this class of ODE have at best first-order accuracy as they do not detect the switches \cite{Acary2008}.
Additionally, in the case of sliding modes the numerical solution obtained by \textcolor{black}{explicit time-stepping methods} tends to \textit{chatter} around the discontinuity.
Convergence of standard \textcolor{black}{time-stepping} discretization methods with order one is studied in \cite{Dontchev1992,Kastner1990,Taubert1981}. 
Most high-accuracy methods include a root-finding procedure for accurate switch location and they usually assume that the trajectory crosses the discontinuity or only passes through two regions at a time \cite{Acary2008}. 
An exception is Stewart's high accuracy method which can deal with almost all switching cases \cite{Stewart1990b,Stewart1996a}. 
However, it also has an \textit{external} switch location routine and is thus difficult to apply in direct optimal control, i.e., in \textit{first discretize, then optimize} approaches to optimal control \cite{Rawlings2017}.

The parametric sensitivities of discontinuous ODE \eqref{eq:pws1} have jump discontinuities \cite{Filippov1988,Stewart2010} when the trajectory passes through or enters a surface of discontinuity, cf. Section \ref{sec:sensitivites}.
Therefore, the application of high-accuracy integrators even in a direct multiple shooting setting \cite{Bock1984}, which hides the external switch detection procedure from the optimizer, will be notoriously difficult, since derivative-based optimization algorithms will likely fail due to the non-Lipschitz sensitivities. 
Other fundamental difficulties within direct optimal control of non-smooth ODE are illuminated in the seminal paper by Stewart and Anitescu \cite{Stewart2010}. 
They show that {direct methods based on time-stepping integration schemes with fixed step sizes} are doomed to fail since the numerical sensitivities obtained by differentiating the results of a simulation are wrong no matter how small the step size is. 
{We refer to this class of methods as standard methods because they are commonly used in optimal control of smooth dynamical systems.}
It is also shown that the numerical sensitivities of a smooth approximation of an ODE with a discontinuous r.h.s. are only correct if the step size approaches zero faster than the smoothing parameter, which makes accurate approximations computationally expensive.
The same effects carry over to many Dynamic Complementarity Systems (DCS) \cite{Nurkanovic2020}.

On the theoretical side, necessary and/or sufficient conditions for optimality are provided in \cite{Colombo2020,Guo2016,Shaikh2007,Vieira2019}. 
Guo and Ye \cite{Guo2016} and Vieira et. al \cite{Vieira2019} study optimal control of a DCS with absolutely continuous solutions. 
The problems from this paper fall into this class. 
A broader overview can be found in \cite{Brogliato2020}.

On the practical side, many authors have developed methods to numerically treat discontinuous ODE in optimal control \cite{Baumrucker2009,Bemporad1999b,Bock2018,Katayama2020,Kirches2006,Nurkanovic2020,Pytlak2021,Stewart2010}. 
Kirches \cite{Kirches2006} develops a direct multiple shooting-based approach with a switch detecting integrator and sensitivity update formulae. Similarly, in \cite{Pytlak2021} a method is developed that can treat sliding modes on a single switching surface. 
Katayama et al. \cite{Katayama2020} fix the switching sequence and optimize the lengths of the phases in a model predictive control loop, where at every sample the switching sequence is updated. 
Bemporad et al. \cite{Bemporad1999b} assign integer variables to every mode and solve mixed-integer optimization problems. 
In \cite{Baumrucker2009,Bock2018,Nurkanovic2020} the ODE is transformed into a DCS, resulting in a Mathematical Program with Complementarity Constraints (MPCC) to be solved after discretization. 
In \cite{Nurkanovic2020} the authors use a standard discretization approach and suggest to use a homotopy to avoid spurious local minima due to wrong sensitivities \cite{Stewart2010}. 
Baumrucker and Biegler \cite{Baumrucker2009} consider systems with a single switching surface (or with multiple independent switching surfaces, cf. Section \ref{sec:cartesian_filippov}) and allow variable step sizes. 
This method yields exact switch detection, higher order integration accuracy, and correct numerical sensitivities. 
The step sizes are left to the optimizer as a degree of freedom, hence it can play with the discretization accuracy, possibly in an undesired way. Unfortunately, a formal proof of the appealing properties of the method is not provided in~\cite{Baumrucker2009}.

As it can be seen from the discussion above, most of the practical methods use only first-order accuracy methods with possibly incorrect sensitivities or treat the discontinuity in the integrator which complicates the use of derivative-based optimization algorithms (except \cite{Baumrucker2009}).
They do not treat sliding modes appropriately or handle only systems with a single switching surface, i.e., only two regions. 
The goal of this paper is to develop a method that resolves these issues and to provide a proper convergence theory. 
Note that the (easier) case of ODE with a continuous but nonsmooth r.h.s. fits into the structure of \eqref{eq:pws1}. 
Conversely, the \textit{time-freezing reformulation} \cite{Nurkanovic2023a,Nurkanovic2022a,Nurkanovic2021} transforms many systems from the (more difficult) case of systems with state jumps into the form of \eqref{eq:pws1}. 
This enables us to treat many classes of nonsmooth systems with the same numerical method in a unified way. FESD uses the DCS representation introduced by Stewart \cite{Stewart1990b,Stewart1996a} and is motivated by the variable step size ideas of Baumrucker and Biegler \cite{Baumrucker2009}.

\paragraph{Contributions:}
In this paper, we develop the FESD method, which can be used both in simulation and optimal control problems. We start with a reformulation of \eqref{eq:pws1} into dynamic complementary systems introduced by Stewart \cite{Stewart1990b,Stewart1996a} and provide a constructive way to pass from more natural definitions of discontinuous ODE to Stewart's form. 
Discretization of the DCS results in a nonlinear complementary problem. 
We build on the ideas of varying the step size and allowing switches only to take place at the boundaries of the finite elements introduced in \cite{Baumrucker2009}. The FESD method can efficiently deal with multiple and simultaneous switches including sliding modes on higher co-dimension surfaces, and thus is more general than \cite{Baumrucker2009}. 
Moreover, in contrast to \cite{Baumrucker2009}, where only Radau-IIA Implicit Runge-Kutta (IRK) methods were considered, in FESD one can use any Runge-Kutta method.

Additionally, we prove that FESD detects the switches \textit{exactly} in time, recovers the high-order accuracy that RK methods enjoy for smooth ODE, and obtains the correct numerical sensitivities even when the solution crosses or stays on a discontinuity. 
To allow switching on the boundaries of the finite elements we introduce the \textit{cross complementarity} formulation. 
Using FESD to discretize OCP with DCS results in mathematical programs with complementarity constraints which can be solved efficiently with smooth optimization techniques by solving few several Nonlinear Programming (NLP) problems. \cite{Anitescu2007,Hall2022,Kirches2022,Leyffer2006,Ralph2004}. Thus, we avoid the nonsmooth difficulty encountered within direct multiple shooting with switch-detecting integrators. 

Since the step sizes $h_n$ are allowed to vary, if no switches occur, we would encounter spurious degrees of freedom which let the optimizer play with the integrator accuracy in a possibly undesired way. 
To avoid this problem we propose additional conditions called \textit{step equilibration} which are based on an indicator function.
This decouples the integrator accuracy from the optimizer and results in piecewise equidistant discretization grids between the switches. 
We illustrate the practicability of FESD on several simulation and optimal control examples and verify the theoretical findings.

The FESD method with its many variations and the remainder of the tool-chain for numerically solving optimal control problems with nonsmooth systems are implemented in the open source MATLAB toolbox \texttt{NOSNOC} \cite{Nurkanovic2022c,Nurkanovic2022b}. 

\paragraph{Organization of the paper:} Section \ref{sec:piece_wise_smooth_systems} gives an introduction to Filippov systems and Stewart's reformulation. 
It provides a practical procedure to construct Stewart's indicator functions and discusses issues with discontinuous sensitivities. 
Section \ref{sec:FESD} introduces the FESD method and discusses the main ideas that lead to it.  
In Section \ref{sec:theory} several relevant theoretical properties of the FESD are studied. 
The sections contain numerical examples that illustrate the theoretical and algorithmic developments. 
Section \ref{sec:optimal_control} shows how to use FESD in numerical optimal control. 
The section finishes with an optimal control example and a benchmark compassion of FESD to the standard approach. 
Finally,  Section \ref{sec:summary} summarizes the paper and outlines future research.

\paragraph{Notation:} The complementary conditions for two vectors  $a,b \in \R^{n}$ \sloppy{read as ${0\leq a \perp b\geq 0}$, where $a \perp b$ means $a^{\top}b =0$}.  

For two scalar variables $a,b$ the so-called C-functions \cite[Section 1.5.1]{Facchinei2003} have the property $\phi(a,b) = 0 \iff a\geq 0, b\geq 0, ab = 0$. Examples are the natural residual functions $\phi_{\mathrm{NR}}(a,b)=\min(a,b)$ or the Fischer-Burmeister function $\phi_{\mathrm{FB}}(a,b) = a+b-\sqrt{a^2+b^2}$. 
If $a,b \in \R^{n}$, we use $\phi(\cdot)$ component-wise and define $\Phi(a,b) = (\phi(a_1,b_1),\dots,\phi(a_{n},b_{n}))$. 
\color{black}
All vector inequalities are to be understood element-wise, $\mathrm{diag}(x)\in \R^{n\times n}$ returns a diagonal matrix with $x \in \R^n$ containing the diagonal entries. 
The concatenation of two column vectors $a\in \R^{n_a}$, $b\in \R^{n_b}$ is denoted by $(a,b)\coloneqq[a^\top,b^\top]^\top$, the concatenation of several column vectors is defined analogously. 
The identity matrix is denoted by $I \in \R^{n \times n}$ and a column vector with all ones is denoted by $e=(1,1,\dots,1) \in \R^n$, their dimension is clear from the context.
The closure of a set $C$ is denoted by $\overline{C}$, its boundary as $ \partial C$ and ${\textmd{conv}}(C)$ is its convex hull. 
Given a matrix $M \in \R^{n \times m}$, its $i$-th row is denoted by $M_{i,\bullet}$ and its $j$-th column is denoted by $M_{\bullet,j}$. For a function $f:\R^{n} \to \R^{m}$ we denote by $\mathrm{D}f(x) = \partialder{f}{x}(x)\in \R^{m\times n}$ the Jacobian matrix and by
$\nabla f(x) \coloneqq \partialder{f}{x}(x)^\top $ its transpose. 

For the left and the right limits, we use the notation ${x(\ts^+)  = \lim\limits_{t\to \ts,\ t>\ts} x(t)}$ and ${x(\ts^-)  = \lim\limits_{t\to \ts,\  t<\ts}x(t)}$, respectively.
\color{black}
When clear from context, we often drop the dependency on time~$t$. 

\section{Piecewise smooth differential equations}\label{sec:piece_wise_smooth_systems}
This section will introduce some necessary assumptions on the systems PSS \eqref{eq:pws1}, its Filippov convexification \cite{Filippov1964}, and Stewart's reformulation into Dynamic Complementarity Systems (DCS) \cite{Stewart1990b,Stewart1996a} to prepare the ground for the novel method presented in Section \ref{sec:FESD}. 
We discuss some properties of the DCS for a fixed active set and active-set changes.
\subsection{Filippov convexification}
Initial value problems arising from the nonsmooth ODE \eqref{eq:pws1} usually fail to have classic Carath\'eodory solutions, for a counterexample, see e.g., \cite[Section 1]{Stewart1990b}.
To have a meaningful solution concept for this class of ODE, the main idea of Filippov was to replace the r.h.s. of \eqref{eq:pws1} with a convex set and to obtain the following Differential Inclusion (DI):
	\begin{align} \label{eq:FilippovDI}
		\dot{x}(t) &\in F_{\mathrm{F}}(x(t),u(t)) \coloneqq  		
		\bigcap_{\delta>0} \bigcap_{\mu(N) = 0}\overline{\mathrm{conv}} f(x+\delta B(x)\setminus N,u(t)),
	\end{align}
where $f(x,u) = f_i(x,u)$ if $x \in R_i$,  $B(x)$ is the Euclidean unit ball at $x$ in $\R^{n_x}$, $\mu(\cdot)$ is the Lebesgue measure on $\R^{n_x}$ and $\overline{\mathrm{conv}}(\cdot)$ maps a subset of $\mathbb{R}^{n_x}$ to its closed convex hull. 
Throughout the paper we assume that the regions $R_i$ are disjoint, connected and open. 
They are assumed to be nonempty and to have piecewise-smooth boundaries $\partial R_i$. 
We assume that $\overline{\bigcup\limits_{i\in \mathcal{J}} R_i} = \R^{n_x}$ and that $\R^{n_x} \setminus \bigcup\limits_{i\in\mathcal{J}} R_i$ is a set of measure zero.

Let $\mathcal{I}(x) \coloneqq \{ i \mid x \in \overline{R}_i\} \subseteq \mathcal{J}$ be the active set at $x \in \R^{n_x}$.
Due to the special structure of \eqref{eq:pws1} the Filippov  DI \eqref{eq:FilippovDI} can be written as
\begin{align*}
 \dot{x} \in \overline{\mathrm{conv}}\{f_i(x,u) \mid i \in \mathcal{I}(x)\}.
\end{align*}
This means that in the interior of the regions, $R_i$ the Filippov set $F_\mathrm{F}(x,u)$ is equal to $\{f_i(x,u)\}$ and on the boundary between regions we have a convex combination of the neighboring vector fields. 
If $\dot{x}$ exists, functions $\theta_i(\cdot)$ which serve as convex multipliers can be introduced and the Filippov DI can be written as 
\begin{align}\label{eq:FilippovDI_with_multiplers}
		\dot{x}  \in    F_{\mathrm{F}}(x,u) = \Big\{ \sum_{i\in \mathcal{J}}
		f_i(x,u) \, \theta_i  \mid \sum_{i\in \mathcal{J}}\theta_i = 1,
		 \ \theta_i \geq 0,
		 \ \theta_i = 0 \  \mathrm{if} \;  x \notin \overline{R_i}, 
		 \forall  i  \in \mathcal{J} \Big\}.
\end{align}
We call the functions $\theta_i(\cdot)$ \textit{Filippov multipliers}.
As it will be seen later the functions $\theta_i(\cdot)$ lack any continuity properties. 
But it can be shown that they are at least measurable \cite{Filippov1962,Stewart1990b}. 
Given \eqref{eq:FilippovDI_with_multiplers}, we will compute \textit{piecewise active} solutions \cite{Stewart1990b}, which are defined as follows.

\begin{definition1}[Piecewise active solution \cite{Stewart1990b}]\label{def:piecewise_active} 
For an initial value $x(0) = x_0$, a given measurable control function $u(t)$ and a compact interval $[0,T]$, a function $x:[0,T] \to \R^{n_x}$ is said to be a solution of \eqref{eq:FilippovDI}, if $\dot{x}(t) \in  F_{\mathrm{F}}(x(t),u(t))$ almost everywhere on $[0,T]$.
This function is called a \textit{piecewise active solution} if the active set $\mathcal{I}(x(t))$ is a piecewise constant function of time and it changes its value only finitely many times on $[0,T]$.  
A time point $\ts \in [0,T]$ is called a \textit{switching point} if $\mathcal{I}(x(t))$ is not constant in any sufficiently small neighborhood of $\ts$.
\end{definition1}
{
Note that this definition assumes a finite number of switches and excludes so-called Zeno solutions, where infinitely many switches can occur in a finite time interval.
Zeno solutions cannot be treated with event-detecting methods.
}

For a constant active set $\mathcal{I}(x)$ one can derive an ODE or Differential Algebraic Equation (DAE) (for sliding modes) from \eqref{eq:FilippovDI_with_multiplers} and apply standard integration methods. 
The overall ODE solution $x(t)$ is continuous and consists of smooth pieces connected by nondifferentiable points ("kinks") at the switching times $\ts$.
\subsection{Stewart's reformulation}\label{sec:stewarts_dcs}
We regard a specific representation of the sets $R_i$ which was introduced by Stewart~\cite{Stewart1990b}. 
The main assumption is that the regions $R_i$ are given as 
\begin{align}\label{eq:stewart_sets}
R_i  = \{ x \in \R^{n_x} \mid  g_i(x) < \min_{j\in \mathcal{J}\!,\, j \neq i } g_j(x)\}.	
\end{align}
It is assumed that the \textit{indicator} functions $g_i(\cdot),\ i \in \mathcal{J}$, are smooth functions. 
Moreover, throughout the paper we assume additionally that $g_i(\cdot), f_i(\cdot)$ and $\nabla g_i(\cdot)$ are Lipschitz continuous.

Note that due to the definition of the sets $R_i$ in \eqref{eq:stewart_sets}, the active set can be defined as 
\begin{align} \label{eq:active_set_defintion}
	\I(x(t)) \coloneqq \Bigl\{ i \in \mathcal{J} \mid g_i(x(t)) = \min_{j\in \mathcal{J}} g_j(x(t)) \Bigr\}.
\end{align}
\color{black}

We define the vectors $\theta = (\theta_1,\ldots,\theta_{\Nsys}) \in \R^{\Nsys}$, $g(x) = (g_1(x),\ldots,g_{\Nsys}(x)) \in \R^{\Nsys}$ and the matrix $F(x) = \begin{bmatrix}f_1(x), \ldots, f_{\Nsys}(x)\end{bmatrix}\in\R^{n_x \times {\Nsys}}$. Using the specific representations \eqref{eq:stewart_sets}, from \cite{Stewart1990b} and equation \eqref{eq:FilippovDI_with_multiplers} one can deduce that the Filippov DI can be written as
\begin{align}\label{eq:FilippovDI_with_multiplers_compact}
	\dot{x}= F(x,u) \theta(x),
\end{align}
where the algebraic variables $\theta(x(t))$ are a solution of the parametric Linear Program (LP)
\begin{align}\label{eq:stewart_lp}
\mathrm{LP}(x): \quad \theta(x)  \in\arg\min_{\tilde{\theta} \in \R^{\Nsys}} \quad & g(x)^\top \, \tilde{\theta} 
\quad	\textrm{s.t.} \quad  e^\top \tilde{\theta} = 1,\ \tilde{\theta}\geq0.
\end{align}
Using the Karush–Kuhn–Tucker (KKT) conditions of LP$(x)$ and \eqref{eq:FilippovDI_with_multiplers_compact}, we obtain the dynamic complementarity system
\begin{subequations}\label{eq:dcs_1}
	\begin{align}
		\dot{x} & = F(x,u)\theta,\\
		0& = g(x) - \lambda - \mu e,  \label{eq:dcs1_lambda} \\
		1& = e^\top \theta ,\\
		0& \leq \theta \perp \lambda  \geq 0,  \label{eq:dcs1_cc}
	\end{align}
\end{subequations}
where the algebraic variables $\lambda  \in \R^{\Nsys}$ and $\mu  \in \R$ are the Lagrange multipliers of the parametric LP \eqref{eq:stewart_lp}. To have an even more compact representation we use a C-function $\Psi$ for the complementarity conditions and rewrite the KKT conditions of the LP \eqref{eq:stewart_lp} as the nonsmooth equation:
\begin{align}\label{eq:dcs_lp}
	G_{\LP}(x,\theta,\lambda,\mu) \coloneqq 
	\begin{bmatrix}
		g(x) - \lambda - \mu e\\
		1- e^\top \theta\\
		\Psi( \theta,\lambda)
	\end{bmatrix} = 0.
\end{align}
It provides $2 \Nsys+1$ conditions for the $2 \Nsys+1$  algebraic variables $\theta,\,\lambda$ and $\mu$.
The DCS reads in compact form as a nonsmooth differential algebraic equation:
\begin{subequations}\label{eq:dcs_2}
	\begin{align}
		\dot{x} & = F(x,u)\theta,\\
		0&= G_{\LP}(x,\theta,\lambda,\mu).
	\end{align}
\end{subequations}
\begin{example1}\label{ex:dcs_defintion}
We illustrate this formulation on the simple example of $\dot{x} \in 2 - \mathrm{sign}(x)$. This ODE is characterized by the regions $R_1 = \{x \in \R \mid  x<0 \}$ and $R_2 = \{x \in \R \mid  x>0 \}$, with $f_1(x) = 3$, $f_2(x) = 1$ and $F(x) = [3 \quad 1]$.
It can be verified that with the functions $g_1(x) = x $ and $ g_2(x) = -x$ we have a representation of the regions as in \eqref{eq:stewart_sets}. Moreover, we have the multipliers $\theta,\lambda \in \R^2$ and $\mu \in \R$. 
Thus, the corresponding DCS reads as:
\begin{subequations}\label{eq:dcs_example}
\begin{align}
\dot{x} &= \begin{bmatrix}
3 &\quad 1
\end{bmatrix} (\theta_1,\theta_2),\\
0 &= x - \lambda_1 - \mu ,\\
0 &= -x - \lambda_2 - \mu ,\\
1 & = \theta_1+\theta_2,\\
0 & = \theta \perp \lambda \geq 0.
\end{align}
\end{subequations}
\end{example1}
\subsubsection{Remark on how to treat switching functions}\label{sec:stewart_indicator_functions}
Definition \eqref{eq:stewart_sets} might not be the most intuitive way to represent the sets $R_i$. 
In many practical examples some smooth scalar functions $c_i(\cdot)$, called \textit{switching functions}, are given. 
Their zero-level sets define the boundaries of the regions $R_i$.
For example, $R'_1 = \{ x\in \R^{n_x} \mid c_1(x)>0,\ldots,c_m(x)>0\}$, $R'_2 = \{ x\in \R^{n_x} \mid c_1(x)>0,\ldots,c_{m-1}(x)>0,c_m(x)<0\}$ and so on. Let $c(x) = (c_1(x),\ldots,c_m(x)) \in \R^m$ and assume that $\nabla c(x) \in \R^{n\times m}$ has rank $m$. 
Thus, we can locally define up to $\Nsys = 2^m$ regions and encode them via a sign matrix  $S \in \R^{2^m \times m}$ defined as 
\begin{align*}
	S = \begin{bmatrix}
		1 & 1 & \dots &1 & 1\\
		1 & 1 & \dots & 1& -1\\
		\vdots &  \vdots & \dots &  \vdots\\
		-1 & -1 & \dots &-1 & -1\\
	\end{bmatrix}.
\end{align*}
Note that the matrix $S$ has no repeating rows.
Moreover, we assume that this matrix has no zero entries. 
 The sets $R'_i$ can be compactly represented using the rows $S_{i,\bullet}$ as 
\begin{align}\label{eq:standard_sets}
	R'_i &= \{ x\in \R^{n_x} \mid \mathrm{diag} (S_{i,\bullet}) c(x)>0\}.
\end{align}
The next proposition provides a constructive way to find the functions $g(\cdot)$ from the more intuitive representation of the regions via $c(\cdot)$.
\begin{proposition1}
	Let the function $g: \R^{n_x}  \to \R^{\Nsys}$ be defined as 
	\begin{align}\label{eq:indicator_func_formula}
		g(x) = - S c(x),
	\end{align} then for all $ x \in R'_i$ the following statements are true:
	\begin{enumerate}[(i)]
		\item $g_i(x) < g_j(x),\ \mbox{for} \ i \neq j$,
		\item the definitions \eqref{eq:stewart_sets} and \eqref{eq:standard_sets} define the same set, i.e., $R_i = R'_i$.
	\end{enumerate}
\end{proposition1}
\textit{Proof.} 
For (i), note that for $x \in R'_i$ all terms in the sum $g_i(x) = -S_{i,\bullet}c(x) = -\sum_k S_{i,k} c_k(x)$ are strictly positive.  
On the other hand, for any $g_j(x)  = -S_{j,\bullet}c(x) =- \sum_k S_{j,k} c_k(x),\ j \neq i$ and  $x \in R'_i$, due to \eqref{eq:standard_sets}, all terms in the sum where $S_{j,k} \neq S_{i,k}$ are strictly negative. 
Therefore  $S_{i,\bullet}c(x) > S_{j,\bullet}c(x)$, thus (i) holds.

For (ii), first regard the rows $S_{j,\bullet}$ that differ from $S_{i,\bullet}$ only in the $k$-th column. Then $g_i(x) - g_j(x) = -(S_{i,k}-S_{j,k}) c_k(x) < 0$. 
If $S_{i,k} = 1$, then $g_i(x) - g_j(x) = -2c_k(x) < 0$. Likewise, for $S_{i,k} = -1$, then $g_i(x) - g_j(x) = 2c_k(x) < 0$. 
Therefore, from \eqref{eq:stewart_sets} we recover the definition of \eqref{eq:standard_sets} by looking at the rows where  $S_{i,k}$ and $S_{j,k}$ differ by one element.
For all rows $j$ that differ from $S_{i,\bullet}$ by more than one column, by similar reasoning, we obtain inequalities that do not tighten \eqref{eq:standard_sets}, since  $g_i(x) - g_j(x)$ consists of a sum of the terms from the inequalities where only one component of $c(x)$ is left. Therefore, statement (ii) holds and this completes the proof. 
\qed
\begin{example1}
For our tutorial example $\dot{x} \in 2-\mathrm{sign}(x)$ and the corresponding DCS \eqref{eq:dcs_example} we have $c(x) = x$ and $S = \begin{bmatrix}
	-1 & 1
\end{bmatrix}^\top$ and we obtain $g(x) = -Sc(x) = (x,-x)$ as used in Example \ref{ex:dcs_defintion}.
\end{example1}

\subsubsection{Fixed active set}\label{sec:stewarts_dcs_active_set_fixed}
\begin{figure}[t]
	\centering
	{\pgfplotsset{compat=1.13}
\setlength{\fwidth}{6.3cm}
\setlength{\fheight}{3.5cm}
\definecolor{mycolor1}{rgb}{0.00000,0.44700,0.74100}%
\definecolor{mycolor2}{rgb}{0.85000,0.32500,0.09800}%
\begin{tikzpicture}

\begin{axis}[%
width=0.951\fwidth,
height=\fheight,
at={(0\fwidth,0\fheight)},
scale only axis,
xmin=-3,
xmax=2,
xtick={\empty},
ymin=-1.5,
ymax=2,
ytick={\empty},
axis background/.style={fill=white},
legend style={legend cell align=left, align=left, draw=white!15!black},
xlabel style={font={\scriptsize}},ylabel style={font=\scriptsize},  ylabel shift={-0cm},ticklabel style={font=\scriptsize}
]
\addplot [color=black]
  table[row sep=crcr]{%
-1.85	-1.85\\
2	2\\
};

\addplot [color=black]
  table[row sep=crcr]{%
0	-0\\
1.85	-1.85\\
};

\addplot [color=black]
  table[row sep=crcr]{%
-3	0.75\\
0	-0\\
};

\addplot [color=mycolor1, line width=2.0pt]
  table[row sep=crcr]{%
-3	2\\
-2.5	0.625\\
};

\addplot [color=mycolor1, line width=2.0pt]
  table[row sep=crcr]{%
-2.5	0.625\\
-1	-1\\
};

\addplot [color=mycolor1, line width=2.0pt]
  table[row sep=crcr]{%
-1	-1\\
0	0\\
};

\addplot [color=mycolor2, line width=2.0pt, draw=none, mark=x, mark options={solid, mycolor2}]
  table[row sep=crcr]{%
-2.75	1.3125\\
};

\addplot [color=mycolor2, line width=2.0pt, draw=none, mark=x, mark options={solid, mycolor2}]
  table[row sep=crcr]{%
-2.5	0.625\\
};

\addplot [color=mycolor2, line width=2.0pt, draw=none, mark=x, mark options={solid, mycolor2}]
  table[row sep=crcr]{%
-1	-1\\
};

\addplot [color=mycolor2, line width=2.0pt, draw=none, mark=x, mark options={solid, mycolor2}]
  table[row sep=crcr]{%
0	0\\
};

\node[right, align=left]
at (axis cs:-2,-0.5) {$R_2$};
\node[right, align=left]
at (axis cs:-1.5,-0.4) {$x(t)$};
\node[right, align=left]
at (axis cs:-0.1,1.2) {$R_1$};
\node[right, align=left]
at (axis cs:-0.1,-0.5) {$R_3$};
\node[right, align=left]
at (axis cs:1.5,-0.5) {$R_4$};
\node[right, align=left]
at (axis cs:-2.65,1.312) {$x(t_1)$};
\node[right, align=left]
at (axis cs:-2.4,0.725) {$x(t_{\mathrm{s},1})$};
\node[right, align=left]
at (axis cs:-0.85,-1) {$x(t_{\mathrm{s},2})$};
\node[right, align=left]
at (axis cs:0.15,0) {$x(t_{\mathrm{s},3})$};
\end{axis}

\begin{axis}[%
width=1.227\fwidth,
height=1.227\fheight,
at={(-0.16\fwidth,-0.135\fheight)},
scale only axis,
xmin=0,
xmax=1,
ymin=0,
ymax=1,
axis line style={draw=none},
ticks=none,
axis x line*=bottom,
axis y line*=left,
legend style={legend cell align=left, align=left, draw=white!15!black},
xlabel style={font={\scriptsize}},ylabel style={font=\scriptsize},  ylabel shift={-0cm},ticklabel style={font=\scriptsize}
]
\end{axis}
\end{tikzpicture}%
}
	\hspace{-1.22cm}
	\caption{Illustration of active sets at different points. It can be seen that $\I(x(t_1)) = \I_0 = \{1\}$. At $x(t_{\mathrm{s},1})$ the trajectory crosses the surface of discontinuity between $R_1$ and $R_2$, hence  $\I(x(t_{\mathrm{s},1})) = \I_1^0 = \{1,2\}$ and later $\I_1 = \{2\}$. The segment between $x(t_{\mathrm{s},2})$ and $x(t_{\mathrm{s},3})$ is a sliding mode and we have $\I_2^0 = \{2,3\}$ and  $\I_2 = \{2,3\}$. Finally we have at $x(t_{\mathrm{s},3})$ that $\I_3^0 = \{1,2,3,4\}$.  }
	\label{fig:active_set_illustration}			
\end{figure}
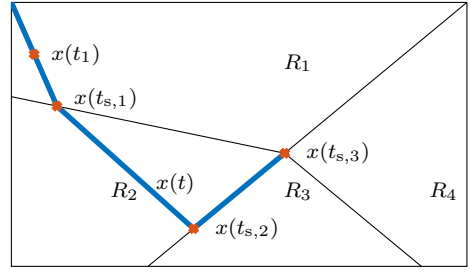

For a given solution $x(\cdot)$ let us denote all switching points by $0= t_{\mathrm{s},0}<t_{\mathrm{s},1}< \dots < t_{\mathrm{s},\Nswitch}=T$. 
The fixed active set between two switches is denoted by $\I_n\coloneqq \I(x(t))  , \ t\in (t_{\mathrm{s},n},t_{\mathrm{s},n+1})\eqqcolon I_n$ and at a switching point $t_{\mathrm{s},n}$ by $\I^0_n\coloneqq \I(x(t_{\mathrm{s},n}))$. 
Note that $\I_n^0 = \I_n \cup \I_{n-1}$. 
These definitions are illustrated in Figure \ref{fig:active_set_illustration}. 
In this subsection, we regard the DCS \eqref{eq:dcs_1} for a single fixed $\I_n$.
For ease of notation, we drop the subscripts in this subsection and denote the fixed active set by $\I$. 
Depending on the active set, the DCS \eqref{eq:dcs_1} reduces either to an ODE or to a DAE.

To simplify our exposition we introduce the following notation. 
For a given vector $a\in\R^n$ and set $\I \subseteq \{1,\ldots,n\}$, we define the projection matrix $P_\I \in \R^{|\I| \times n}$ which has zeros or ones as entries. 
It selects all component $a_i, i\in \I$ from the vector $a$, i.e., $a_\I = P_\I a \in \R^{|\I|}$ and $a_\I = [a_i \mid i\in \I]$.

In the DAE case, $x$ is on the boundary of one or more regions $R_i$, we speak of sliding modes \cite{Filippov1988}, i.e., $|\I|>1$ and typically obtain an index 2 differential algebraic equation.
In this case, two or more equal entries of $g(x)$ are the smallest components of this vector, and the solution $\theta$ of the LP$(x)$ is not unique and lies on a facet of the unit simplex. 
To compute the values of $\theta$, we must treat the DCS as a DAE. We define $F_{\I}(x,u) \coloneqq F(x,u) P_{\I}^\top$, which selects the appropriate columns of $F(x,u)$.
For $t \in I$, we have $\theta_i =0, i \notin \I$ and $\lambda_i = 0, i\in \I$, thus the DCS \eqref{eq:dcs_1} reduces to the DAE
\begin{subequations}\label{eq:dcs_dae}
	\begin{align}
		\dot{x} & = F_{\I}(x,u)\theta_{\I},\\
		0& = g_{\I}(x) - \mu e, \label{eq:dcs_dae_switching} \\
		1& = e^\top \theta_{\I}.
	\end{align}
\end{subequations}
There are $|\I|+1$ nontrivial algebraic equations and $|\I|+1$ unknown algebraic variables, namely $\mu$ and $\theta_i$ for $i \in \I$, since we consider $\theta_i(t)=0,\ i \notin \I$ as fixed.

In the ODE case, $x$ is in the interior of some region $R_i$, we have $|\I|=1$. 
The algebraic variables $\mu$ and $\theta_i$ can be computed explicitly from \eqref{eq:dcs_dae} and we have $\theta_i = 1$ and $\mu = g_i(x)$. 
Thus, the DCS reduces to the ODE $\dot{x} = f_i(x)$.

Next, we provide sufficient conditions for solution uniqueness of the DAE \eqref{eq:dcs_dae} for a given $|\I|\geq 1$. 
We define the matrix
\begin{align}\label{eq:stewart_matrix}
M_{\I}(x) &= \nabla g_\I(x)^\top F_\I(x,u) \in \R^{|\I| \times |\I|}.
\end{align}
Note that entries of this matrix arise by taking the total time derivative of~\eqref{eq:dcs_dae_switching}.
\begin{assumption}\label{ass:solution_existence}
	Given a fixed active set $\I(x(t)) = \I$ for $t \in I$,  it holds that the matrix $M_{\I}(x(t))$ \sloppy{is invertible and $e^\top M_{\I}(x(t))^{-1} e \neq 0$ for all $t\in I$}.
\end{assumption}
\begin{proposition1}\label{prop:solution_existence}
	Suppose that Assumption \ref{ass:solution_existence} holds. Given the initial value $x(t_{\mathrm{s},n})$, then the DAE \eqref{eq:dcs_dae} has a unique solution for all $t\in I$.
\end{proposition1}
\textit{Proof.}
For a given $x(\cdot)$ we can differentiate equation \eqref{eq:dcs_dae_switching} w.r.t. $t$ and obtain the following index 1 DAE
\begin{subequations}
\begin{align}
&	\dot{x}  = F_{\I}(x,u)\theta_{\I},\; \dot{\mu} = -v,\\
&\begin{bmatrix}
	M_{\I}(x) &e \\
	e^\top &0
\end{bmatrix}
\begin{bmatrix}
	\theta_{\I} \\
	v
\end{bmatrix}
=
\begin{bmatrix}
	0 \\
	1
\end{bmatrix} \label{eq:stewart_lin_system}.
\end{align}
\end{subequations}
with the algebraic variables $\theta_{\I}$ and $v\in \R$.  
For a given initial 
condition $x(t_{\mathrm{s,n}})$, $\mu(t_{\mathrm{s,n}})$ can be directly computed from any component of \eqref{eq:dcs_dae_switching}.
Using the Schur complement and Assumption \ref{ass:solution_existence}, we conclude that we can find unique $\theta_{\I}$ and $v$ by solving the linear system \eqref{eq:stewart_lin_system}. 
Therefore, the DAE \eqref{eq:dcs_dae} can be reduced to an ODE. 
Since the functions $f_i$ are assumed to be Lipschitz the resulting ODE has a unique solution $x(t), t\in I$. \qed

A similar result, with a more complicated proof but different assumptions can be found in \cite[Section 2]{Stewart1990b}.
\color{black}
Note that even though the DAE has a unique solution for a given active set $\I$, there might be multiple $\I$ that give a well-defined ODE, as we discuss in the subsequent sections.
We do not know a priori whether we need to treat an ODE or a DAE, but for both cases, we will use Runge-Kutta methods within FESD to provide high-accuracy solutions. 
The crucial part of FESD is the automatic active set and switching time detection so that sliding modes and crossings of region boundaries can be treated in a unified way.

\subsubsection{Active-set changes and continuity of $\lambda$ and $\mu$}\label{sec:stewarts_dcs_active_set_changes}
Every active-set change in \eqref{eq:dcs1_cc} corresponds to crossing a discontinuity, entering or leaving a sliding mode, or a spontaneous leaving of a surface of discontinuity. 
These events in time are called \textit{switches}.

From \eqref{eq:FilippovDI_with_multiplers}, Eq. \eqref{eq:dcs1_lambda} and the complementarity conditions \eqref{eq:dcs1_cc} for $i \in \I(x)$ it follows that $\theta_i \geq 0$ and $\lambda_i =0$. 
Likewise, for $i \notin \I(x)$ it follows that $\theta_i = 0$ and $\lambda_i \geq 0$. 
Hence,  for $i  \in \I(x)$ from \eqref{eq:dcs1_lambda} and \eqref{eq:active_set_defintion} we conclude that $\mu =  \min_{j\in \mathcal{J}} g_j(x)$.
\begin{lemma1}\label{lem:lambda_cont}
	The functions $\lambda(t)$ and $\mu(t)$ in \eqref{eq:dcs_1} are continuous in time.
\end{lemma1}
\textit{Proof}: The function $\mu(t)$ is a minimum of continuous functions and is thus continuous. 
Therefore, continuity of $\lambda(t) = g(x(t)) - \mu(t)e$ follows from the continuity of $x(t)$ and $g(x)$ and Equation \eqref{eq:dcs1_lambda}. \qed 

\begin{remark1}\label{rem:lambda_at_switch}
Continuity of $\lambda(t)$ implies that at an active-set change of a component $i$ at $t_{\mathrm{s},n+1}$ some $\lambda_i(t_{\mathrm{s},n+1})$ must be zero. 
Moreover, for some $i \notin \I_n$, in the case of crossing a discontinuity or entering a sliding mode with $i\in \I_{n+1}$, it holds that {the left time derivative of $\lambda_i$ is negative, i.e.,} $\dot{\lambda}_i(t_{\mathrm{s},n+1}^-) < 0$. 
Likewise, in the case of leaving a sliding mode or a spontaneous switch, with $i \in \I_n$ and $i \notin \I_{n+1}$, it follows that {the right time derivative of $\lambda_i$ is positive, i.e.,} $\dot{\lambda}_i(t_{\mathrm{s},n+1}^+) > 0$. 
If some of the first-order one-sided derivatives of ${\lambda}(\cdot)$ are zero at a switching point $t_{\mathrm{s},n+1}$, then one must look at higher-order derivatives to determine if it stays active or not.
\end{remark1}
We exploit the continuity of $\lambda(\cdot)$ and $\mu(\cdot)$ later in the derivation of the FESD method. 
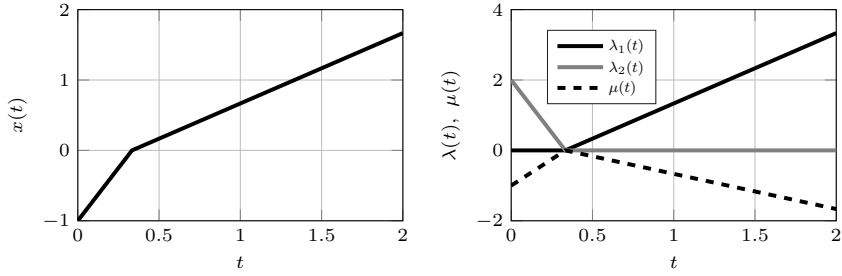
\begin{figure}[t]
	\centering
	{\pgfplotsset{compat=1.13}
\setlength{\fwidth}{9.5cm}
\setlength{\fheight}{2.8cm}
\begin{tikzpicture}

\begin{axis}[%
width=0.45\fwidth,
height=\fheight,
at={(-0.04\fwidth,0\fheight)},
scale only axis,
xmin=0,
xmax=2,
xtick= {0, 0.5,1,1.5,2},
 x tick label style={
	/pgf/number format/.cd,
	fixed,
	/tikz/.cd
},
xlabel style={font=\color{white!15!black}},
xlabel={$t$},
ymin=-1,
ymax=2,
ylabel style={font=\color{white!15!black}},
ylabel={$x(t)$},
axis background/.style={fill=white},
xmajorgrids,
ymajorgrids,
legend style={legend cell align=left, align=left, draw=white!15!black},
xlabel style={font={\scriptsize}},ylabel style={font=\scriptsize},  ylabel shift={-0cm},ticklabel style={font=\scriptsize}
]
\addplot [color=black, line width=1.5pt]
  table[row sep=crcr]{%
0	-1\\
0.333333333333333	0\\
2	1.66666666666711\\
};
\end{axis}

\begin{axis}[%
width=0.45\fwidth,
height=\fheight,
at={(0.56\fwidth,0\fheight)},
scale only axis,
xmin=0,
xmax=2,
 x tick label style={
	/pgf/number format/.cd,
	fixed,
	/tikz/.cd
},
xlabel style={font=\color{white!15!black}},
xlabel={$t$},
ymin=-2,
ymax=4,
ylabel style={font=\color{white!15!black}},
ylabel={$\lambda(t),\ \mu(t)$},
axis background/.style={fill=white},
xmajorgrids,
ymajorgrids,
legend style={at={(0.11,0.552)}, anchor=south west, legend cell align=left, align=left, draw=white!15!black,nodes={scale=0.70, transform shape}},
xlabel style={font={\scriptsize}},ylabel style={font=\scriptsize},  ylabel shift={-0cm},ticklabel style={font=\scriptsize}
]
\addplot [color=black, line width=1.5pt]
  table[row sep=crcr]{%
0	0\\
0.333333333333333	0\\
2	3.3333333333333\\
};
\addlegendentry{$\lambda_1(t)$}

\addplot [color=gray, line width=1.5pt]
  table[row sep=crcr]{%
0	2\\
0.333333333333333	0\\
2	0\\
};
\addlegendentry{$\lambda_2(t)$}

\addplot [color=black, dashed, line width=1.5pt]
  table[row sep=crcr]{%
0	-1\\
0.333333333333334	0\\
2	-1.66666666666711\\
};
\addlegendentry{$\mu(t)$}
\end{axis}

\begin{axis}[%
width=1.227\fwidth,
height=1.227\fheight,
at={(-0.16\fwidth,-0.135\fheight)},
scale only axis,
xmin=0,
xmax=1,
ymin=0,
ymax=1,
axis line style={draw=none},
ticks=none,
axis x line*=bottom,
axis y line*=left,
legend style={legend cell align=left, align=left, draw=white!15!black},
xlabel style={font={\scriptsize}},ylabel style={font=\scriptsize},  ylabel shift={-0cm},ticklabel style={font=\scriptsize}
]
\end{axis}
\end{tikzpicture}%
}
	\hspace{-1.22cm}
	\caption{{Illustration of the arguments of Lemma \ref{lem:lambda_cont} and Remark \ref{rem:lambda_at_switch} on the Example \ref{ex:switching_cases} (a) and its corresponding DCS \eqref{eq:dcs_example} for $t\in[0,2]$ and $x(0) = -1$ with $\ts = -\frac{1}{3}$. The functions $\mu = \min(-x,x)$, $\lambda_1 = x-\mu$ and $\lambda_2 =- x-\mu$ are continuous. 
			At the switching point $\ts$ we have $\dot{\lambda}_1(\ts^-)=0$, $\dot{\lambda}_1(\ts^-)>0$ and $\dot{\lambda}_2(\ts^-)<0$, $\dot{\lambda}_1(\ts^-)=0$.}}
	\label{fig:active_set}	
	\vspace{-0.45cm}					
\end{figure}
The next example discusses the difference between the possible switching cases.
\begin{example1}\label{ex:switching_cases}
There are four possible switching cases which we illustrate with the following examples:
\begin{enumerate}[(a)]
	\item crossing a surface of discontinuity, $\dot{x}(t) \in 2-\mathrm{sign}(x(t))$,
	\item sliding mode, $\dot{x}(t) \in -\mathrm{sign}(x(t))$,
	\item leaving sliding mode $\dot{x}(t) \in -\mathrm{sign}(x(t))+t$.
	\item spontaneous switch, $\dot{x}(t) \in \mathrm{sign}(x(t))$,
\end{enumerate}
In case (a), for $x(0) < 0$ the trajectory reaches $x = 0$ and crosses it. 
In example (b), for any finite $x(0)$, the trajectory reaches $x=0$ and stays there. 
On the other hand, in example (c), for $x(0) =0$ the DI has a unique solution and leaves $x=0$ at $t=1$. 
In the last example, for $x(0) =0$, the DI has infinitely many solutions, and $x(t)$ can spontaneously leave $x = 0$ at any $t\geq0$.
Note that there is a qualitative difference between leaving a sliding mode (c) and spontaneous switch (d).   
The arguments of Lemma \ref{lem:lambda_cont} and Remark \ref{rem:lambda_at_switch} are illustrated in Figure \ref{fig:active_set} for the Example \ref{ex:switching_cases} (a).
\end{example1}

\subsubsection{Predicting the new active set}\label{sec:stewarts_lcp}
In this subsection, we restate a more technical result from \cite{Stewart1990b} which is later needed in the convergence proof. 
The reader not interested in the proofs may skip this part.

As already noted in Remark \ref{rem:lambda_at_switch}, switches are characterized by the time derivative of $\lambda(\cdot)$. 
Note that, e.g., for crossing a discontinuity or entering a sliding mode at a switching point and for the subsequent interval it holds that $\I_n \subseteq \I_n^0$. 
Moreover, one can construct a Linear Complementarity Problem (LCP) with the data at $x(t_{\mathrm{s},n})$ and predict $\I_{n+1}$. 

We define the vector $w_{\I}(t) \coloneqq \frac{\dd}{\dd t}\lambda_{\I}(t) = M_{\I}(x(t)) \theta_{\I}(t) - \dot{\mu}(t) e$. 
One can construct the following mixed LCP between $\dot{\lambda}_{\I_n^0}$ and $\theta_{\I_n^0}$ at $t_{\mathrm{s},n}$:
\vspace{-0.1cm}
\begin{subequations}\label{eq:mixed_lcp_switching}
	\begin{align}
		w_{\I^0_n}&=  M_{\I_n^0}(x)\theta_{\I_n^0} - \dot{\mu}e, \\
		1 &=e^\top \theta_{\I_n^0},\\
		0 &\leq w_{\I_n^0} \perp \theta_{\I_n^0} \geq 0.
	\end{align}	
\end{subequations}
For a sufficiently large $\alpha>0$ all entries of the matrix  $M_{\I,\alpha}(x) = M_{\I}(x)+ \alpha e e^\top$ are strictly positive. 
This means the matrix $M_{\I,\alpha}(x)$ is strictly copositive, i.e., for any $a\geq0, a\neq 0$ it holds that $a^\top M_{\I,\alpha}(x) a >0$ \cite{Facchinei2003}.
One can derive an LCP equivalent to \eqref{eq:mixed_lcp_switching} \cite[Lemma 3.3]{Stewart1990b}:
\begin{align}\label{eq:lcp_switching}
	0 & \leq \tilde{w}_{\I_n^0} = 	M_{\I_n^0,\alpha}(x)\tilde{\theta}_{\I_n^0}  -e \perp \tilde{\theta}_{\I_n^0} \geq 0.
\end{align}
The motivation for rewriting \eqref{eq:mixed_lcp_switching} as \eqref{eq:lcp_switching} is twofold. 
It is both easier to prove solution existence and to compute a solution for an LCP with a strictly copositve matrix than for the initial mixed LCP \cite{Stewart1990b}.
The solution of the initial LCP $(w_{\I_n^0},\theta_{\I_n^0})$ can be reconstructed via $\theta_{\I_n^0} = {\tilde{\theta}_{\I_n^0}}/{e^\top \tilde{\theta}_{\I_n^0}}$ and $w_{\I_n^0} = {\tilde{w}_{\I_n^0}}/{e^\top \tilde{\theta}_{\I_n^0}}$, for further details cf. \cite[Lemma 3.3]{Stewart1990b}. 
There is a one-to-one correspondence between the active set in a neighborhood of a switching point $t_{\mathrm{s},n}$ and the solutions of the LCP \eqref{eq:lcp_switching}. 
This is summarized in the next theorem proved by Stewart \cite{Stewart1990b}.
\begin{theorem1}[Theorem 3.2 \cite{Stewart1990b}] \label{th:index_sets}
	Let $x(t)$ be a solution in the sense of Definition \ref{def:piecewise_active} for $t \in [t_a,t_b]$, with $\I^0 = \I(x(t_a))$ and $\I = \I(x(t))$ for all $t\in (t_a,t_b)$. 
	Suppose Assumption \ref{ass:solution_existence} holds for all  $t\in (t_a,t_b)$. Then for each $t \in (t_a,t_b)$ there is a solution of the LCP \eqref{eq:lcp_switching} such that
	\begin{align*}
		\{i \mid \tilde{\theta}_i>0\} &\subseteq \I \subseteq \{i \mid \tilde{w}_i = 0 \}.
	\end{align*}
	Conversely, let $x_0\in \R^{n_x}$ and $t_a$ be given with $\I^0 = \I(x_0)$. 
	Then if $(\tilde{w}_{\I^0},\tilde{\theta}_{\I^0})$ is a solution of the LCP \eqref{eq:lcp_switching} such that
	\begin{align*}
		\{i \mid \tilde{\theta}_i>0\} &=\I = \{i \mid \tilde{w}_i = 0 \}
	\end{align*}
	and the conditions of Assumptions \ref{ass:solution_existence} are satisfied for $\nabla g_i(x)$ and $f_i(x,u)$, $i \in \I$,
	 then there is a $t_b > t_a$ and a solution $x(\cdot)$ in the sense of Definition \ref{def:piecewise_active} on $[t_a,t_b]$ such that $x(t_a) = x_0$ and $\I(x(t)) = \I$ for all $t\in(t_a,t_b)$. 
\end{theorem1}
Regard an LCP 
\begin{align}\label{eq:generic_lcp}
	0&\leq M \theta + q \perp \theta\geq 0,
\end{align}
with $M \in \R^{l \times l}$ and $q\in\R^{l}$. 
The given LCP \eqref{eq:generic_lcp} is compactly denoted by $\mathrm{LCP}(M,q)$ and its set of solutions is denoted by $\mathrm{SOL}(M,q)\subseteq \R^l$. If a solution satisfies $(M \theta + q) +\theta >0$, we say that strict complementarity holds.

To show convergence we will require the solutions of the LCP to be \textit{strongly stable} \cite{Facchinei2003,Stewart1990b}. 
A solution $(w^*,\theta^*)\in \mathrm{SOL}(M,q)$ of a given $\mathrm{LCP}(M,q)$ is said to be strongly stable if there is a neighborhood $U$ of $\theta^*$ and a neighborhood $V$ of the problem data $M\in \R^{l\times l}$ and $q\in \R^{l}$, such that the intersection of $U$ with the solution set of an LCP constructed from the data from one point in $V$ is a singleton. 
We state a regularity assumption about the LCP \eqref{eq:lcp_switching}.
\begin{assumption}\label{ass:lcp_switching}
Consider a solution $x(t)$ in the sense of Definition \ref{def:piecewise_active} for $t\in[0, T]$, {and let $\mathcal{S}=\{ t_{\mathrm{s},0},\ldots, t_{\mathrm{s},\Nswitch}\}$ be the set of switching points.
The solutions of the LCP \eqref{eq:lcp_switching} are strongly stable and satisfy strict complementarity for all  
$t\in [\ts-\epsilon, \ts+\epsilon]\cap[0,T]$, $\ts \in \mathcal{S}$, for a sufficiently small $\epsilon>0$.}
\end{assumption}

The strict complementarity assumption is needed to obtain a tight prediction of the next active set $\I$, cf. first part of Theorem \ref{th:index_sets}. 
From the proof of~\cite[Theorem 3.2]{Stewart1990b} it follows that the strict complementarity condition implies that the one-sided time derivatives of $\lambda_i(t), i \notin \I(x(t))$ are nonzero, see also Remark \ref{rem:lambda_at_switch}.
Without this assumption, one can obtain only an over-approximation of $\I$.
However, it can be relaxed at the cost of looking at higher-order time derivatives of $\lambda_i(t_{\mathrm{s},n}), i \in \I_{n}^0$ and constructing an appropriate LCP for determining the active sets past some switching point, cf.~\cite[Section 4.2]{Stewart1990c} for derivations.
Note that the strict complementarity is needed only in a neighborhood of the switching points. 
Strong stability is assumed in order to apply some results on parametric LCPs. 
In our case, we will use it to draw the same conclusions from LCPs constructed at $t$ and $t'$, where $t$ and $t'$ are sufficiently close.
\color{black}
\begin{example1}
We briefly illustrate Theorem \ref{th:index_sets} on our example $\dot{x}  = 2-\mathrm{sign}(x)$ with $x(0) = -1$, cf. Figure \ref{fig:active_set}. It is easy to see that $t_{\mathrm{s},1} = -\frac{1}{3}$ and that the relevant active sets are $\I_0 = \{1\}$, $\I_1 =\{2\}$ and $\I^0_1(t_{\mathrm{s},1}) = \{1,2\}$. 
The LCP \eqref{eq:lcp_switching} for our example at $t_{\mathrm{s},1}$ reads as 
\begin{align*}
	0 & \leq \tilde{w}_{\I^0_1} = 	\begin{bmatrix}
		3+\alpha &&1+\alpha\\
		-3+\alpha && -1+\alpha
	\end{bmatrix}\tilde{\theta}_{\I_1^0}  -e \perp \tilde{\theta}_{\I^0_1} \geq 0.
\end{align*}
With $\alpha = 5$ this LCP has the unique solution $\tilde{\theta}_{\I_1^0} = (0,\frac{1}{4})$ and $\tilde{w}_{\I_1^0} = (\frac{1}{2},0)$ and according to the last theorem it correctly predicts $\I_1 = \{2\}$.
\end{example1}

\subsection{Remark on Cartesian products of Filippov systems} \label{sec:cartesian_filippov}
The reformulation from the last subsection given by the DCS \eqref{eq:dcs_1} fails on some simple examples such as: $\dot{x}_1 \in -\sign(x_1), \ \dot{x}_2\in -\sign(x_2)$, $x \in \R^2$. 
This example satisfies the one-sided Lipschitz condition and has a unique Filippov solution \cite{Filippov1988,Stewart1996a}.
However, as shown in \cite{Stewart1996a} at $(0,0)$ the DAE arising from \eqref{eq:dcs_1} fails to have a unique solution. 
One can see that $\dot{x}_1 \in -\sign(x_1)$ and $\dot{x}_2 \in -\sign(x_2)$ are completely independent and thus they should be treated in such a way. 

Stewart introduced a generalization of his reformulation for such cases in \cite{Stewart1996a}. 
One should identify the $\Nsubsys$ {independent subsystems} with index $k = 1,\ldots,\Nsubsys$, where each subsystem has $n_f^k$ modes.
We equip all variables related to the $k-$th subsystem with the superscript $k$.
Instead of \eqref{eq:FilippovDI_with_multiplers} one can write
\begin{align}\label{eq:filippov_di_with_multiplers_multiple}
	\dot{x}& \in \Big\{ \sum_{k =1}^{\Nsubsys} \! \sum_{i= 1}^{n_f^k} \theta_{i}^k f_{i}^k (x,u) \ \mid\ \sum_{i= 1}^{n_f^k} \theta_{i}^k =1,\;
	\theta^k \geq 0, k = 1,\dots, \Nsubsys \Big\}.
\end{align}
Finding the functions ${g}^{k}(\cdot) \in \R^{n_f^k}$ from ${c}^{k}(\cdot) \in \R^{{n_f^k}}$ for every {subsystem} works the same way as in Section \ref{sec:stewart_indicator_functions}. 
Thereby, the regions of every subsystem are defined via the matrix $S^k$ and the switching functions $c^k(x) \in \R^{n_c^k}$.
Every mode's convex combination is encoded by its parametric linear program \eqref{eq:stewart_lp}, constructed with the $k$-th modes' switching functions $g^k(x) \in \R^{n_f^k}$.
Thus, we can derive the DCS
\begin{subequations}\label{eq:dcs_multiple}
	\begin{align}
		\dot{x}  & = \sum_{k=1}^{\Nsubsys} F^k(x,u ){\theta}^k ,\\
		0& = {g}^k(x ) - {\lambda}^k  - {\mu}^k e, &\mbox{for all} \ k \in \{1,\dots \Nsubsys\},\\
		1& = e^\top {\theta}^k, &\mbox{for all} \ k \in \{1,\dots \Nsubsys\},\\
		0& \leq {\theta}^k  \perp {\lambda}^k  \geq 0, &\mbox{for all} \ k \in \{1,\dots \Nsubsys\},
	\end{align}
\end{subequations}
where $F^k(x,u) = [f_{1}^{k}(x,u),\ldots,f_{n_f^k}^{k}(x,u)] \in \R^{n_x \times n_f^k}$ and  ${g}^{k}(x) \in \R^{n_f^k}$, ${\theta}^k\in \R^{n_f^k}$, ${\lambda}^k  \in \R^{n_f^k}$ and ${\mu}^k \in \R$, for all $k \in \{ 1,\dots, \Nsubsys\}$. 
For ease of notation, in the remainder of the paper we treat the case with $\Nsubsys = 1$, as all extensions are straightforward. 

To the best of the authors' knowledge, there are no general conditions known which identify when the r.h.s. of \eqref{eq:pws1} is {partially separable} as in \eqref{eq:filippov_di_with_multiplers_multiple} and there might even be multiple ways to write it in this form. 
However, in practice, it is usually easy to identify the structure of \eqref{eq:filippov_di_with_multiplers_multiple} by inspection. 
For example, this occurs if we have multiple surfaces with friction, or multiple objects touching the same frictional surface \cite{Stewart1996a}.

\subsection{Sensitivities with respect to parameters and initial values} \label{sec:sensitivites}
Correct calculation of derivatives of solutions w.r.t. parameters (e.g., discretized control functions) and initial values is crucial for efficient numerical optimal control algorithms and verifying the optimality of a solution. 
This is not straightforward for ODE with a discontinuous r.h.s., as the sensitivity usually exhibits jumps when switches occur. 
As any constant parameter $\hat{p}$ can be modeled via adding the state $\dot{p} = 0$ and $p(0) = \hat{p}$, we restrict our attention to sensitivities w.r.t. initial values.

Regard the DCS given by Eq. \eqref{eq:dcs_1} on a time interval $[0,T]$ with the initial condition $x(0) = x_0$. 
Assume that the surface $\partial R_j$ is reached at $\ts(x_0)\in(0,T)$ and that $x_0 \in R_i$. 
We consider the case where the solution crosses a co-dimension one surface of discontinuity $\partial R_j$. 
Other cases are where the trajectory:
(a) slides on the surface of discontinuity after reaching it,
(b) starts on a surface of discontinuity and stays on it or leaves it, or
(c) goes from one surface to another.
They can be analyzed with the same arguments as below, but we omit these cases here for brevity, cf. \cite[Section 2.11]{Filippov1988}.

In the case of crossing, we have for $t \in \left[0,\ts\right)$ that $\I(x(t)) = \{i\}$ and from \eqref{eq:dcs_1} it follows that $\dot{x} = f_i(x)$. 
After crossing $\partial R_j$ at $\ts$ we have $\I(x(t)) = \{ j\}$ for $t\in \left(\ts,T\right]$ and $\dot{x} = f_j(x)$.
At $\ts$ it holds that $\psi_{i,j}(x(\ts))=0$ with
\begin{align}\label{eq:switching_function_psi}
\psi_{i,j}(x(t))&\coloneqq g_i(x(t)) -  g_j(x(t)).
\end{align}
  Thus, the system can be compactly represented by
\begin{align}\label{eq:compact_sens}
	\dot{x}(t)  &=  \begin{cases}
		f_i(x(t)),\ &\psi_{i,j}(x(t)) <0, \\   	
		f_j(x(t)),\ &\psi_{i,j}(x(t)) \geq 0.
	\end{cases}
\end{align}
We are interested in the exact sensitivity matrix $X(t,0;x_0) = \frac{\partial x(t;x_0)}{\partial x_0} \in \R^{n_x \times n_x}$ of a solution $x(t;x_0)$ of the system \eqref{eq:compact_sens}. 
The function $X(t,0;x_0)$ obeys smooth linear variational differential equations on both sides of $\ts$, but exhibits a jump at $\ts$~\cite{Filippov1988}.
{The statement of the next proposition is adapted from \cite[Section 3.3]{Stewart2010}.}
\begin{proposition1}\label{prop:sensitivity_exact}
Regard the system \eqref{eq:compact_sens} with $x(0) = x_0 \in R_i$ on an interval $[0,T]$ with a switch at $\ts \in (0,T)$. 
Assume that the functions $f_i(x),\; f_j(x),\psi_{i,j}(x)$ are continuously differentiable along $x(t), t\in [0,T]$. 
Assume the solution $x(t)$ reaches the surface of discontinuity transversally, i.e., $\nabla\psi_{i,j} (x(\ts))^\top f_i (x(\ts))>0$. 
Then the sensitivity $X(T,0;x_0)$ of a solution $x(t;x_0)$ of the system described by the ODE \eqref{eq:compact_sens} is given by
\begin{align}\label{eq:exact_parametric_sensitivites}
\begin{split}
		X(T,0;x_0) &= X(T,\ts^+;x(\ts)) J(x(\ts;x_0)) X(\ts^-,0;x_0)\ \text{with}\\
		J(x(\ts;x_0)) &\coloneqq  I + \frac{(f_j(x(\ts;x_0))-f_i(x(\ts;x_0)))\nabla\psi_{i,j} (x(\ts;x_0))^\top}{\nabla\psi_{i,j} (x(\ts;x_0))^\top f_i (x(\ts;x_0))} .
\end{split}
\end{align}
\end{proposition1}
This proposition can also be adapted to the case of sliding modes. We obtain similar expressions for the sensitivity jump formula as in \eqref{eq:exact_parametric_sensitivites}. 
The only change needed to be made is to replace $f_j(x)$ with $f^*(x)$, where $f^*(x)$ defines the sliding vector field \cite{Filippov1964}.

Since numerical sensitivities obtained via {standard time-stepping methods} fail to converge to their correct values \eqref{eq:exact_parametric_sensitivites} \cite{Nurkanovic2021,Stewart2010}, artificial local minima arbitrarily close to the initialization point may arise in the context of optimization and impair the progress of the optimizer. 
This is resolved within FESD, where the convergence of the discrete-time sensitivities is recovered, cf. Section \ref{sec:sensitivity_convergence}.

\section{Finite Elements with Switch Detection} \label{sec:FESD}
This section introduces the main algorithmic ingredients of the FESD method. 
The goal of the method is to:
(a) detect exactly the time of reaching or leaving the region boundaries which is necessary for high accuracy of integration methods,
(b) exactly compute the sensitivities across regions in order to correctly treat the nonsmoothness and 
(c) appropriately treat the possible evolution on the boundary that is present in sliding modes.
 
In this section, we regard a single control interval $[0,T]$ with a constant externally chosen control input $q \in \R^{n_u}$, i.e., we set $u(t) = q$ for $ t\in [0,T]$. Extensions with more complex smooth parametrizations of the control function are straightforward.
\subsection{Standard Runge-Kutta discretization}\label{sec:dcs_irk}
As a starting point in our analysis, we regard a standard Runge-Kutta (RK) discretization of the DCS~\eqref{eq:dcs_1}. 
In the nonsmooth ODE community, these schemes are known as \textit{time-stepping} methods. 
Opposed to \textit{event-based/switch-detection} methods, they assume fixed step sizes $h_n$ and do not try to detect the switches. 
As a consequence, they have in general only first-order accuracy \cite{Acary2008}. 
The theoretical properties of RK methods for DI and DCS have been studied by many authors, e.g., \cite{Dontchev1992,Kastner1990,Stewart2011,Taubert1981}. 

Suppose the initial value $x(0) = s_0$ is given. We divide the control interval into $\NFE$ \textit{finite elements} (i.e., integration intervals) $[t_{n},t_{n+1}]$ via the grid points $0= t_0 < t_1 < \ldots <t_{\NFE} = T$. 
On each of the finite elements we consider an $\Nstg$-stage Runge-Kutta method which is characterized by the Butcher tableau entries $a_{i,j} ,b_i$ and $c_i$ with $i,j\in\{1,\ldots,\Nstg\}$ \cite{Hairer1991}.
The fixed step size reads as $h_{n} = t_{n+1} - t_{n},\; n = 0, \ldots,\NFE-1$. 
The approximation of the differential state at the grid points $t_n$ is denoted by $x_n \approx x(t_n)$. 
We regard a \textit{differential} representation of the Runge-Kutta method where the derivatives of states at the stage points $t_{n,i} \coloneqq t_n + c_i h_n,\; i = 1,\ldots, \Nstg$, are degrees of freedom. 
For a single finite element, they are summarized in the vector $V_n \coloneqq (v_{n,1}, \ldots, v_{n,\Nstg}) \in \R^{\Nstg n_x}$. 
The stage values for the algebraic variables are collected in the vectors: 
$\Theta_n \coloneqq (\theta_{n,1}, \ldots, \theta_{n,\Nstg} )\in \R^{\Nstg \cdot \Nsys}$,
$\Lambda_n \coloneqq (\lambda_{n,1}, \ldots, \lambda_{n,\Nstg} )\in \R^{\Nstg \cdot \Nsys}$ and
$M_n \coloneqq (\mu_{n,1}, \ldots, \mu_{n,\Nstg} )\in \R^{\Nstg}$. 
We also define the vector $Z_n =(x_n,\Theta_n,\Lambda_n,M_n,V_n)$ which collects all internal variables. 
With $x_n^{\mathrm{next}}$ we denote the value at $t_{n+1}$, which is obtained after a single integration step. 
Finally, the RK equations for a single finite element for the DCS \eqref{eq:dcs_1} are given by:
\begin{align}\label{eq:dcs_irk_single}
	&G_{\irk}(x_n^{\mathrm{next}}\!,Z_n,h_n,q)\! \!\coloneqq\!\!
	\begin{bmatrix}
		\! v_{n,1}\! -\!  F(x_n +h_n \sum_{j=1}^{\Nstg} a_{1,j} v_{n,j},q)\theta_{n,1}\\
		\vdots\\
		v_{n,\Nstg} \! -\! F(x_n +h_n \sum_{j=1}^{\Nstg} a_{\Nstg,j} v_{n,j},q)\theta_{n,\Nstg}\\
		G_{\LP}(x_n + h_n\sum_{j=1}^{\Nstg} a_{1,j} v_{n,j},\theta_{n,1},\lambda_{n,1},\mu_{n,1})\\
		\vdots\\
		G_{\LP}(x_n + h_n\sum_{j=1}^{\Nstg} a_{\Nstg,j} v_{n,j},\theta_{n,\Nstg},\lambda_{n,\Nstg},\mu_{n,\Nstg})\\
		x_n^{\mathrm{next}} - x_n - h_n \sum_{i=1}^{\Nstg} b_i v_{n,i}
	\end{bmatrix}=0.
\end{align}

Next, we summarize the equations for all $\NFE$ finite elements over the whole interval $[0,T]$ in a discrete-time system manner. 
For this purpose, we introduce some additional shorthands. 
All variables of all finite elements for a single control interval are collected in the vectors 
$\mathbf{x}= (x_0,\ldots,x_{\NFE}) \in \R^{(\NFE+1)n_x}$,
$\mathbf{V} = (V_0,\ldots,V_{\NFE-1}) \in \R^{\NFE \Nstg n_x}$ and
$\mathbf{h}\coloneqq (h_0,\ldots,h_{\NFE-1})\in \R^{\NFE}$.
Note that the simple continuity condition $x_{n+1} = x_{n}^{\mathrm{next}}$ holds. 
We collect all stage values of the Filippov multipliers in the vector
$\mathbf{\Theta} = ({\Theta}_0,\ldots,\Theta_{\NFE-1})\in \R^{n_{{\theta}}}$ and 
$n_{{\theta}}= \NFE\Nstg\Nsys$.
The vectors 
$\mathbf{\Lambda}\in \R^{n_{\theta}},\ 
\mathbf{M}\in \R^{n_{\mu}}$
for the stage values of the Lagrange multipliers are defined accordingly, with $n_{\mu} = \frac{n_{{\theta}}}{\Nsys}$. The vector 
$\mathbf{Z} = (\textbf{x},\mathbf{V},\mathbf{\Theta},\mathbf{\Lambda},\mathbf{M})\in \R^{n_{\mathbf{Z}}}$
collects all \textit{internal} variables and $n_{\mathbf{Z}} = (\NFE+1)n_x + \NFE \Nstg n_x + 2n_{\theta}+n_{\mu}$.

All computations over a single control interval which we call here the \textit{standard discretization} are summarized in the following equations which resemble a discrete-time system:
\begin{subequations}\label{eq:dcs_irk}
	\begin{align}
		{s}_1 \! &= \!F_{\mathrm{std}}(\textbf{Z}),\\
		0\! &=\! {G}_{\mathrm{std}}(\mathbf{Z},\mathbf{h},s_0,q),
	\end{align}
\end{subequations}
where $s_1\in \R^{n_x}$ is the approximation of $x(T)$ and
\begin{align*}
	F_{\mathrm{std}}(\textbf{Z})  &= x_{\NFE},\\
	G_{\mathrm{std}}(\mathbf{Z},\mathbf{h},s_0,q)
	\coloneqq	&
	\begin{bmatrix}
		x_0- s_0\\
		G_{\irk}(x_1,Z_0,h_0,q)\\
		\vdots\\
		G_{\irk}(x_{\NFE},Z_{\NFE-1},h_{\NFE-1},q)
	\end{bmatrix}.
\end{align*}
Note that $\mathbf{h}$ are given parameters, implicitly fixed by the chosen discretization grid. 
It is usually impossible to obtain high-accuracy solutions with this method, as this can only happen if active-set changes occur coincidentally at $t_n$. 
Despite the high accuracy in this unlikely case, the numerical sensitivities would still be wrong \cite{Nurkanovic2020,Stewart2010}. 
When active-set changes happen within a finite element, the IRK method tries to approximate a nonsmooth trajectory by a smooth polynomial, cf. the left plot in Figure \ref{fig:irk_vs_fesd}, which results in a poor approximation.

\subsection{Algorithmic ingredients of the FESD method}
\begin{figure}[t]
	\centering
	{ \pgfplotsset{compat=1.13}
\setlength{\fwidth}{4.5cm}
\setlength{\fheight}{3cm}
\pgfplotsset{
	every axis plot/.append style={line width=1.0pt},
	every axis plot post/.append style={
		every mark/.append style={line width=0.0pt}
	}
}
%
%
\definecolor{mycolor1}{rgb}{0.49400,0.18400,0.55600}%
\definecolor{mycolor2}{rgb}{0.46600,0.67400,0.18800}%
\begin{tikzpicture}

\begin{axis}[%
width=0.951\fwidth,
height=\fheight,
at={(0\fwidth,0\fheight)},
scale only axis,
xmin=0.26,
xmax=0.45,
xtick={0,0.1,0.2,0.3,0.4,0.5},
xticklabels={{},{},{},{},{},{}},
xlabel style={font=\color{white!15!black}},
xlabel={$t$},
ymin=-3.2,
ymax=-2.5,
ytick={-4,-3,-2,-1,0,1,2,3,4},
yticklabels={{},{},{},{},{},{},{},{},{}},
ylabel style={font=\color{white!15!black}},
ylabel={$x(t)$},
axis background/.style={fill=white},
axis x line*=bottom,
axis y line*=left,
legend style={at={(0.542,0.022)}, anchor=south west, legend cell align=left, align=left, draw=white!15!black,nodes={scale=0.80, transform shape}},
legend columns=2
xlabel style={font={\scriptsize}},ylabel style={font=\scriptsize},ylabel shift={-0cm},ticklabel style={font=\scriptsize}
]
\addplot [color=black,line width=0.8pt, opacity=0.3]
  table[row sep=crcr]{%
0.259183673469388	-2.5425918398291\\
0.3	-2.94300000479129\\
0.304081632653061	-3.00261228134354\\
0.308163265306122	-3.05008320637936\\
0.312244897959184	-3.08676767331227\\
0.316326530612245	-3.11393861888749\\
0.318367346938776	-3.12433190387233\\
0.320408163265306	-3.13278702318175\\
0.322448979591837	-3.13944004359584\\
0.324489795918367	-3.14442190960301\\
0.326530612244898	-3.14785844339963\\
0.328571428571429	-3.14987034489053\\
0.330612244897959	-3.15057319168858\\
0.33265306122449	-3.15007743911502\\
0.33469387755102	-3.14848842019902\\
0.336734693877551	-3.14590634567836\\
0.338775510204082	-3.14242630399869\\
0.342857142857143	-3.13312706148656\\
0.346938775510204	-3.12124895439302\\
0.351020408163266	-3.10736828777971\\
0.357142857142857	-3.08385146584322\\
0.36530612244898	-3.04941342246181\\
0.377551020408163	-2.99468982030223\\
0.387755102040817	-2.94684889406324\\
0.395918367346939	-2.9052966490452\\
0.4	-2.88250483561509\\
0.408163265306122	-2.83328048251374\\
0.416326530612245	-2.78040282939393\\
0.426530612244898	-2.70929212288369\\
0.436734693877551	-2.63276680970992\\
0.446938775510204	-2.55100313834503\\
0.451020408163265	-2.51687221166528\\
};
\addlegendentry{$x(t)$}

\addplot [color=black, dashed,line width=0.8pt]
  table[row sep=crcr]{%
0.259853228934982	-2.54916017581808\\
0.264399970762481	-2.59376371314588\\
0.269388649456113	-2.64270265113044\\
0.274377328149744	-2.69164158911501\\
0.279366006843376	-2.74058052709958\\
0.284354685537008	-2.78951946508414\\
0.28934336423064	-2.83845840306871\\
0.294332042924271	-2.88739734105328\\
0.299320721617903	-2.93633627903786\\
0.304309400311535	-2.98527521702243\\
0.309298079005166	-3.03421415500701\\
0.314286757698798	-3.08315309299159\\
0.31927543639243	-3.13209203099295\\
0.3232295867396	-3.12789198749425\\
0.326351138725635	-3.1238892647057\\
0.331371265492512	-3.11618430949195\\
0.336391392259389	-3.10692152782178\\
0.341411519026266	-3.09610792754326\\
0.346431645793143	-3.08375128700916\\
0.35145177256002	-3.06986015000077\\
0.356471899326897	-3.05444382027378\\
0.363390600874616	-3.03071634874601\\
0.369921426847973	-3.00570456232776\\
0.376552406394405	-2.97773891163839\\
0.384537084361831	-2.9406693323971\\
0.391612786695036	-2.90476236318418\\
0.400042187449235	-2.85829664639028\\
0.408571741776509	-2.80726966300944\\
0.416713420529421	-2.75487837740923\\
0.426753674063175	-2.6854321083158\\
0.436793927596929	-2.61078941963021\\
0.446834181130683	-2.53111616745583\\
0.450243455118005	-2.50294901807495\\
};
\addlegendentry{$x^*(t)$}

\addplot [color=black, dashdotted, line width=0.8pt, forget plot]
  table[row sep=crcr]{%
0.31927542840705	-3.27\\
0.31927542840705	-2.43\\
};
\addplot [color=mycolor1, draw=none, mark size=1.3pt, mark=*, mark options={solid, black}, forget plot, fill opacity=0.3]
  table[row sep=crcr]{%
0.3	-2.94300000479128\\
0.4	-2.88250483561514\\
};
\addplot [color=mycolor2, draw=none, mark size=1.3pt, mark=*, mark options={solid, black}, forget plot,fill opacity=0.3]
  table[row sep=crcr]{%
0.278765946176085	-2.73469393583587\\
0.3	-2.94300000479129\\
0.308858795951271	-3.05706071221748\\
0.340946686444073	-3.13783891358996\\
0.378765946176085	-2.98912753898159\\
0.4	-2.88250483561514\\
0.408858795951271	-2.82891643257465\\
0.440946686444073	-2.59964127056315\\
};
\addplot [color=black, dotted, line width=0.7pt, forget plot]
  table[row sep=crcr]{%
0.3	-3.27\\
0.3	-2.43\\
};
\addplot [color=black, dotted, line width=0.7pt, forget plot]
  table[row sep=crcr]{%
0.4	-3.27\\
0.4	-2.43\\
};
\end{axis}

\begin{axis}[%
width=1.227\fwidth,
height=1.227\fheight,
at={(-0.16\fwidth,-0.135\fheight)},
scale only axis,
xmin=0,
xmax=1,
ymin=0,
ymax=1,
axis line style={draw=none},
ticks=none,
axis x line*=bottom,
axis y line*=left,
legend style={legend cell align=left, align=left, draw=white!15!black},
xlabel style={font={\scriptsize}},ylabel style={font=\scriptsize},ylabel shift={-0cm},ticklabel style={font=\scriptsize}
]
\end{axis}
\end{tikzpicture}%
}
	\hspace{-1.22cm}
	\centering
	{ \pgfplotsset{compat=1.13}
\setlength{\fwidth}{4.5cm}
\setlength{\fheight}{3cm}
\pgfplotsset{
	every axis plot/.append style={line width=1.0pt},
	every axis plot post/.append style={
		every mark/.append style={line width=0.0pt}
	}
}
\definecolor{mycolor1}{rgb}{0.49400,0.18400,0.55600}%
\definecolor{mycolor2}{rgb}{0.46600,0.67400,0.18800}%
\begin{tikzpicture}

\begin{axis}[%
width=0.951\fwidth,
height=\fheight,
at={(0\fwidth,0\fheight)},
scale only axis,
xmin=0.26,
xmax=0.45,
xtick={0,0.1,0.2,0.3,0.4,0.5},
xticklabels={{},{},{},{},{},{}},
xlabel style={font=\color{white!15!black}},
xlabel={$t$},
ymin=-3.2,
ymax=-2.5,
ytick={-4,-3,-2,-1,0,1,2,3,4},
yticklabels={{},{},{},{},{},{},{},{},{}},
ylabel style={font=\color{white!15!black}},
axis background/.style={fill=white},
axis x line*=bottom,
axis y line*=left,
legend style={at={(0.727,0.824)}, anchor=south west, legend cell align=left, align=left, draw=white!15!black},
xlabel style={font={\scriptsize}},ylabel style={font=\scriptsize},ylabel shift={-0cm},ticklabel style={font=\scriptsize}
]
\addplot [color=black,line width=0.8pt, opacity=0.3]
  table[row sep=crcr]{%
0.258461061013155	-2.53550301028279\\
0.319275428310368	-3.13209195354015\\
0.322963684875463	-3.12820535838396\\
0.326651941440558	-3.12347216481139\\
0.330340198005652	-3.11789469018172\\
0.335872582853294	-3.10795104087074\\
0.341404967700936	-3.09612292545699\\
0.346937352548577	-3.08242084133052\\
0.352469737396219	-3.06685635640232\\
0.358002122243861	-3.0494421091044\\
0.363534507091503	-3.0301918083897\\
0.370911020221692	-3.00169435663739\\
0.378287533351881	-2.96999720117467\\
0.38566404648207	-2.93514171865666\\
0.393040559612259	-2.89717266911342\\
0.400417072742448	-2.85613819594994\\
0.409637714155184	-2.80061317250596\\
0.418858355567921	-2.74048782022218\\
0.428078996980657	-2.67587364842215\\
0.437299638393393	-2.60689448757036\\
0.448364408088677	-2.51854114870517\\
0.450208536371224	-2.50323912639127\\
};

\addplot [color=black, dashed,line width=0.8pt]
  table[row sep=crcr]{%
0.259853228934982	-2.54916017581808\\
0.264399970762481	-2.59376371314588\\
0.269388649456113	-2.64270265113044\\
0.274377328149744	-2.69164158911501\\
0.279366006843376	-2.74058052709958\\
0.284354685537008	-2.78951946508414\\
0.28934336423064	-2.83845840306871\\
0.294332042924271	-2.88739734105328\\
0.299320721617903	-2.93633627903786\\
0.304309400311535	-2.98527521702243\\
0.309298079005166	-3.03421415500701\\
0.314286757698798	-3.08315309299159\\
0.31927543639243	-3.13209203099295\\
0.3232295867396	-3.12789198749425\\
0.326351138725635	-3.1238892647057\\
0.331371265492512	-3.11618430949195\\
0.336391392259389	-3.10692152782178\\
0.341411519026266	-3.09610792754326\\
0.346431645793143	-3.08375128700916\\
0.35145177256002	-3.06986015000077\\
0.356471899326897	-3.05444382027378\\
0.363390600874616	-3.03071634874601\\
0.369921426847973	-3.00570456232776\\
0.376552406394405	-2.97773891163839\\
0.384537084361831	-2.9406693323971\\
0.391612786695036	-2.90476236318418\\
0.400042187449235	-2.85829664639028\\
0.408571741776509	-2.80726966300944\\
0.416713420529421	-2.75487837740923\\
0.426753674063175	-2.6854321083158\\
0.436793927596929	-2.61078941963021\\
0.446834181130683	-2.53111616745583\\
0.450243455118005	-2.50294901807495\\
};

\addplot [color=black, dashdotted, line width=0.5pt, forget plot]
  table[row sep=crcr]{%
0.31927542840705	-3.27\\
0.31927542840705	-2.43\\
};
\addplot [color=mycolor1, draw=none, mark size=1.4pt, mark=*, mark options={solid, black}, forget plot,fill opacity=0.3]
  table[row sep=crcr]{%
0.319275428310368	-3.13209195354015\\
0.409637714155184	-2.80061317250596\\
};
\addplot [color=mycolor2, draw=none, mark size=1.4pt, mark=*, mark options={solid, black}, forget plot,fill opacity=0.3]
  table[row sep=crcr]{%
0.296677056212382	-2.91040192322881\\
0.319275428310368	-3.13209195354015\\
0.327280438830265	-3.12258134493403\\
0.356275790158942	-3.05507381088572\\
0.390450137742376	-2.91085809112557\\
0.409637714155184	-2.80061317250596\\
0.41764272467508	-2.74867411695469\\
0.446638076003758	-2.53271817022897\\
};
\addplot [color=black, dotted, line width=0.7pt, forget plot]
  table[row sep=crcr]{%
0.319275428310368	-3.27\\
0.319275428310368	-2.43\\
};
\addplot [color=black, dotted, line width=0.7pt, forget plot]
  table[row sep=crcr]{%
0.409637714155184	-3.27\\
0.409637714155184	-2.43\\
};
\end{axis}

\begin{axis}[%
width=1.227\fwidth,
height=1.227\fheight,
at={(-0.16\fwidth,-0.135\fheight)},
scale only axis,
xmin=0,
xmax=1,
ymin=0,
ymax=1,
axis line style={draw=none},
ticks=none,
axis x line*=bottom,
axis y line*=left,
legend style={legend cell align=left, align=left, draw=white!15!black},
xlabel style={font={\scriptsize}},ylabel style={font=\scriptsize},ylabel shift={-0cm},ticklabel style={font=\scriptsize}
]
\end{axis}
\end{tikzpicture}%
}
	\caption{Illustration of the analytic solution and a polynomial solution approximation to a PSS via an IRK Radau-IIA method of order 7. 
		The left plot shows an approximation with a fixed step size where an active-set change happens on a stage point. 
		The right plot shows an approximation obtained with FESD (based on the same IRK method) where the switch happens on the boundary.
		The circles represent the stage values, the vertical dotted lines the finite elements boundaries, and the vertical dashed line the switching time $\ts$.}
	\label{fig:irk_vs_fesd}			
	\vspace{-0.45cm}								
\end{figure}
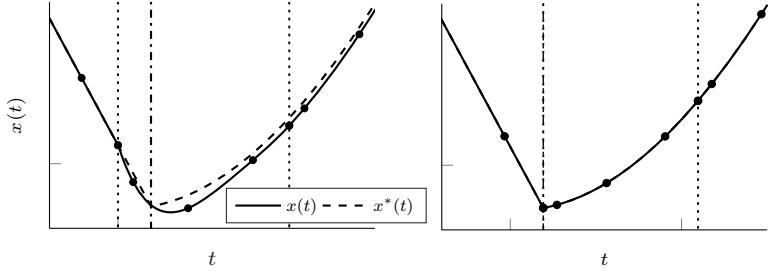

To ensure high-accuracy solutions of FESD, we allow the optimization routine to vary the lengths $h_n$ of the finite elements such that all switching points coincide with grid points $t_n$. 
Consequently, active-set changes cannot happen in the interior of each finite element, and smooth functions are approximated by smooth polynomials within a finite element, cf. the right plot in Figure \ref{fig:irk_vs_fesd}. 
Thus, the active set $\mathcal{I}(x(t))$ changes its value only at some grid point $t_n$ and is constant in the interior of all intervals $(t_{n},t_{n+1})$.  
A key assumption in any event-based method is that there are finitely many switches in finite time.
We also assume that there are enough finite elements to capture every switch that occurs in the time interval $[0,T]$. 
\subsubsection{The step sizes as degrees of freedom} 
To capture the switches with the discretization grid points $t_n$, the step sizes $h_n$ are left to be degrees of freedom in the RK method in the remainder of this paper. 
Additionally, the condition $\sum_{n=0}^{\NFE-1}h_n = T$ ensures that we regard a time interval of unaltered length.

\subsubsection{Cross complementarity}\label{sec:cross_comp}
We want to prohibit active-set changes on stage points inside a finite element. 
To achieve this, next to the complementarity conditions for every stage point $0=\Phi(\theta_{n,m},\lambda_{n,m})$ we include additional conditions on the variables $\mathbf{\Theta}$ and $\mathbf{\Lambda}$. These conditions ensure that the variable step size $h_n$ adapts so that the switching times are indeed captured, as shown below.

For ease of exposition, we assume that the underlying RK scheme satisfies $c_{\Nstg} =1$ (e.g., Radau and Lobatto methods \cite{Hairer1991}). 
This means that the right boundary point of a finite element is a stage point, since $t_{n+1} = t_n+c_{\Nstg} h_n$ for $c_{\Nstg}=1$. 
At the end of the section, we detail how to treat the case with $c_{\Nstg} \neq 1$ (e.g., Gauss-Legendre methods). 
\paragraph{Continuity of $\lambda(\cdot)$ and $\mu(\cdot)$.} 

The boundary values of the approximation of $\lambda(\cdot)$ and $\mu(\cdot)$ on an interval $[t_n,t_{n+1}]$ play a crucial role in FESD. 
Therefore, we regard their values at $t_n$ and $t_{n+1}$ which are denoted by $\lambda_{n,0},\; \mu_{n,0}$ and $\lambda_{n,\Nstg},\; \mu_{n,\Nstg}$, respectively.
We exploit the continuity of $\lambda(\cdot)$ and $\mu(\cdot)$ (cf. Lemma \ref{lem:lambda_cont}) and impose for their discrete-time counterparts for $n = 0,\ldots, \NFE-1$:
\begin{align}\label{eq:continuity_of_lambda}
	\lambda_{n,\Nstg}= \lambda_{n+1,0},\; \mu_{n,\Nstg}= \mu_{n+1,0}.
\end{align}
Therefore, in the sequel we use only the right boundary points $\lambda_{n,\Nstg}$ and $\mu_{n,\Nstg}$ which are degrees of freedom in the RK equations \eqref{eq:dcs_irk}.

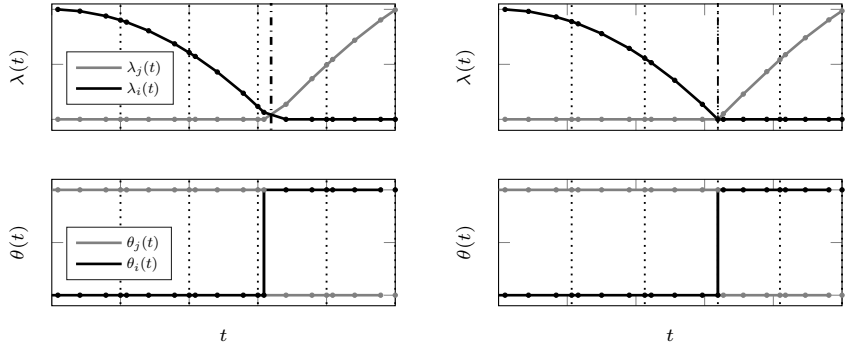
\begin{figure}[t]
	\centering
	{ \pgfplotsset{compat=1.13}
\setlength{\fwidth}{4.75cm}
\setlength{\fheight}{4cm}
\definecolor{mycolor1}{rgb}{0.92900,0.69400,0.12500}%
\definecolor{mycolor2}{rgb}{0.49400,0.18400,0.55600}%
\begin{tikzpicture}

\begin{axis}[%
width=0.951\fwidth,
height=0.419\fheight,
at={(0\fwidth,0.581\fheight)},
scale only axis,
xmin=0,
xmax=0.50,
xtick={0,0.1,0.2,0.3,0.4,0.5},
xticklabels={{},{},{},{},{},{}},
xlabel style={font=\color{white!15!black}},
ymin=-0.1,
ymax=1.04919162937417,
ytick={0,0.5,1},
yticklabels={{},{},{}},
ylabel style={font=\color{white!15!black}},
ylabel={$\lambda(t)$},
axis background/.style={fill=white},
legend style={at={(0.042,0.202)}, anchor=south west, legend cell align=left, align=left, draw=white!15!black,nodes={scale=0.70, transform shape}},
xlabel style={font={\scriptsize}},ylabel style={font=\scriptsize},ylabel shift={-0cm},ticklabel style={font=\scriptsize}
]
\addplot [color=gray, line width=1.0pt]
  table[row sep=crcr]{%
0.00885879595127037	0\\
0.308858795951271	0\\
0.340946686444074	0.136685473479308\\
0.378765946176085	0.369321591694656\\
0.408858795951271	0.544535745791255\\
0.440946686444074	0.719016057021289\\
0.478765946176085	0.903385600973894\\
0.5	0.995040564573716\\
};
\addlegendentry{$\lambda_j(t)$}

\addplot [color=black, line width=1.0pt]
  table[row sep=crcr]{%
0.00885879595127037	0.999230123213493\\
0.0409466864440735	0.983552243608376\\
0.0787659461760847	0.939138026342854\\
0.1	0.901899994996438\\
0.10885879595127	0.883749165550498\\
0.140946686444073	0.805114844776975\\
0.178765946176085	0.686499239863904\\
0.20885879595127	0.572068207741883\\
0.240946686444073	0.430477445711883\\
0.278765946176085	0.237660453021521\\
0.308858795951271	0.0635922450174007\\
0.340946686444074	0\\
0.5	0\\
};
\addlegendentry{$\lambda_i(t)$}

\addplot [color=mycolor1, draw=none, mark size=0.8pt, mark=*, mark options={solid, gray}, forget plot]
  table[row sep=crcr]{%
0.00885879595127037	0\\
0.0409466864440735	0\\
0.0787659461760847	0\\
0.1	0\\
0.10885879595127	0\\
0.140946686444073	0\\
0.178765946176085	0\\
0.2	0\\
0.20885879595127	0\\
0.240946686444073	0\\
0.278765946176085	0\\
0.3	0\\
0.308858795951271	0\\
0.340946686444074	0.136685473479308\\
0.378765946176085	0.369321591694656\\
0.4	0.493933816103132\\
0.408858795951271	0.544535745791255\\
0.440946686444074	0.719016057021289\\
0.478765946176085	0.903385600973894\\
0.5	0.995040564573716\\
};
\addplot [color=mycolor2, draw=none, mark size=0.8pt, mark=*, mark options={solid, black}, forget plot]
  table[row sep=crcr]{%
0.00885879595127037	0.999230123213493\\
0.0409466864440735	0.983552243608376\\
0.0787659461760847	0.939138026342854\\
0.1	0.901899994996438\\
0.10885879595127	0.883749165550498\\
0.140946686444073	0.805114844776975\\
0.178765946176085	0.686499239863904\\
0.2	0.607599994868011\\
0.20885879595127	0.572068207741883\\
0.240946686444073	0.430477445711883\\
0.278765946176085	0.237660453021521\\
0.3	0.117099994291006\\
0.308858795951271	0.0635922450174007\\
0.340946686444074	0\\
0.378765946176085	0\\
0.4	0\\
0.408858795951271	0\\
0.440946686444074	0\\
0.478765946176085	0\\
0.5	0\\
};
\addplot [color=black, dashdotted, line width=1.0pt, forget plot]
  table[row sep=crcr]{%
0.31927542840705	-0.214919162937417\\
0.31927542840705	1.16411079231158\\
};
\addplot [color=black, dotted, line width=0.7pt, forget plot]
  table[row sep=crcr]{%
0.1	-0.214919162937417\\
0.1	1.16411079231158\\
};
\addplot [color=black, dotted, line width=0.7pt, forget plot]
  table[row sep=crcr]{%
0.2	-0.214919162937417\\
0.2	1.16411079231158\\
};
\addplot [color=black, dotted, line width=0.7pt, forget plot]
  table[row sep=crcr]{%
0.3	-0.214919162937417\\
0.3	1.16411079231158\\
};
\addplot [color=black, dotted, line width=0.7pt, forget plot]
  table[row sep=crcr]{%
0.4	-0.214919162937417\\
0.4	1.16411079231158\\
};
\addplot [color=black, dotted, line width=0.7pt, forget plot]
  table[row sep=crcr]{%
0.5	-0.214919162937417\\
0.5	1.16411079231158\\
};
\end{axis}

\begin{axis}[%
width=0.951\fwidth,
height=0.419\fheight,
at={(0\fwidth,0\fheight)},
scale only axis,
xmin=0,
xmax=0.50,
xtick={0,0.1,0.2,0.3,0.4,0.5},
xticklabels={{},{},{},{},{},{}},
xlabel style={font=\color{white!15!black}},
xlabel={$t$},
ymin=-0.1,
ymax=1.1,
ytick={0,0.5,1},
yticklabels={{},{},{}},
ylabel style={font=\color{white!15!black}},
ylabel={$\theta(t)$},
axis background/.style={fill=white},
legend style={at={(0.042,0.202)}, anchor=south west, legend cell align=left, align=left, draw=white!15!black,nodes={scale=0.70, transform shape}},
xlabel style={font={\scriptsize}},ylabel style={font=\scriptsize},ylabel shift={-0cm},ticklabel style={font=\scriptsize}
]
\addplot[const plot, color=gray, line width=1.0pt] table[row sep=crcr] {%
0	1.000000005\\
0.2	1.000000005\\
0.20885879595127	1.000000005\\
0.240946686444073	1.00000000500001\\
0.278765946176085	1.00000000500004\\
0.3	1.00000000500006\\
0.308858795951271	0\\
0.478765946176085	0\\
};
\addlegendentry{$\theta_j(t)$}

\addplot[const plot, color=black, line width=1.0pt] table[row sep=crcr] {%
0	0\\
0.308858795951271	1.00000000500001\\
0.340946686444074	1.000000005\\
0.378765946176085	1.000000005\\
0.40885879595127	1.000000005\\
0.478765946176085	1.000000005\\
};
\addlegendentry{$\theta_i(t)$}

\addplot [color=mycolor1, draw=none, mark size=0.8pt, mark=*, mark options={solid, gray}, forget plot]
  table[row sep=crcr]{%
0.00885879595127048	1.000000005\\
0.0409466864440735	1.000000005\\
0.0787659461760848	1.000000005\\
0.1	1.000000005\\
0.10885879595127	1.000000005\\
0.140946686444074	1.000000005\\
0.178765946176085	1.000000005\\
0.2	1.000000005\\
0.20885879595127	1.000000005\\
0.240946686444073	1.000000005\\
0.278765946176085	1.00000000500001\\
0.3	1.00000000500004\\
0.308858795951271	1.00000000500006\\
0.340946686444074	0\\
0.378765946176085	0\\
0.4	0\\
0.40885879595127	0\\
0.440946686444073	0\\
0.478765946176085	0\\
0.5	0\\
};
\addplot [color=mycolor2, draw=none, mark size=0.8pt, mark=*, mark options={solid, black}, forget plot]
  table[row sep=crcr]{%
0.00885879595127048	0\\
0.0409466864440735	0\\
0.0787659461760848	0\\
0.1	0\\
0.10885879595127	0\\
0.140946686444074	0\\
0.178765946176085	0\\
0.2	0\\
0.20885879595127	0\\
0.240946686444073	0\\
0.278765946176085	0\\
0.3	0\\
0.308858795951271	0\\
0.340946686444074	1.00000000500001\\
0.378765946176085	1.000000005\\
0.4	1.000000005\\
0.40885879595127	1.000000005\\
0.440946686444073	1.000000005\\
0.478765946176085	1.000000005\\
0.5	1.000000005\\
};
\addplot [color=black, dotted, line width=0.7pt, forget plot]
  table[row sep=crcr]{%
0.1	-0.22\\
0.1	1.22\\
};
\addplot [color=black, dotted, line width=0.7pt, forget plot]
  table[row sep=crcr]{%
0.2	-0.22\\
0.2	1.22\\
};
\addplot [color=black, dotted, line width=0.7pt, forget plot]
  table[row sep=crcr]{%
0.3	-0.22\\
0.3	1.22\\
};
\addplot [color=black, dotted, line width=0.7pt, forget plot]
  table[row sep=crcr]{%
0.4	-0.22\\
0.4	1.22\\
};
\addplot [color=black, dotted, line width=0.7pt, forget plot]
  table[row sep=crcr]{%
0.5	-0.22\\
0.5	1.22\\
};
\end{axis}
\end{tikzpicture}%
}
	{ \pgfplotsset{compat=1.13}
\setlength{\fwidth}{4.75cm}
\setlength{\fheight}{4cm}
\definecolor{mycolor1}{rgb}{0.92900,0.69400,0.12500}%
\definecolor{mycolor2}{rgb}{0.49400,0.18400,0.55600}%
\begin{tikzpicture}

\begin{axis}[%
width=0.951\fwidth,
height=0.419\fheight,
at={(0\fwidth,0.581\fheight)},
scale only axis,
xmin=0,
xmax=0.50,
xtick={0,0.1,0.2,0.3,0.4,0.5},
xticklabels={{},{},{},{},{},{}},
xlabel style={font=\color{white!15!black}},
ymin=-0.1,
ymax=1.04908441530619,
ytick={0,0.5,1},
yticklabels={{},{},{}},
ylabel style={font=\color{white!15!black}},
ylabel={$\lambda(t)$},
axis background/.style={fill=white},
legend style={at={(0.751,0.665)}, anchor=south west, legend cell align=left, align=left, draw=white!15!black},
xlabel style={font={\scriptsize}},ylabel style={font=\scriptsize},ylabel shift={-0cm},ticklabel style={font=\scriptsize}
]
\addplot [color=gray, line width=1.0pt]
  table[row sep=crcr]{%
0.00942798623885344	0\\
0.319275428310368	1.15273945944239e-09\\
0.356275790158943	0.229446560538081\\
0.409637714155184	0.543371616508079\\
0.41764272467508	0.587797833778914\\
0.446638076003758	0.741140464404376\\
0.480812423587192	0.904084889580441\\
0.5	0.985860435912629\\
};

\addplot [color=black, line width=1.0pt]
  table[row sep=crcr]{%
0.00942798623885344	0.99912801457732\\
0.0435775695077397	0.981370760533064\\
0.0838267706721363	0.931065836714802\\
0.106425142770123	0.888888881803094\\
0.115853129008976	0.868330691171923\\
0.150002712277863	0.779267013825529\\
0.190251913442259	0.644919290585132\\
0.222278271779099	0.515311149617434\\
0.256427855047985	0.354941048905875\\
0.296677056212382	0.136550526239938\\
0.319275428310368	1.21291277022095e-09\\
0.5	0\\
};

\addplot [color=mycolor1, draw=none, mark size=0.8pt, mark=*, mark options={solid, gray}, forget plot]
  table[row sep=crcr]{%
0.00942798623885344	0\\
0.0435775695077397	0\\
0.0838267706721363	0\\
0.106425142770123	0\\
0.115853129008976	0\\
0.150002712277863	0\\
0.190251913442259	0\\
0.212850285540246	0\\
0.222278271779099	1.15275333723019e-09\\
0.256427855047985	1.1527534482525e-09\\
0.296677056212382	1.15275400336401e-09\\
0.319275428310368	1.15273945944239e-09\\
0.327280438830264	0.0500734148229947\\
0.356275790158943	0.229446560538081\\
0.390450137742376	0.433716708445978\\
0.409637714155184	0.543371616508079\\
0.41764272467508	0.587797833778914\\
0.446638076003758	0.741140464404376\\
0.480812423587192	0.904084889580441\\
0.5	0.985860435912629\\
};
\addplot [color=mycolor2, draw=none, mark size=0.8pt, mark=*, mark options={solid, black}, forget plot]
  table[row sep=crcr]{%
0.00942798623885344	0.99912801457732\\
0.0435775695077397	0.981370760533064\\
0.0838267706721363	0.931065836714802\\
0.106425142770123	0.888888881803094\\
0.115853129008976	0.868330691171923\\
0.150002712277863	0.779267013825529\\
0.190251913442259	0.644919290585132\\
0.212850285540246	0.555555549624267\\
0.222278271779099	0.515311149617434\\
0.256427855047985	0.354941048905875\\
0.296677056212382	0.136550526239938\\
0.319275428310368	1.21291277022095e-09\\
0.327280438830264	1.2129015569684e-09\\
0.356275790158943	1.21289955856696e-09\\
0.390450137742376	1.21289933652236e-09\\
0.409637714155184	0\\
0.41764272467508	0\\
0.446638076003758	0\\
0.480812423587192	0\\
0.5	0\\
};
\addplot [color=black, dashdotted, line width=0.7pt, forget plot]
  table[row sep=crcr]{%
0.31927542840705	-0.214908441530619\\
0.31927542840705	1.1639928568368\\
};
\addplot [color=black, dotted, line width=0.7pt, forget plot]
  table[row sep=crcr]{%
0.106425142770123	-0.214908441530619\\
0.106425142770123	1.1639928568368\\
};
\addplot [color=black, dotted, line width=0.7pt, forget plot]
  table[row sep=crcr]{%
0.212850285540246	-0.214908441530619\\
0.212850285540246	1.1639928568368\\
};
\addplot [color=black, dotted, line width=0.7pt, forget plot]
  table[row sep=crcr]{%
0.319275428310368	-0.214908441530619\\
0.319275428310368	1.1639928568368\\
};
\addplot [color=black, dotted, line width=0.7pt, forget plot]
  table[row sep=crcr]{%
0.409637714155184	-0.214908441530619\\
0.409637714155184	1.1639928568368\\
};
\addplot [color=black, dotted, line width=0.7pt, forget plot]
  table[row sep=crcr]{%
0.5	-0.214908441530619\\
0.5	1.1639928568368\\
};
\end{axis}

\begin{axis}[%
width=0.951\fwidth,
height=0.419\fheight,
at={(0\fwidth,0\fheight)},
scale only axis,
xmin=0,
xmax=0.50,
xtick={0,0.1,0.2,0.3,0.4,0.5},
xticklabels={{},{},{},{},{},{}},
xlabel style={font=\color{white!15!black}},
xlabel={$t$},
ymin=-0.1,
ymax=1.1,
ytick={0,0.5,1},
yticklabels={{},{},{}},
ylabel style={font=\color{white!15!black}},
ylabel={$\theta(t)$},
axis background/.style={fill=white},
legend style={at={(0.746,0.191)}, anchor=south west, legend cell align=left, align=left, draw=white!15!black},
xlabel style={font={\scriptsize}},ylabel style={font=\scriptsize},ylabel shift={-0cm},ticklabel style={font=\scriptsize}
]
\addplot[const plot, color=gray, line width=1.0pt] table[row sep=crcr] {%
0	1.000000005\\
0.00942798623885333	1.000000005\\
0.0435775695077396	1.000000005\\
0.115853129008976	1.000000005\\
0.190251913442259	1.000000005\\
0.212850285540246	1.00000000013278\\
0.319275428310368	0\\
0.480812423587192	0\\
};

\addplot[const plot, color=black, line width=1.0pt] table[row sep=crcr] {%
0	0\\
0.319275428310368	1.00000000011366\\
0.409637714155184	1.000000005\\
0.480812423587192	1.000000005\\
};

\addplot [color=mycolor1, draw=none, mark size=0.8pt, mark=*, mark options={solid, gray}, forget plot]
  table[row sep=crcr]{%
0.00942798623885333	1.000000005\\
0.0435775695077396	1.000000005\\
0.0838267706721363	1.000000005\\
0.106425142770123	1.000000005\\
0.115853129008976	1.000000005\\
0.150002712277862	1.000000005\\
0.190251913442259	1.000000005\\
0.212850285540246	1.000000005\\
0.222278271779099	1.00000000013278\\
0.256427855047985	1.00000000013278\\
0.296677056212382	1.00000000013278\\
0.319275428310368	1.00000000013278\\
0.327280438830264	0\\
0.356275790158943	0\\
0.390450137742376	0\\
0.409637714155184	0\\
0.41764272467508	0\\
0.446638076003758	0\\
0.480812423587192	0\\
0.5	0\\
};
\addplot [color=mycolor2, draw=none, mark size=0.8pt, mark=*, mark options={solid, black}, forget plot]
  table[row sep=crcr]{%
0.00942798623885333	0\\
0.0435775695077396	0\\
0.0838267706721363	0\\
0.106425142770123	0\\
0.115853129008976	0\\
0.150002712277862	0\\
0.190251913442259	0\\
0.212850285540246	0\\
0.222278271779099	0\\
0.256427855047985	0\\
0.296677056212382	0\\
0.319275428310368	0\\
0.327280438830264	1.00000000011366\\
0.356275790158943	1.00000000011366\\
0.390450137742376	1.00000000011366\\
0.409637714155184	1.00000000011366\\
0.41764272467508	1.000000005\\
0.446638076003758	1.000000005\\
0.480812423587192	1.000000005\\
0.5	1.000000005\\
};
\addplot [color=black, dotted, line width=0.7pt, forget plot]
  table[row sep=crcr]{%
0.106425142770123	-0.22\\
0.106425142770123	1.22\\
};
\addplot [color=black, dotted, line width=0.7pt, forget plot]
  table[row sep=crcr]{%
0.212850285540246	-0.22\\
0.212850285540246	1.22\\
};
\addplot [color=black, dotted, line width=0.7pt, forget plot]
  table[row sep=crcr]{%
0.319275428310368	-0.22\\
0.319275428310368	1.22\\
};
\addplot [color=black, dotted, line width=0.7pt, forget plot]
  table[row sep=crcr]{%
0.409637714155184	-0.22\\
0.409637714155184	1.22\\
};
\addplot [color=black, dotted, line width=0.7pt, forget plot]
  table[row sep=crcr]{%
0.5	-0.22\\
0.5	1.22\\
};
\end{axis}

\begin{axis}[%
width=1.227\fwidth,
height=1.227\fheight,
at={(-0.16\fwidth,-0.135\fheight)},
scale only axis,
xmin=0,
xmax=1,
ymin=0,
ymax=1,
axis line style={draw=none},
ticks=none,
axis x line*=bottom,
axis y line*=left,
legend style={legend cell align=left, align=left, draw=white!15!black},
xlabel style={font={\scriptsize}},ylabel style={font=\scriptsize},ylabel shift={-0cm},ticklabel style={font=\scriptsize}
]
\end{axis}
\end{tikzpicture}%
}
	\caption{
		An illustration of the standard complementarity conditions $\Psi(\mathbf{\Theta},\mathbf{\Lambda}) =0$ (left plot) and the standard complementarity conditions augmented by $0=G_{\mathrm{cross}}(\mathbf{\Theta},\mathbf{\Lambda})$ (right plot). 
		The dots represent the stage values. 
		The vertical dotted line marks the finite element boundaries, and the vertical dashed line marks the switching time $\ts$. In the standard case (left plot), an active-set change can happen at any complementarity pair. 
		With the cross complementarities \eqref{eq:cross_comp} (right plot) an active-set change can only happen on the boundaries of a finite element.}
	\label{fig:cross_comp}
\end{figure}

\paragraph{Moving the switching points to the boundary.} 
Since $\lambda(\cdot)$ is continuous, on some interval $(t_n,t_{n+1})$ with a fixed active set $\I_n$, in the interior of the regarded interval its components are either zero or positive on the whole interval. 
The stage values $\lambda_{n,i}$ of the discrete-time counterpart should satisfy this property as well. 
This is achieved by the \textit{cross complementarity} conditions, which read for all  $n \in \{0,\ldots,\NFE\!-\!1\}$ as
\begin{align}\label{eq:cross_cc_true}
	0\! &  =\!  \mathrm{diag}(\theta_{n,m})\lambda_{n,m'},\, m \in \{1,\ldots, \Nstg\},\, m' \!\in \{0,\ldots, \Nstg\}\!, \, m \! \neq m'.
\end{align}
Some of the appealing properties of the constraints \eqref{eq:cross_cc_true} are given by the next lemma. 
In our notation $\theta_{n,m,i}$ is the $i$-th component of the vector $\theta_{n,m}$. 
\begin{lemma1}\label{lem:cross_cc_statemnt}
Regard a fixed $n \in \{0,\ldots,\NFE\!-\!1\}$ and a fixed $i \in \mathcal{J}$. If any $\theta_{n,m,i}$ with $m \in \{1,\ldots, \Nstg\}$ is positive, then all $\lambda_{n,m',i}$ with $m'\in \{0,\ldots, \Nstg\}$ must be zero. Conversely, if any $\lambda_{n,m',i}$ is positive, then all $\theta_{n,m,i}$ are zero.
\end{lemma1}
\textit{Proof.} Let $\theta_{n,m,i}$ be positive, and suppose $\lambda_{n,j,i} = 0 $ and $\lambda_{n,k,i} >0$ for some $k,j\in \{0,\ldots, \Nstg\}, k\neq j$, then $\theta_{n,m,i}\lambda_{n,k,i} >0$ which violates \eqref{eq:cross_cc_true}, thus all $\lambda_{n,m',i}=0,\ m'\in \{0,\ldots, \Nstg\}$. 
The converse is proven similarly. \qed
A consequence of this lemma is that at the boundary points $t_{n+1}$ for active-set changes we have $\lambda_{n,\Nstg,i} = \lambda_{n+1,\Nstg,i} = 0$. 
This is important for the switch detection as we discuss below.
The results of the last lemma is illustrated in Figure \ref{fig:cross_comp}.
Note that in contrast to the left plot illustrating the standard complementary conditions, in the right plot, all stage points inside a finite element have the same active set and on the finite element boundary we have $\lambda_{n,\Nstg,i} = 0$.


Note that $\lambda_{0,0}$ and $\mu_{0,0}$ are not defined via Eq.~\eqref{eq:continuity_of_lambda}, as we do not have a preceding finite element.
However, they are crucial for the statement of the last lemma, especially, if the boundary point is the only stage point, as is the case for the implicit Euler method.
This can be resolved by pre-computing $\lambda_{0,0}$ explicitly and using it in \eqref{eq:cross_cc_true}. 
Note that $\lambda_{0,0}$ is not a degree of freedom.
Since $x_0$ is known, we obtain $\mu_{0,0} = \min_i g_i(x_0)$ and thus we have $\lambda_{0,0} = g(x_0) - \mu_{0,0}$.
\color{black}

The conditions \eqref{eq:cross_cc_true} are given in their sparsest form. 
Due to the non-negativity of $\Lambda_{n}$ and $\Theta_{n}$ there are many equivalent formulations of this condition, e.g., all conditions above can be summed up for a single finite element or even for all finite elements on the regarded control interval. 
Moreover, instead of the component-wise products in $\theta_{n,m}$ and $\lambda_{n, m'}$ we can use also inner products of these vectors.
Thus, we use a more compact form of \eqref{eq:cross_cc_true} where we combine the conditions for two neighboring finite elements. 
The motivation for this form is that we end up with the same number of new conditions as we have new degrees of freedom by varying $h_n$. 
The conditions read as:
\begin{align}\label{eq:cross_comp}
	\begin{split}
		&\! G_{\mathrm{cross}}(\mathbf{\Theta},\mathbf{\Lambda}) \! \coloneqq \!\\
		&\!\!\!\!
		\begin{bmatrix}
			\sum_{m = 1}^{\Nstg}\! \sum_{\substack{m'=0,\\m'\neq m}}^{\Nstg}\theta_{0,m}^\top \lambda_{0,m'}+\sum_{m = 1}^{\Nstg}\! \sum_{\substack{m'=0,\\m'\neq m}}^{\Nstg}\theta_{1,m}^\top \lambda_{1,m'}\\
			\vdots\\
			\sum_{m = 1}^{\Nstg}\! \sum_{\substack{m'=0,\\m'\neq m}}^{\Nstg}\theta_{\NFE\!-\!2,m}^\top \lambda_{\NFE\!-\!2,m'}+\sum_{m = 1}^{\Nstg}\! \sum_{\substack{m'=0,\\m'\neq m}}^{\Nstg}\theta_{\NFE\!-\!1,m}^\top \lambda_{\NFE\!-\!1,m'}
		\end{bmatrix}.\!\!
	\end{split}
\end{align}
\paragraph{Implicit switch detection.} 
We briefly explain how the switch detection for the solution approximation is realized and formalize it later in Section \ref{sec:theory}. 
Note that for $x_n^{\mathrm{next}} = x_{n+1}$ we have from the KKT conditions of the $\mathrm{LP}(x_{n+1})$ (cf. Eq.\eqref{eq:dcs_lp}) that $\mu_{n,\Nstg} = \min_j g_j(x_{n+1})$.
Moreover, if the active-set changes between the $n$-th and $n+1$-st finite element in the $i$-th component, then from Lemma \ref{lem:cross_cc_statemnt} it follows that $\lambda_{n,\Nstg,i}=0$.
Therefore, we obtain from \eqref{eq:dcs_irk_single} implicitly the condition
\begin{align*}
	g_i(x_{n+1}) &= \lambda_{n,\Nstg,i} - \mu_{n,\Nstg},
\end{align*}
which is equal to
\begin{align}\label{eq:switching_condition}
	0 &= g_i(x_{n+1})-g_j(x_{n+1})= \psi_{i,j}(x_{n+1}),
\end{align} 
where $\psi_{i,j}(x_{n+1}) = 0$ defines the switching surface between $R_i$ and $R_j$.
This condition forces $h_n$ to adapt such that the switch is detected exactly. 
Note that the condition \eqref{eq:switching_condition} appears only if active-set changes happen, hence the whole switch detection procedure is implicit.
\subsubsection{Step equilibration} \label{sec:step_equilibration}
If no switches occur, i.e., the active sets $\I_n$ do not change between two neighboring finite elements, then the cross complementarity conditions in \eqref{eq:cross_comp} are trivially satisfied. 
This yields spurious degrees of freedom in the step sizes $h_{n}$ and the optimizer can adapt the grid in an undesirable way and harm the discretization accuracy. 
Also, the path-constraint discretization can be exploited unfavorably, just to decrease the objective value. To resolve this problem we introduce \textit{step equilibration} conditions. 

The step size should only change if a switch occurs and otherwise be constant. 
This results in a piecewise uniform discretization grid for the differential and algebraic states on the regarded control interval. 
To accomplish this, we derive an indicator function that is zero only if a switch occurs otherwise its value is strictly positive. 

If some $\lambda_i(t_n)$ is equal to zero and its left or right time derivative is nonzero, then an active-set change has occurred. 
Instead of looking at the time derivatives, in the discrete-time case, we exploit the non-negativity of $\lambda_{n,m}$ and the fact that the active set is fixed for the whole finite element (due to cross complementarity, cf. Lemma \ref{lem:cross_cc_statemnt}). 
For $n \in \{1,\ldots,\NFE-1\}$, we define the following backward and forward sums of the stage values over the neighboring finite elements $[t_{n-1},t_n]$ and $[t_{n},t_{n+1}]$:
\begin{align*}
	\sigma_{n}^{\lambda,\mathrm{B}} &=  \sum_{m=0}^{\Nstg}  \lambda_{n-1,m},\
	\sigma_{n}^{\lambda,\mathrm{F}} =  \sum_{m=0}^{\Nstg}  \lambda_{n,m}.
\end{align*}

The components of $\sigma_{n}^{\lambda,\mathrm{B}}$ and $\sigma_{n}^{\lambda,\mathrm{F}}$ are zero if the left and right time derivatives of the corresponding components of $\lambda_{n,m}$ are zero.
\color{black}
Likewise, they are positive when the left and right time derivatives are nonzero.
Analogously, the sums for $\theta_{n,m}$ are defined as:
\begin{align*}
	\sigma_{n}^{\theta,\mathrm{B}} &=  \sum_{m=1}^{\Nstg}  \theta_{n-1,m},\	
	\sigma_{n}^{\theta,\mathrm{F}} =  \sum_{m=1}^{\Nstg}   \theta_{n,m}.
\end{align*}
Additionally, we define the following vectors for all $n \in \{ 1,\dots,\NFE-1$\}:
\begin{align*}
	\pi_{n}^{\lambda}  =\mathrm{diag}(\sigma_{n}^{\lambda,\mathrm{B}})	\sigma_{n}^{\lambda,\mathrm{F}},\
	\pi_{n}^{\theta}   = \mathrm{diag}(\sigma_{n}^{\theta,\mathrm{B}})	\sigma_{n}^{\theta,\mathrm{F}}.
\end{align*}
If there is an active-set change in the $i$-th complementarity pair, then at most one of the $i$-th components of $\sigma_{n}^{\lambda,\mathrm{B}}$ and $\sigma_{n}^{\lambda,\mathrm{F}}$ is nonzero, hence their product, i.e., the  $i$-th component of $\pi_{n}^{\lambda}$, is zero. Due to complementarity, the same holds for $\pi_{n}^{\theta}$. 
For sliding modes the corresponding components of $\pi_{n}^{\lambda}$ are zero and of $\pi_{n}^{\theta}$ they are positive (due to complementarity). Thus, the $i$-th component of 
\begin{align*}
	\upsilon_n  = 	\pi_{n}^{\lambda} +	\pi_{n}^{\theta},
\end{align*}
is only zero if there is an active-set change in the $i$-th complementarity pair at $t_n$. A function that has the desired properties is defined as:
\begin{align*}
	\eta_n(\mathbf{\Theta},\mathbf{\Lambda}) &\coloneqq \prod_{i=1}^{\Nsys} (\upsilon_n)_i.
\end{align*}
This scalar function summarizes the effects of all components. It is zero only if an active-set change happens at the boundary point $t_{n}$, otherwise, it is strictly positive. 
Finally, the constraints that remove possible spurious degrees of freedom in $h_n$ read as: 
\begin{align}\label{eq:step_eq}
	0\!&=\!G_{\mathrm{eq}}(\mathbf{h},\mathbf{\Theta},\mathbf{\Lambda}) 
	\!\!	\coloneqq\!\!\!
	\begin{bmatrix}
		(h_{1}-h_{0})\eta_1(\mathbf{\Theta},\!\mathbf{\Lambda}) \\
		\vdots\\
		(h_{\NFE\!-\!1}-h_{\NFE\!-\!2})\eta_{\NFE\!-\!1}(\mathbf{\Theta},\!\mathbf{\Lambda}) 
	\end{bmatrix}\!\!.
\end{align} 
Since many products are involved in $\eta_n(\mathbf{\Theta},\mathbf{\Lambda})$, one can replace it by $\tilde{\eta}_n(\mathbf{\Theta},\mathbf{\Lambda}) \coloneqq  \tanh(\eta_n(\mathbf{\Theta},\mathbf{\Lambda}))$  to have a better scaling. 
An example for {step equilibration} is studied in Subsection \ref{sec:step_equilb_example} numerically.  

\subsubsection{The FESD discretization}
We have now all the ingredients to extend the standard RK discretization \eqref{eq:dcs_irk} to the \textit{FESD discretization}. 
We use again the same discrete-time representation
 \begin{subequations} \label{eq:fesd_compact}
 	\begin{align}
 		s_{1} \!&=\! F_{\fesd}(\mathbf{Z}), \label{eq:fesd_compact_state_transition}\\ 
 		0 \!&= \!G_{\fesd}(\mathbf{Z},\mathbf{h},s_0, q , T),
 	\end{align}
 \end{subequations}
 where $F_{\fesd}(\mathbf{x})\!=x_{\NFE}$ is the state transition map and  $G_{\fesd}(\mathbf{x},\mathbf{h},\mathbf{Z},q, T)$ collects all other internal computations including all RK steps within the regarded control interval:
 \begin{align*}
 	&G_{\fesd}(\mathbf{Z},\mathbf{h},s_0,q, T)\! \coloneqq\!\!
 	\begin{bmatrix}
 		{G}_{\mathrm{std}}(\mathbf{Z},\mathbf{h},s_0,q)\\
 		G_{\mathrm{cross}}(\mathbf{\Theta},\mathbf{\Lambda})\\
 		G_{\mathrm{eq}}(\mathbf{h},\mathbf{\Theta},\mathbf{\Lambda})\\
		\sum_{n=0}^{\NFE-1} h_n - T			
 	\end{bmatrix}.
 \end{align*}
For a fixed control function $q$, horizon length $T$ and initial value $s_0$, the formulation \eqref{eq:fesd_compact} can be used as an integrator with exact switch detection for PSS~\eqref{eq:pws1}.
Since Filippov DI does not always have unique solutions, one cannot expect uniqueness of solutions for their discrete-time counterparts \eqref{eq:fesd_compact} in all cases. 
In simulation methods, a common approach is to either pick one \textit{local} solution obtained by the solver for the nonlinear complementarity problem \eqref{eq:fesd_compact} or to enumerate all possible solutions at an active-set change \cite{Acary2014,Stewart1990b}. 
In this paper, we consider only the first option. 
Note that in sliding modes, we implicitly obtain differential algebraic equations of index 2, cf. Section \ref{sec:stewarts_dcs_active_set_fixed}. To achieve good accuracy in practice it is usually required to use stiffly accurate methods, e.g., Radau-IIA methods \cite{Hairer1991}.

\subsubsection{Remark on RK methods with $c_{\Nstg}\neq 1$}
We outline how to extend the FESD method when an RK scheme with $c_{\Nstg}\neq1$ is regarded. 
In contrast to the developments so far, with $c_{\Nstg} \neq 1$ the variables $\lambda_{n,\Nstg},\; \mu_{n,\Nstg}$ do not correspond the boundary values $\lambda(t_{n+1})$ and $\mu(t_{n+1})$ anymore (since $t_n+c_{\Nstg} h_n < t_{n+1}$). 
We denote the boundary points in this case by  $\lambda_{n,\Nstg+1},\; \mu_{n,\Nstg+1}$. 
They are computed by solving $\LP(x_{n+1})$ for $n=0,\ldots
\NFE-2$:
\begin{align}\label{eq:lp_boundary_points}
	0 &= G_{\mathrm{LP}}(x_{n+1},\theta_{n,\Nstg+1},\lambda_{n,\Nstg+1},\mu_{n,\Nstg+1}).
\end{align}
We still exploit the continuity of $\lambda(\cdot)$ and $\mu(\cdot)$ (cf. Lemma \ref{lem:lambda_cont}), by replacing~\eqref{eq:continuity_of_lambda} with the following continuity conditions for their discrete-time counterparts for $n = 0,\ldots, \NFE-1$:
\begin{align}\label{eq:continuity_of_lambda_2}
	\lambda_{n,\Nstg+1}= \lambda_{n+1,0},\; \mu_{n,\Nstg+1}= \mu_{n+1,0}.
\end{align}
With slight abuse of notation, we add the new variables $\theta_{n,\Nstg+1},\lambda_{n,\Nstg+1}$ and $\mu_{n,\Nstg+1}$ to the vectors $\mathbf{\Theta}$, $\mathbf{\Lambda}$ and $\mathbf{M}$, respectively. 
The vector $\mathbf{Z}$ is redefined accordingly. 
The cross complementarity conditions are now modified such that next to the stage points we include the boundary points with the index $\Nstg+1$:
 \begin{align*}
 	\begin{split}
 		&\! G_{\mathrm{cross}}(\mathbf{\Theta},\mathbf{\Lambda}) \! \coloneqq \!\\
 		&\!\!\!\!
 		\begin{bmatrix}
 			\sum_{m = 1}^{\Nstg}\! \sum_{m'=1,m'\neq m}^{\Nstg+1}\theta_{0,m}^\top \lambda_{0,m'}+\sum_{m = 1}^{\Nstg}\! \sum_{m'=0,m'\neq m}^{\Nstg+1}\theta_{1,m}^\top \lambda_{1,m'}\\
 			\vdots\\
 			\sum_{m = 1}^{\Nstg}\! \sum_{\substack{m'=0,\\m'\neq m}}^{\Nstg+1}\theta_{\NFE\!-\!2,m}^\top \lambda_{\NFE\!-\!2,m'}+\sum_{m = 1}^{\Nstg}\! \sum_{\substack{m'=0,\\m'\neq m}}^{\Nstg}\theta_{\NFE\!-\!1,m}^\top \lambda_{\NFE\!-\!1,m'}
 		\end{bmatrix}.\!\!
 	\end{split}
 \end{align*}
For the whole control time we have in total $(\NFE-1)(2\Nsys+1)$ new variables.
\section{Convergence theory}\label{sec:theory}
In this section we present the main convergence result of the FESD method. 
First, we prove that even though the FESD system \eqref{eq:fesd_compact} \sloppy{is always over-determined it still has a locally isolated solution}. 
Second, we show that the numerical solution approximation $\hat{x}_h(\cdot)$ generated by FESD converges to a solution $x(\cdot)$ in the sense of Definition \ref{def:piecewise_active},  with the same order that the underlying RK method has for smooth ODE. 
Additionally, we prove that the numerical sensitivities converge to their correct values with high accuracy.
\subsection{Main assumptions}\label{sec:assumptions}
We start by introducing some notation and stating some assumptions related to the FESD formulation \eqref{eq:fesd_compact}, which are important for our theoretical study in this section. 

\begin{assumption}(Runge-Kutta method)
	\label{ass:irk_scheme}
	A Butcher tableau with the entries $a_{i,j} ,b_i$ and $c_i$, $i,j\in\{1,\ldots,\Nstg\}$ related to an $\Nstg$-stage Runge-Kutta (RK) method 
	is used in the FESD \eqref{eq:fesd_compact}. Moreover, we assume that:
	\begin{enumerate}[(a)]
		\item If the same RK method is applied to the differential algebraic equation \eqref{eq:dcs_dae} on an interval $[t_a,t_b]$, it has a global accuracy of $O(h^p)$ for the differential states.
		\item The RK equations applied to \eqref{eq:dcs_dae} have a locally isolated solution for a sufficiently small $h_n>0$. 
	\end{enumerate}
\end{assumption}
This assumption aims to consider a broad class of RK methods, and both assumptions are standard assumptions \cite{Hairer1991}.
\begin{assumption}(Solution existence)
	\label{ass:solution_existence_fesd}
	For given parameters $s_0,q$ and $T$, there exists a solution to the FESD problem \eqref{eq:fesd_compact}, such that for all $n \in \{0,\ldots,\NFE-1\}$ it holds that {${h}_n\geq0$}.
\end{assumption}
This assumption means that there exists a solution and that we can compute it. 
If the FESD method is used in direct optimal control, non-negativity of the step sizes can easily be achieved by adding box constraints on $h_n$.
{
This is the strongest assumption we make in this paper. 
Ideally, one would prove the existence of solutions. 
Since the system is over-determined this cannot be done straightforwardly by applying standard existence results~\cite{Facchinei2003}. 
As we will show below, in practice numerical solvers have no trouble computing such solutions. 
}

We state a technical assumption that ensures regularity of the FESD problem \eqref{eq:fesd_compact}. 
\begin{assumption}(Regularity)
	\label{ass:regularity}
	Given the complementarity pairs $\Psi(\theta_{n,m},\lambda_{n,m})=0$, for all $n = 0,\ldots\NFE-1$ there exists an $m\in\{1,\dots,\Nstg\}$ and $i \in \{1,\ldots,\Nsys\}$, such that the strict complementarity property holds, i.e., $\theta_{n,m,i}+\lambda_{n,m,i}>0$. 
	Moreover, for the RK equations \eqref{eq:dcs_irk_single} it holds for all $n = 0,\ldots\NFE-1$, that at least one entry of the vector $\nabla_{h_n} G_{\irk}(x_{n+1},Z_n,h_n,q)$ is nonzero.
\end{assumption}

Once all stage values are computed by solving \eqref{eq:fesd_compact}, we can use some interpolation method to construct the solution approximation candidate in continuous time, cf. Assumption \ref{ass:irk_scheme}. 
For example, if we use a collocation-based IRK method continuous-time approximation $\hat{x}_n(t;h_{n})$ on every finite element is easily obtained via Lagrange polynomials \cite{Hairer1991}. We append the approximation for every finite element and write
\begin{align}\label{eq:continious_time_fesd}
	\hat{x}_h(t) &= \hat{x}_n(t;h_{n})\ \text{if}\ t\in [t_n,t_{n+1}],
\end{align}
where $h = \max_{n\in\{0,\ldots \NFE-1\}} h_n$. 
Similarly, continuous-time representations can be found for the algebraic variables, and we denote them compactly as $\hat{\lambda}_h(t)$, $\hat{\theta}_h(t)$ and $\hat{\mu}_h(t)$.
Similar to the definitions in Section \ref{sec:stewarts_dcs_active_set_changes}, the fixed active set in this case is denoted by $\I(\hat{x}_h(t)) = \hat{\I}_n,\ t\in(\hat{t}_{\mathrm{s},n},\hat{t}_{\mathrm{s},n+1})$ and the active set at switching point $\hat{t}_{\mathrm{s},n}$ by $\I(\hat{x}_h(\hat{t}_{\mathrm{s},n})) = \hat{\I}_n^0$.

\subsection{Solutions of the FESD problem are locally isolated}

In this subsection, we analyze some properties of solutions of the FESD problem~\eqref{eq:fesd_compact}.
For the  convenience of the reader, we restate the problem but discard the trivial state transition map $s_1 = F_{\fesd}(\mathbf{Z}) = x_{\NFE}$:
 \begin{align}\label{eq:fesd_equation}
	G_{\fesd}(\mathbf{Z},\mathbf{h},s_0,q,T)&= 
	\begin{bmatrix}
		{G}_{\mathrm{std}}(\mathbf{Z},\mathbf{h},s_0,q,T)\\
		G_{\mathrm{cross}}(\mathbf{\Theta},\mathbf{\Lambda})\\
		G_{\mathrm{eq}}(\mathbf{h},\mathbf{\Theta},\mathbf{\Lambda})\\
		\sum_{n=0}^{\NFE-1} h_n - T			
	\end{bmatrix} = 0.
\end{align}
Recall that $\mathbf{Z} = (\textbf{x},\mathbf{V},\mathbf{\Theta},\mathbf{\Lambda},\mathbf{M})\in \R^{n_{\mathbf{Z}}}$.
Additionally, we have that
${G}_{\mathrm{std}}: \R^{n_{\mathbf{Z}}} \times \R^{\NFE} \times \R^{n_x} \times \R^{n_u} \times \R \rightarrow \R^{n_{\mathbf{Z}}}$, 
$G_{\mathrm{cross}}:  \R^{n_{\theta}} \times \R^{n_{\theta}} \rightarrow \R^{\NFE-1}$ and
$G_{\mathrm{eq}}: \R^{n_{\NFE}} \times  \R^{n_{\theta}} \times \R^{n_{\theta}} \rightarrow \R^{\NFE-1}$. 
Finally, we have that 
$G_{\fesd}: \R^{n_{\mathbf{Z}}} \times \R^{\NFE} \times \R^{n_x} \times \R^{n_u} \times \R \to \R^{n_{\mathbf{Z}} + 2\NFE-1}$. 
Again, for ease of exposition, we regard $c_{\Nstg} = 1$ and give the extensions later with $n_{\theta}= \NFE\Nstg\Nsys$ and  $n_{\mu} =  \NFE\Nstg$.

The vectors $s_0 \in \R^{n_x}$, $q \in \R^{n_u}$ and $T\in \R$ are given parameters, hence we have $n_{\mathbf{Z}} + \NFE$ unknowns and $n_{\mathbf{Z}} + 2\NFE-1$ equations. 
Consequently, for $\NFE>1$, which we always assume in FESD, the system \eqref{eq:fesd_equation} is over-determined.
However, we show in the next theorem that for a given active set $\NFE-1$ equations in \eqref{eq:fesd_equation} are implicitly satisfied, and we always end up with a square system. 
As a consequence, Eq. \eqref{eq:fesd_equation} has under reasonable assumptions a locally unique solutions. 
Nevertheless, since we do not know the active set a priori, we can also not know which equations are binding and which are implicitly satisfied.

\begin{lemma1}[Corollary 6.1 in \cite{Matsaglia1974}]\label{lem:rank_of_product}
	Let $A_1 \in \R^{k \times m}$ and $A_2 \in \R^{m \times q}$, then
	\begin{align*}
		\mathrm{rank}(A_1)+	\mathrm{rank}(A_2) - m \leq 	\mathrm{rank}(A_1 A_2) \leq \min(\mathrm{rank}(A_1),	\mathrm{rank}(A_2)).
	\end{align*}
\end{lemma1}

\begin{theorem1}\label{th:fesd_locally_unique_solution}
	Suppose that Assumptions \ref{ass:irk_scheme}, \ref{ass:solution_existence_fesd} and \ref{ass:regularity} hold. 
	Let ${s}_0$, ${q}_0$  and ${T}>0$ be some fixed parameters such that $G_\fesd(\mathbf{Z}^*,\mathbf{h}^*,{s}_0,{q},{T}) =0$.
	Let $P^*\subseteq \R^{n_x}\times \R^{n_u} \times \R $ be the set of all parameters $(\hat{s}_0,\hat{q},\hat{T})$ such that $\mathbf{Z} \in \R^{n_\mathbf{Z}}$, which is the solution of $G_\fesd(\mathbf{Z},\mathbf{h},\hat{s}_0,\hat{q},\hat{T})=0$, has the same active set as $\mathbf{Z}^*$. 
	Additionally, suppose that $G_\fesd(\cdot)$ is continuously differentiable in $s_0,q$ and $T$ for all $(s_0,q,T) \in P^*$. 
	Then there exists a neighborhood ${P} \subseteq P^*$ of $({s}_0,{q}_0,{T})$ such that there exist continuously differentiable single valued functions  $\mathbf{Z}^{*}: {P} \to  \R^{n_{\mathbf{Z}}}$ and $\mathbf{h}^*: {P} \to \R^{\NFE}$.
\end{theorem1}
\textit{Proof.}
We regard the active sets for every finite element $\hat{\I}_n$ for all $n \in \{0\ldots,\NFE-1\}$ that correspond to the solution $(\mathbf{Z}^*,\mathbf{h}^*)$.
First, we look closer at the equations $G_{\mathrm{cross}}(\mathbf{\Theta}^*,\mathbf{\Lambda}^*)=0$ and $G_{\mathrm{eq}}(\mathbf{h}^*,\mathbf{\Theta}^*,\mathbf{\Lambda}^*)=0$.
If two neighboring finite elements have the same active set, i.e., $\hat{\I}_{n} = \hat{\I}_{n+1}$, then the $(n+1)$-th entry of $G_{\mathrm{cross}}(\mathbf{\Theta}^*,\mathbf{\Lambda}^*)$ is implicitly satisfied due to the point-wise complementarity conditions $\Psi(\Theta_n,\Lambda_n)=0$ and $\Psi(\Theta_{n+1},\Lambda_{n+1})=0$. Moreover, by construction we have $\eta_{n+1}(\mathbf{\Theta}^*,\mathbf{\Lambda}^*)>0$ and the $(n+1)$-th entry of $G_{\mathrm{eq}}(\mathbf{h}^*,\mathbf{\Theta}^*,\mathbf{\Lambda}^*,{T})=0$ is binding, i.e., it implies $h^*_{n+1}=h^{*}_{n}$.
On the other hand, if $\hat{\I}_n \neq \hat{\I}_{n+1}$, we have by construction that $\eta_{n+1}(\mathbf{\Theta}^*,\mathbf{\Lambda}^*)=0$ and then $(n+1)$-th entry of $G_{\mathrm{eq}}(\mathbf{h}^*,\mathbf{\Theta}^*,\mathbf{\Lambda}^*,{T})=0$ vanishes, i.e., is satisfied for any $h^{*}_{n}$ and $h^{*}_{n+1}$.
However, the $(n+1)$-th entry of $G_{\mathrm{cross}}(\mathbf{\Theta}^*,\mathbf{\Lambda}^*)=0$ is now binding, cf. Lemma \ref{lem:cross_cc_statemnt}. 

We collect the binding $n_1$ cross complementarity conditions, with $0\leq n_1 \leq \NFE-1$, in the equation $G^*_{\mathrm{cross}}(\mathbf{\Theta}^*,\mathbf{\Lambda}^*)=0$, and the $\NFE-1-n_1$ implicitly satisfied into $G^{\mathrm{res}}_{\mathrm{cross}}(\mathbf{\Theta}^*,\mathbf{\Lambda}^*) = 0$. 
Likewise, we collect the binding $n_2$ step equilibration conditions, with $1\leq n_2 \leq \NFE-1$, in  $G^*_{\mathrm{eq}}(\mathbf{h}^*,\mathbf{\Theta}^*,\mathbf{\Lambda}^{*}) = 0$.
The remaining $\NFE-1-n_2$ conditions are implicitly satisfied and are collected in $G^{\mathrm{res}}_{\mathrm{eq}}(\mathbf{h}^*,\mathbf{\Theta}^*,\mathbf{\Lambda}^*)=0$. 
Note that $n_1+n_2 = \NFE-1$.
We highlight that $\sum_{n=0}^{\NFE-1} h_n - T$ is always binding. 

We can further simplify our system of equations by eliminating some degrees of freedom using 
$G^*_{\mathrm{eq}}(\mathbf{h}^*,\mathbf{\Theta}^*,\mathbf{\Lambda}^*) = 0$.
All components of this vector are of the form $\eta_n (h_{n}-h_{n+1})$ with $\eta_n >0$. 
Therefore, we have $n_2$ equations of the form of $h_{n} = h_{n+1}$ and can remove $n_2$ degrees of freedom. 
Furthermore, we can express any $h_{j} = T- \sum_{i=0, i\neq j}^{\NFE-1} h_n$ and remove another degree of freedom.
In total we removed $n_2+1$ degrees of freedom and can regard a reduced number of unknown step-sizes, which we denote by $\mathbf{\tilde{h}}^{*} \in \R^{n_1}, n_1 = \NFE-n_2-1$.
With a slight abuse of notation, we redefine the standard RK equations accordingly and obtain
${G}_{\mathrm{std}}(\mathbf{Z}^*,\mathbf{\tilde{h}}^*,{s}_0,{q},T) = 0$ with
${G}_{\mathrm{std}}: \R^{n_{\mathbf{Z}}} \times \R^{n_1} \times \R^{n_x} \times \R^{n_u} \times \R \to
\R^{n_{\mathbf{Z}}}$.

To summarize, for a fixed active set we can rewrite \eqref{eq:fesd_equation} in a reduced form as
\begin{align}\label{eq:fesd_reduced_binding_stewart}
	G_{\fesd}^*(\mathbf{Z}^*,\mathbf{\tilde{h}}^*,{s}_0,{q},{T})
	& \coloneqq
	\begin{bmatrix}
		{G}_{\mathrm{std}}(\mathbf{Z}^*,\mathbf{\tilde{h}}^*,{s}_0,{q},T)\\
		G^*_{\mathrm{cross}}(\mathbf{\Theta}^*,\mathbf{\Lambda}^*)
	\end{bmatrix} = 0,
\end{align}
with $G_{\fesd}^*(\mathbf{Z}^*,\mathbf{h}^*,{s}_0,{q},{T}) \in \R^{n_{\mathbf{Z}}+n_1}$.
These conditions imply 
\begin{align}\label{eq:fesd_reduced_implied_stewart}
	G_{\fesd}^{\mathrm{res}}(\mathbf{h}^*,\mathbf{\Theta}^*,\mathbf{\Lambda}^*)&\coloneqq 
	\begin{bmatrix}
		G^{\mathrm{res}}_{\mathrm{cross}}(\mathbf{\Theta}^*,\mathbf{\Lambda}^*)\\
		G^{\mathrm{res}}_{\mathrm{eq}}(\mathbf{h}^*,\mathbf{\Theta}^*,\mathbf{\Lambda}^*)
	\end{bmatrix}=0,
\end{align}
with $G_{\fesd}^{\mathrm{res}}(\mathbf{h}^*,\mathbf{\Theta}^*,\mathbf{\Lambda}^*) \in \R^{\NFE-1}$.
Thus, for a given active set we can discard \eqref{eq:fesd_reduced_implied_stewart} and regard only the equivalent reduced problem \eqref{eq:fesd_reduced_binding_stewart}, which is a square system of equations.

Next, we show that the Jacobian matrix $\nabla_{(\mathbf{Z},\mathbf{\tilde{h}})} G_{\fesd}^*(\mathbf{Z}^*,\mathbf{\tilde{h}}^*,{s}_0,{q},{T})^\top$ has full rank. 
This enables us to apply the implicit function theorem (cf. \cite[Theorem 1B.1]{Dontchev2014}) and establish the result of this theorem. 
We take a closer look at the matrix:
\begin{align*}
	\begin{split}
		&\nabla_{(\mathbf{Z},\mathbf{\tilde{h}})} G_{\fesd}^*(\mathbf{Z}^*,\mathbf{\tilde{h}}^*,{s}_0,{q},{T})^{\top}\\
		&\!=\!\begin{bmatrix}
			\nabla_{\mathbf{Z}} G_{\mathrm{std}}(\mathbf{Z}^*,\mathbf{\tilde{h}}^*,{s}_0,{q},{T})^{\top} 
			& \nabla_{\mathbf{\tilde{h}}} G_{\mathrm{std}}(\mathbf{Z}^*,\mathbf{\tilde{h}}^*,{s}_0,{q},{T})^{\top}\\
			\nabla_{\mathbf{Z}} G_{\mathrm{cross}}^*(\mathbf{Z}^*,\mathbf{\tilde{h}}^*,{s}_0,{q},{T})^{\top}
			& \nabla_{\mathbf{\tilde{h}}} G_{\mathrm{cross}}^*(\mathbf{Z}^*,\mathbf{\tilde{h}}^*,{s}_0,{q},{T})^{\top}
		\end{bmatrix}.
	\end{split}
\end{align*}
Under Assumption \ref{ass:solution_existence_fesd}, for a fixed active set and a fixed $h_n^{*}$ the equation ${G}_{\mathrm{std}}(\mathbf{Z}^*,\mathbf{\tilde{h}}^*,{s}_0,{q},{T}) = 0$ boils down to the RK equations for the differential algebraic equation \eqref{eq:dcs_dae}. 
\sloppy{Due to Assumption \ref{ass:irk_scheme} the RK system ${G}_{\mathrm{std}}(\mathbf{Z}^*,\mathbf{\tilde{h}}^*,{s}_0,{q},{T}) = 0$ has a locally isolated solution}. 
A necessary and sufficient condition for this property is the invertibility of the Jacobian 
$ \nabla_{\mathbf{Z}} G_{\mathrm{std}}(\mathbf{Z}^*,\mathbf{\tilde{h}}^*,{s}_0,{q},{T})^\top$  \cite[Theorem 1B.8]{Dontchev2014}. 
Thus, we have that
$\mathrm{rank}( \nabla_{\mathbf{Z}} G_{\fesd}^*(\mathbf{Z}^*,\mathbf{\tilde{h}}^*,{s}_0,{q},{T})^\top) = n_{\mathbf{Z}}$. 
Second, due to the block diagonal structure of $ \nabla_{\mathbf{\tilde{h}}} G_{\mathrm{std}}^*(\mathbf{Z}^*,\mathbf{\tilde{h}}^*,{s}_0,{q},{T})$ and Assumption \ref{ass:regularity} we can deduce that 
$\mathrm{rank}( \nabla_{\mathbf{\tilde{h}}} G_{\fesd}^*(\mathbf{Z}^*,\mathbf{\tilde{h}}^*,{s}_0,{q},{T})^\top) = n_1$.
Third, due to the nonnegativity of $(\mathbf{\Theta},\mathbf{\Lambda})$ and Assumption \ref{ass:regularity} by direct computation it can be verified that $\mathrm{rank}(\nabla_{\mathbf{Z}} G^*_{\mathrm{cross}}(\mathbf{\Theta},\mathbf{\Lambda})^\top) = n_1$ and 
$\nabla_{\mathbf{\tilde{h}}} G^*_{\mathrm{cross}}(\mathbf{\Theta},\mathbf{\Lambda})^\top = 0$.

We introduce more compact notation and summarize the results so far with:
\begin{itemize}
	\item $M_1 = \nabla_{\mathbf{Z}} G_{\mathrm{std}}(\mathbf{Z}^*,\mathbf{\tilde{h}}^*,{s}_0,{q},{T})^\top \in \R^{n_{\mathbf{Z}}  \times n_{\mathbf{Z}}}$ with $\mathrm{rank}(M_1)  = n_{\mathbf{Z}}$
	\item $M_2 = \nabla_{\mathbf{\tilde{h}}} G_{\mathrm{std}}(\mathbf{Z}^*,\mathbf{\tilde{h}}^*,{s}_0,{q},{T})^\top \in \R^{n_{\mathbf{Z}}  \times n_1 }$ with $\mathrm{rank}(M_2)  = n_{1}$ and
	\item 	$M_3 = \nabla_{\mathbf{Z}} G_{\mathrm{cross}}(\mathbf{\Theta},\mathbf{\Lambda})^\top \in \R^{n_1 \times n_{\mathbf{Z}}}$ with $\mathrm{rank}(M_3)  = n_{1}$.
\end{itemize}

To show that $\nabla_{(\mathbf{Z},\mathbf{\tilde{h}})} G_{\fesd}^*(\mathbf{Z}^*,\mathbf{\tilde{h}}^*,{s}_0,{q},{T})^\top$ has a rank of $n_{\mathbf{Z}}+n_1$, we show that the  linear system
\begin{align*}
	\begin{bmatrix}
		M_1 &M_2 \\
		M_3 &0
	\end{bmatrix} 
	\begin{bmatrix}
		v \\ w
	\end{bmatrix} = 0,
\end{align*}
with $v \in \R^{n_\mathbf{Z}}$ and $w \in \R^{n_1}$ has zero as the only solution. 

From the first line in this linear system, we have that $v = -M_1^{-1} M_2 w$.
Since $ n_{\mathbf{Z}} > n_1$, from Lemma \ref{lem:rank_of_product}, we conclude that 
$\mathrm{rank}(M_1^{-1} M_2 ) = n_1$.
Next, from the second part of our linear system, we have that 
$ - M_3 M_1^{-1} M_2 w = 0$. 
Again, using Lemma \ref{lem:rank_of_product}, we conclude that $\mathrm{rank}(M_3 M_1^{-1} M_2 ) = n_1$.
Hence, we have $w = 0$ and $v=0$ to be the only solution of the regarded linear system.
This completes the proof.\qed 

\color{black}
\begin{remark1}
We note that one cannot apply more general forms of implicit function theorems for generalized and nonsmooth equations \cite{Dontchev2014}. 
They usually require Lipschitz continuity of the solution map to reason about local uniqueness, but the solution map for FESD is not continuous in general, but only piecewise continuous.
\end{remark1}
\begin{example1}
To illustrate the discontinuity of the solution map, we look at the example of $\dot{x} \in  2-\mathrm{sign}(x) +x^2$, with $\NFE = 2$, $T = 0.2$ and vary $x_0\in [-0.7,0.1]$. 
A solution approximation is obtained via FESD based on the Radau-IIA method of order 3.
Consider an initial value $x_0$ such that no switch occurs and a perturbed initial value $x_0+\epsilon$ where a single switch occurs on the time interval of interest. 
Clearly, in the first case, we have an equidistant grid with $h_0=h_1$, and in the second case $h_0$ jumps to $\hat{t}_{\mathrm{s},1}$. 
We conclude that $h_0(x_0)$ is not a Lipschitz function, see Figure \ref{fig:non_lipsitz} for an illustration.
\end{example1}
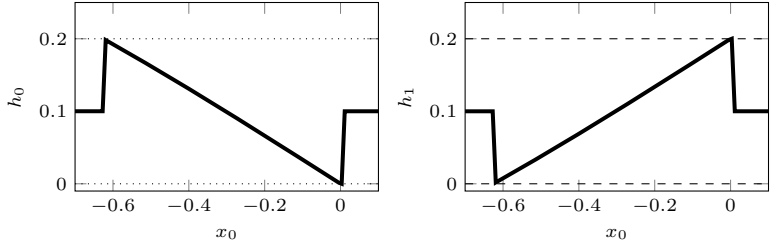
\begin{figure}[t]
	\centering
	{ \pgfplotsset{compat=1.13}
\setlength{\fwidth}{9.5cm}
\setlength{\fheight}{2.5cm}
\begin{tikzpicture}

\begin{axis}[%
width=0.420\fwidth,
height=\fheight,
at={(0\fwidth,0\fheight)},
scale only axis,
xmin=-0.7,
xmax=0.1,
xlabel style={font=\color{white!15!black}},
xlabel={$x_0$},
ymin=-0.01,
ymax=0.25,
ylabel style={font=\color{white!15!black}},
ylabel={$h_0$},
axis background/.style={fill=white},
legend style={legend cell align=left, align=left, draw=white!15!black},
xlabel style={font={\scriptsize}},ylabel style={font=\scriptsize},  ylabel shift={-0cm},ticklabel style={font=\scriptsize}
]
\addplot [color=black, line width = 1.5]
  table[row sep=crcr]{%
-0.7	0.100043026374266\\
-0.627272727272727	0.100027563283989\\
-0.619191919191919	0.198177627324903\\
-0.514141414141414	0.166578804816956\\
-0.401010101010101	0.131350190388944\\
-0.271717171717172	0.0898393813869566\\
-0.11010101010101	0.0366510198833068\\
-0.00505050505050508	0.00168350190095767\\
0.00303030303030305	0\\
0.0111111111111111	0.100000309289804\\
0.1	0.10000093738894\\
};

\addplot [color=black, dotted]
  table[row sep=crcr]{%
-0.7	0\\
0.1	0\\
};

\addplot [color=black, dotted]
  table[row sep=crcr]{%
-0.7	0.2\\
0.1	0.2\\
};

\end{axis}

\begin{axis}[%
width=0.420\fwidth,
height=\fheight,
at={(0.54\fwidth,0\fheight)},
scale only axis,
xmin=-0.7,
xmax=0.1,
xlabel style={font=\color{white!15!black}},
xlabel={$x_0$},
ymin=-0.01,
ymax=0.25,
ylabel style={font=\color{white!15!black}},
ylabel={$h_1$},
axis background/.style={fill=white},
legend style={legend cell align=left, align=left, draw=white!15!black},
xlabel style={font={\scriptsize}},ylabel style={font=\scriptsize},  ylabel shift={-0cm},ticklabel style={font=\scriptsize}
]
\addplot [color=black, line width = 1.5]
  table[row sep=crcr]{%
-0.7	0.0999569736257337\\
-0.627272727272727	0.0999724367160109\\
-0.619191919191919	0.00182237267509722\\
-0.514141414141414	0.0334211951830439\\
-0.401010101010101	0.0686498096110562\\
-0.271717171717172	0.110160618613043\\
-0.11010101010101	0.163348980116693\\
-0.00505050505050508	0.198316498099042\\
0.00303030303030305	0.2\\
0.0111111111111111	0.0999996907101962\\
0.1	0.0999990626110602\\
};

\addplot [color=black, dashed]
  table[row sep=crcr]{%
-0.7	0\\
0.1	0\\
};

\addplot [color=black, dashed]
  table[row sep=crcr]{%
-0.7	0.2\\
0.1	0.2\\
};

\end{axis}

\begin{axis}[%
width=1.227\fwidth,
height=1.227\fheight,
at={(-0.16\fwidth,-0.135\fheight)},
scale only axis,
xmin=0,
xmax=1,
ymin=0,
ymax=1,
axis line style={draw=none},
ticks=none,
axis x line*=bottom,
axis y line*=left,
legend style={legend cell align=left, align=left, draw=white!15!black},
xlabel style={font={\scriptsize}},ylabel style={font=\scriptsize},  ylabel shift={-0cm},ticklabel style={font=\scriptsize}
]
\end{axis}
\end{tikzpicture}%
}
	\caption{Illustration of the discontinuity of the solution map of \eqref{eq:fesd_equation} for the PSS $\dot{x} \in  2-\mathrm{sign}(x) +x^2$ for $T = 0.2$ and $\NFE =2$.}
	\label{fig:non_lipsitz}
		\vspace{-0.45cm}					
\end{figure}
\subsubsection{Extension for the case of $c_{\Nstg}\neq1$}
In the case of $c_{\Nstg} \neq 1$, we need to solve the additional LP \eqref{eq:lp_boundary_points} to obtain the boundary points. 
Note that if we have $g_i(x_{n+1}) = \min_{j\in\mathcal{J}}g_j(x_{n+1})$ for more than one $i$, the variables $\theta_{n,\Nstg+1}$ are not unique and the $\LP(x_{n+1})$ has infinitely many solutions. However, these variables are neither used in the cross complementarities nor step equilibration. 
Therefore, we can discard  $\theta_{n,\Nstg+1}$ and simplify \eqref{eq:lp_boundary_points} to:
\begin{align*}
	\lambda_{n,\Nstg+1} &= g(x_{n+1}) - \mu_{n,\Nstg+1} e,\\
	\lambda_{n,\Nstg+1} &\geq 0,
\end{align*}
which has $\Nsys+1$ unknowns and $\Nsys$ equalities and $\Nsys$ inequalities for a given $x_{n+1}$. Now suppose that the first $m_1$ components of $\lambda_{n,\Nstg+1}$ are zero (e.g., implied by cross complementarity) and the remaining $m_2$ are strictly positive, with $m_1+m_2 = \Nsys$. We have that 
\begin{subequations}\label{eq:lp_boundary_points_simpler}
	\begin{align}
		g_i(x_{n+1})  &=  \mu_{n,\Nstg+1} ,\; i = 1,\ldots,m_1, \\
		\lambda_{n,\Nstg+1,i} &= 0,\; i = 1,\ldots,m_1, \\
		\lambda_{n,\Nstg+1,j} &= g_j(x_{n+1}) - \mu_{n,\Nstg+1},\; j = \Nsys-m_2+1,\ldots,\Nsys.
	\end{align}
\end{subequations}
As the first $m_1$ relations all assign the same value to $\mu_{n,\Nstg+1}$, we can discard $m_1-1$ of them and thus we end up with a system of $m_1+m_2+1= \Nsys+1$ equations and $\Nsys+1$ unknowns. 
This system has still the important property that $\mu_{n,\Nstg+1} = \min_i g_i(x_{n+1})$. With this simplification for an RK method with $c_{\Nstg}\neq1$ we have $n_{{\theta}}= \NFE\Nstg\Nsys$, $n_{\lambda}= n_{{\theta}}+ (\NFE-1)\Nsys$. The new variables are determined by the square linear system \eqref{eq:lp_boundary_points_simpler}. 
Hence, it is straightforward to extend Theorem \ref{th:fesd_locally_unique_solution} for the case of $c_{\Nstg} \neq 1$.

\subsection{Convergence and order of FESD}
In this subsection we prove that under reasonable assumptions the sequence of approximations $\hat{x}_h(\cdot)$ generated by the FESD method converges with high order to a solution of \eqref{eq:pws1} in the sense of Definition \ref{def:piecewise_active}.
\sloppy{Recall that $h = \max_{n\in\{0,\ldots \NFE-1\}} h_n$.} The proof is inspired by the proof of Theorem 4.3 in \cite{Stewart1990b}.
We consider also $t_{\mathrm{s},0} = 0$ as a switching point, since at this time point the active set for the first interval $(t_{\mathrm{s},0},t_{\mathrm{s},1})$ is determined.

Note that for generating solution approximations with FESD it is sufficient to consider only two finite elements at a time, i.e., $\NFE =2$ in Eq.~\eqref{eq:fesd_equation}, and then to append the solutions in order to construct $\hat{x}_h(t),\ t\in[0,T]$ via Eq. \eqref{eq:continious_time_fesd}. 
This requires of course to have only one switch in the regarded time interval, which can always be achieved with a sufficiently small $h$. 
We define the set of all discretization grid points as $\mathcal{G} = \{ t_0,t_1,\dots,t_{\NFE}\}$. 
We treat the cases of crossing a discontinuity or entering a sliding mode, i.e., the case of $\I_n \subset \I_{n+1}^0$ and $\I_{n+1} \subseteq \I_{n+1}^0$.
\color{black}
\begin{theorem1}\label{th:integration_order}
	Suppose that $x(t)$ is a solution of \eqref{eq:pws1} in the sense of Definition \ref{def:piecewise_active} for $t\in [0,T]$ with $x(0) =x_0$. Suppose the following is true:
	\begin{enumerate}[(a)]
		\item the Assumptions \ref{ass:solution_existence} and \ref{ass:lcp_switching} are satisfied and
		$\I_n \subset \I_{n+1}^0$ and $\I_{n+1} \subseteq \I_{n+1}^0$ holds for all	$\ n = 0,\dots,\Nswitch$,
		\item the Assumption \ref{ass:irk_scheme}, \ref{ass:solution_existence_fesd} and \ref{ass:regularity} hold for the FESD problem \eqref{eq:fesd_compact}.
	\end{enumerate}
Then $x(\cdot)$ is a limit point of the sequence of approximations $\hat{x}_h(\cdot)$, defined in Eq. \eqref{eq:continious_time_fesd} as $h\downarrow 0$.
Moreover, for sufficiently small ${h}>0$, the solution of \eqref{eq:fesd_compact} generates a solution approximation $\hat{x}_h(t)$ on $[0,T]$ such that:
\begin{subequations}\label{eq:fesd_convergence}
	\begin{align}
   		|\tsnhat- \tsn | &= O({h}^p)\ \text{for every } n=0,\dots,\Nswitch, \label{eq:fesd_convergence_t}\\
		\| \hat{x}_h(t_n) - x(t_n) \| &= O({h}^p), \ \text{for all } t_{n} \in\mathcal{G}.  \label{eq:fesd_convergence_x}
	\end{align}
\end{subequations}
\end{theorem1}
\textit{Proof.} 
The proof will be carried out by induction, where we consider the switching events $n = 0,\dots, \Nswitch$ and the corresponding time intervals $(\tsn,\tsnn)$, with a slight abuse of notation where $t_{\mathrm{s},\Nswitch+1}=T$ is not necessarily a switching point. 
Regard $n = 0$, where we have trivially that $ {t}_{0} = 0$, thus
\begin{align*}
	| \hat{t}_{0}- {t}_{0} |  & = 0 = O({h}^p),\  \| \hat{x}_{h}(0)- x(0) \|   = 0 = O({h}^p).
\end{align*}
Moreover, $\I(x_0) = \I(\hat{x}_{h}(0)) = \I_0^0$.

Now we suppose \eqref{eq:fesd_convergence} is true for $n$, i.e., at $t = \tsn$.
We show that the same statements are true for $n+1$. 
By the induction hypothesis and due to continuity of $g_i(x), i \in \mathcal{J}$, we have that for sufficiently small $h$ the equality $\I(\hat{x}(\tsnhat)) = \I(x(\tsn)) = \I_n^0$ holds.
Moreover, by Lipschitz continuity of $f_i(x)$ and $\nabla g_i(x), i \in \mathcal{J}$, it follows that (cf. Section \ref{sec:stewarts_lcp}) 

\begin{align*}
	 M_{ \I_n^0}(\hat{x}_h(\tsnhat))  \to  M_{ \I_n^0}(x(\tsn)) \text{ as } h\downarrow 0.
\end{align*}
According to Theorem \ref{th:index_sets} the solution of the LCP \eqref{eq:lcp_switching} \sloppy{corresponding to $M_{ \I_n^0}(x(\tsn))$} determines the new index set $\I_n = \{i \in \I_n^0 \mid \tilde{\theta}_i > 0 \}$. 
Moreover, by  Assumption \ref{ass:lcp_switching} this LCP is strongly stable and due to Lemma \ref{lem:lcp_active_sets} (cf. Appendix \ref{app:switch_detection_lcp}), for sufficiently small $h>0$ the solution of the LCP corresponding $M_{ \I_n^0}(\hat{x}_h(\tsnhat))$ has a solution such that  $\hat{\I}_n=\{i \mid \tilde{\theta}_i > 0 , i \in \I_n^0\} = \{i \mid \theta_i > 0 , i \in \I_n^0\} = \I_n$. Thus, we conclude that both the solution approximation and the solution \textit{predicted} the same active set $\I_n$ in a neighborhood of $\tsnhat$ and $\tsn$, respectively.

It is left to verify that such an active set $\I_n$ predicted by the solution approximation is indeed feasible for the FESD problem. 
Note that by the induction hypothesis and the reasoning above the solution approximation and $x(\cdot)$ have the same corresponding active set in a neighborhood of $\tsn$. 
For a fixed active set, as a consequence of Proposition \ref{prop:solution_existence} the arising DAE \eqref{eq:dcs_dae} has a unique solution. Finally, under this setting with the given active sets in a neighborhood of $\tsn$, according to Theorem \ref{th:fesd_locally_unique_solution} there is a locally unique solution to a FESD problem, thus we can construct an appropriate $\hat{x}_h(\cdot)$.

Note that one can make arbitrarily many integration steps with a fixed $\I_n$ before the next switch in time is reached.
Again, due to Theorem \ref{th:fesd_locally_unique_solution} the corresponding FESD problem has a locally unique solution.

Now we provide an error estimate for the solution approximation until the next switching point. 
First, we define $\tilde{x}(t)$ to be the exact \textit{extended} solution of the DAE \eqref{eq:dcs_dae} with the fixed active set $\I_n$ on the interval $t\in [\tsn,\tilde{t}]$,  with $\tilde{x}(\tsn) = x(\tsn)$ and $\tilde{t} > \tsnn$. 
Obviously, it holds that $x(t) = \tilde{x}(t)$ for all $t \in [\tsn,\tsnn]$.
Second, from the discussions in Section \ref{sec:cross_comp} we know that active-set changes can only happen at boundaries of the finite elements, thus it holds that $\tsnhat \in \mathcal{G}$ for all $n = 0,\ldots,\hat{N}_{\mathrm{sw}}$. 
Third, by the induction hypothesis we have $\| \hat{x}_h(\tsnhat) - x(\tsnhat) \| = O({h}^p)$. 
As noted above, for a fixed active set $\I_n$ and fixed $h_n$ the FESD equations boil down to standard RK equations applied to \eqref{eq:dcs_dae}. 
Thus, from Assumption \ref{ass:irk_scheme} we have the estimate

\begin{align}\label{eq:irk_solver_accuracry}
	\| \hat{x}_h(\tsnnhat) - \tilde{x}(\tsnnhat) \|_{}  &= O({h}^p).
\end{align}

With the help of this estimate, in the next few steps we prove that $| \tsnnhat- \tsnn | = O(h^p)$. 
It is assumed that we regard crossing a discontinuity or entering a sliding mode, i.e., the case of $\I_n \subset \I(x(\tsnn))$ and $\I_{n+1} \subseteq \I(x(\tsnn))$.
We need to distinguish the two scenarios: I. $\tsnnhat > \tsnn$ and II. $\tsnnhat \leq \tsnn$.\\
\textbf{Case I.} 
Regard the following indices $j \in \I_n$ and $i\in \I(x(\tsnn)) \setminus \I_n$. 
This means that $\min_k g_k(x((\tsnn)) = g_i(x(\tsnn)) = g_j(x(\tsnn)) =\mu(\tsnn)$ holds and one can locally regard the following \textit{switching} function
\begin{align*}
	\psi_{i,j}(x(t)) = g_i(x(t)) - g_j(x(t)) = \lambda_i(t)-\lambda_j(t).
\end{align*}
Note that this function is Lipschitz continuous.
It must by definition become zero when an active-set change happens. 

Due to the strict complementarity assumed in Assumption \ref{ass:lcp_switching} (see also part 9 of the proof of \cite[Theorem 4.3]{Stewart1990b}, and the remarks after Assumption \ref{ass:lcp_switching}) we have at $t_{\mathrm{s},n+1}^-$ that $\dot{\lambda}_j(t_{\mathrm{s},n+1}^-)=0$ and $\dot{\lambda}_i(t_{\mathrm{s},n+1}^-)<0$. 
Therefore, it holds that:
\begin{align}\label{eq:derivative_of_psi}
	\psi_{i,j}(x(\tsnn)) = 0,\ \frac{\dd}{\dd t} \psi_{i,j}(x(t_{\mathrm{s},n+1}^-)) < 0.
\end{align}
Obviously, the same assertion holds for $\tilde{x}(t)$. 
Moreover, due to the smoothness of $\tilde{x}(t)$, we have $\psi_{i,j}(\tilde{x}(t)) < 0$ for $t\in(\tsnn,\tsnn+\epsilon)$ for some $\epsilon>0$.
\color{black}

\begin{figure}[t]
	\centering
	{ \setlength{\fwidth}{9.5cm}
\setlength{\fheight}{3.2cm}
\definecolor{mycolor1}{rgb}{0.00000,0.44700,0.74100}%
\definecolor{mycolor2}{rgb}{0.49400,0.18400,0.55600}%
\definecolor{mycolor3}{rgb}{0.85000,0.32500,0.09800}%
\definecolor{mycolor4}{rgb}{0.92900,0.69400,0.12500}%
\begin{tikzpicture}

\begin{axis}[%
width=0.411\fwidth,
height=\fheight,
at={(0\fwidth,0\fheight)},
scale only axis,
xmin=0.9,
xmax=1.2,
xtick={1,1.1},
xticklabels={{$t_{\mathrm{s},n+1}$},{$\hat{t}_{\mathrm{s},n+1}$}},
xlabel style={font=\color{white!15!black}},
xlabel={$t$},
ymin=-0.2,
ymax=0.2,
ytick={\empty},
ylabel style={font=\color{white!15!black}},
ylabel={$\psi_{i,j}(\cdot)$},
axis background/.style={fill=white},
legend style={at={(0.358,0.65)}, anchor=south west, legend cell align=left, align=left, draw=white!15!black,nodes={scale=0.70, transform shape}},
xlabel style={font={\scriptsize}},ylabel style={font=\scriptsize},  ylabel shift={-0cm},ticklabel style={font=\scriptsize}
]
\addplot [color=mycolor1]
  table[row sep=crcr]{%
0.895	0.0820934232844315\\
0.995	0.00392695044435576\\
1	0\\
};
\addlegendentry{$\psi_{i,j}(x(t))$}

\addplot [color=mycolor2]
  table[row sep=crcr]{%
1	0\\
1.125	-0.0975451610080642\\
1.19	-0.147020162616152\\
1.205	-0.158238483590793\\
};
\addlegendentry{$\psi_{i,j}(\tilde{x}(t))$}

\addplot [color=mycolor3]
  table[row sep=crcr]{%
0.895	0.144288040601379\\
0.965	0.0957938454970677\\
1	0.0711574191366426\\
1.1	0\\
};
\addlegendentry{$\psi_{i,j}(\hat{x}_h(t))$}

\addplot [color=mycolor4, line width=1.2pt]
  table[row sep=crcr]{%
0.895	0.0519656256987491\\
1	0\\
1.205	-0.101456697792796\\
};
\addlegendentry{$\psi_{i,j}({x}(\tsnn))-a_{I}(t-\tsnn)$}

\addplot [color=mycolor2, mark size=1.5pt, mark=*, mark options={solid, fill=mycolor2, mycolor2}, forget plot]
  table[row sep=crcr]{%
1	0\\
};
\addplot [color=mycolor3, dashed, forget plot]
  table[row sep=crcr]{%
1.105	-0.00356996131923459\\
1.205	-0.0746892296249124\\
};
\addplot [color=black, forget plot]
  table[row sep=crcr]{%
0.895	0\\
1	0\\
1.205	0\\
};
\addplot [color=black, dotted, forget plot]
  table[row sep=crcr]{%
1	-0.24\\
1	0.24\\
};
\addplot [color=black, dotted, forget plot]
  table[row sep=crcr]{%
1.1	-0.24\\
1.1	0.24\\
};
\end{axis}

\begin{axis}[%
width=0.411\fwidth,
height=\fheight,
at={(0.54\fwidth,0\fheight)},
scale only axis,
xmin=0.9,
xmax=1.2,
xtick={1,1.1},
xticklabels={{$\hat{t}_{\mathrm{s},n+1}$},{$t_{\mathrm{s},n+1}$}},
xlabel style={font=\color{white!15!black}},
xlabel={$t$},
ymin=-0.2,
ymax=0.2,
ytick={\empty},
axis background/.style={fill=white},
legend style={at={(0.378,0.65)}, anchor=south west, legend cell align=left, align=left, draw=white!15!black,nodes={scale=0.70, transform shape}},
xlabel style={font={\scriptsize}},ylabel style={font=\scriptsize},  ylabel shift={-0cm},ticklabel style={font=\scriptsize}
]
\addplot [color=mycolor1]
  table[row sep=crcr]{%
0.895	0.144288040601379\\
0.965	0.0957938454970677\\
1.055	0.0321078169300275\\
1.1	0\\
};
\addlegendentry{$\psi_{i,j}(x(t))$}

\addplot [color=mycolor3]
  table[row sep=crcr]{%
0.895	0.0820934232844315\\
0.995	0.00392695044435554\\
1	0\\
};
\addlegendentry{$\psi_{i,j}(\hat{x}_h(t))$}

\addplot [color=mycolor4, line width=1.2pt]
  table[row sep=crcr]{%
0.895	0.123011057017888\\
1.205	-0.030080635774359\\
};
\addlegendentry{$\psi_{i,j}({x}(\tsnnhat))-a_{II}(t-\tsnnhat)$}

\addplot [color=mycolor1, mark size=1.5pt, mark=*, mark options={solid, fill=mycolor1, mycolor1}, forget plot]
  table[row sep=crcr]{%
1	0.0711574191366426\\
};
\addplot [color=black, forget plot]
  table[row sep=crcr]{%
0.895	0\\
1.205	0\\
};
\addplot [color=mycolor1, dashed, forget plot]
  table[row sep=crcr]{%
1.105	-0.00356996131923459\\
1.205	-0.0746892296249124\\
};
\addplot [color=mycolor3, dashed, forget plot]
  table[row sep=crcr]{%
1.005	-0.00392695044435576\\
1.125	-0.0975451610080642\\
1.19	-0.147020162616152\\
1.205	-0.158238483590793\\
};
\addplot [color=black, dotted, forget plot]
  table[row sep=crcr]{%
1.1	-0.24\\
1.1	0.24\\
};
\addplot [color=black, dotted, forget plot]
  table[row sep=crcr]{%
1	-0.24\\
1	0.24\\
};
\end{axis}

\begin{axis}[%
width=1.227\fwidth,
height=1.227\fheight,
at={(-0.16\fwidth,-0.135\fheight)},
scale only axis,
xmin=0,
xmax=1,
ymin=0,
ymax=1,
axis line style={draw=none},
ticks=none,
axis x line*=bottom,
axis y line*=left,
legend style={legend cell align=left, align=left, draw=white!15!black},
xlabel style={font={\scriptsize}},ylabel style={font=\scriptsize},  ylabel shift={-0cm},ticklabel style={font=\scriptsize}
]
\end{axis}
\end{tikzpicture}%
}
	\caption{The left plots shows an illustration of the argument of Case I: $\tsnnhat > \tsnn$ and the right plot shows an illustration of the argument of Case II: $\tsnnhat \leq t_{n+1}$, for establishing $|\tsnnhat -\tsnn | = O(h^p)$.} 
	\label{fig:switching_arg}
	\vspace{-0.45cm}					
\end{figure}
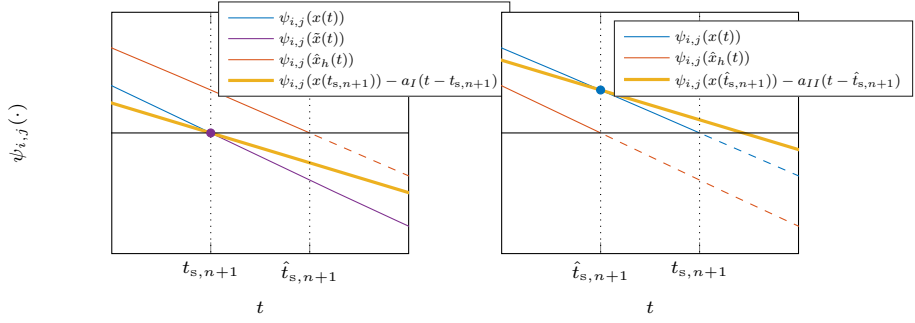

Similarly, for the solution approximation we have $\psi_{i,j}(\hat{x}_h(t)) = g_i(\hat{x}_h( t)) - g_j(\hat{x}_h(t))$.
Since  $\tsnnhat > \tsnn$, due to continuity of $\psi_{i,j}(\cdot)$ and $\hat{x}_h(\cdot)$ it follows that 
\begin{align*}
	 \psi_{i,j}(\hat{x}_h(\tsnn)) >0
\text{ and }
	\frac{\dd}{\dd t} \psi_{i,j}(\hat{x}_h( \tsnn)) < 0.
\end{align*}
Now from Lipschitz continuity of $\psi_{i,j}(\cdot)$ and \eqref{eq:irk_solver_accuracry} we can establish that
\begin{align}\label{eq:first_estimate_psi}
	\begin{split}
		| \underbrace{\psi_{i,j}(\hat{x}_h( \tsnnhat))}_{=0}-\underbrace{\psi_{i,j}(\tilde{x}(\tsnnhat))}_{<0} |  &\leq L_{\psi} \| 	\hat{x}_h(\tsnnhat) - \tilde{x}(\tsnnhat) \|_{}, \\
		|\psi_{i,j}(\tilde{x}(\tsnnhat)) | &\stackrel{}{=} O(h^p).
	\end{split}
\end{align}
Note that in contrast to $\psi_{i,j}(x(t))$, the function $\psi_{i,j}(\tilde{x}(t))$ is smooth in a neighborhood of $\tsnn$.
Thus, we look at the first-order Taylor approximation of $\psi_{i,j}(\cdot)$ at $\tilde{x}(\tsnn)$.
\begin{align*}
	\psi_{i,j}(\tilde{x}(t)) &= \psi_{i,j}(\tilde{x}(\tsnn))+\frac{\dd}{\dd t}\psi_{i,j}(\tilde{x}(\tsnn)) (t-\tsnn)  + o(|t-\tsnn|). 
\end{align*}
Due to continuity, $\psi_{i,j}(\tilde{x}(t))$ is decreasing for $t\in [\tsnn,\tsnnhat]$\color{black}, there exists some positive constant $a_{I}$ with $0< a_{I} <  | \frac{\dd}{\dd t} \psi_{i,j}(\tilde{x}(\tsnn)| $ such that for sufficiently small $h$ and $t\in [\tsnn,\tsnnhat]$ it holds that
\begin{align*}
	\psi_{i,j}(\tilde{x}(t)) \leq  \psi_{i,j}(\tilde{x}(\tsnn))-a_{I}(t-\tsnn).
\end{align*}
The arguments above are illustrated in the left plot of Figure \ref{fig:switching_arg}. 
From the last inequality and \eqref{eq:first_estimate_psi} at  $t = \tsnnhat$ we have that
$\psi_{i,j}(\tilde{x}(\tsnnhat)) <0$, 
$\psi_{i,j}(\tilde{x}(\tsnn)) = 0$,
and conclude that
\begin{align}\label{eq:time_bound1}
	\tsnnhat -\tsnn \leq O(h^p).
\end{align}
This completes the consideration of case I.\\

\textbf{Case II.} 
We apply similar arguments as for I. 
Under the assumption of $\tsnnhat < \tsnn$, following similar lines as in the proof of \cite[Theorem 4.3]{Stewart1990b}, we first prove $\tsnnhat \to \tsnn$ and establish subsequently the convergence rate.

Regard the set $H = \{ h >0 \mid \tsnnhat < \tsnn\}$. 
By the assumption of case II and the induction hypothesis, it holds that $\tsnnhat \in [\tsnhat,\tsnn]$. 
Since this is a bounded set, there must be a subsequence $H'\subseteq H$ with $h\downarrow0$ such that $\tsnnhat \to \bar{t}$. 
We show now that $\bar{t} = \tsnnhat$. 
We consider again the function $\psi_{i,j}(\cdot)$ for some $j \in \I_n$ and $i\in \I(x(\tsnn)) \setminus \I_n$
\begin{align*}
	\psi_{i,j}(x(t)) = g_i(x(t)) - g_j(x(t)) = \lambda_i-\lambda_j,
\end{align*}
which becomes zero at an active-set change and is positive before.  
Similarly, active-set changes for the solution approximation happen when ${\psi}_{i,j}(\hat{x}_h(\tsnnhat)) = 0$. 
We remind the reader that earlier it was shown that $\hat{\I}_n = \I_n$.

By taking $h \downarrow 0$ and $h \in H'$ from \eqref{eq:irk_solver_accuracry} it follows that ${\psi}_{i,j}({x}(\bar{t}))  = 0$.  
By the definition of a switching point, there must be a $i \notin \I_n$, but $i \in \I(x(\bar{t}))$. 
However, this contradicts the assumption that $\I(x(t)) = \I_n$ for $t \in (\tsn,\tsnn)$ and we conclude that $\bar{t} \notin (\tsn,\tsnn)$.

On the other hand at $\tsn$, {due to strict complementarity in the active-set determining LCP (cf. Theorem\ref{th:index_sets} and Assumption \ref{ass:lcp_switching}}), if some $i \in \I(x(\tsn))\setminus \I_n$ and $j\in \I_n$, we know that 
\begin{align*}
\frac{\dd}{\dd t}\psi_{i,j}(x(t_{\mathrm{s,n}}^+)) =  g_i(x(t_{\mathrm{s,n}}^+))-g_j(x(t_{\mathrm{s,n}}^+)) > 0.
\end{align*}
\color{black}
Due to continuity of the functions $g_i(\cdot), i\in \mathcal{J}$, and the induction hypothesis, there exists some $\epsilon>0$ such that
\begin{align*}
 \frac{\dd}{\dd t}(g_i(\hat{x}_h(t))-g_j(\hat{x}_h(t)) > 0 \ \text{for}\ t \in [\tsnhat,\tsnhat+\epsilon].
\end{align*}

However, when a switch happens the derivative in the last line must be negative at $t$ (cf. Remark \ref{rem:lambda_at_switch}), thus $\tsnnhat > \tsnhat+\epsilon$, i.e., with $h\downarrow0,  h\in H'$, $\bar{t} > \tsn+\epsilon$. 
This means that $\bar{t} \neq \tsn$ and the only option that is left is $\tsnnhat \to \bar{t} = \tsnnhat$ as $h\downarrow0, h\in H'$.
\color{black}

Now we continue with establishing the convergence rate for $\tsnnhat \to \tsnn$.  
From $\tsnnhat \leq \tsnn$ we have from the definition of $\psi_{i,j}(\cdot)$ that $\psi_{i,j}(x(\tsnnhat))>0$ and $\psi_{i,j}(\hat{x}_h(\tsnnhat))=0$. 
Using the fact that $\tilde{x}(\tsnnhat) = {x}(\tsnnhat)$  and \eqref{eq:irk_solver_accuracry} we have
\begin{align*}
|\psi_{i,j}(x(\tsnnhat))|&= |\psi_{i,j}(x(\tsnnhat))-\psi_{i,j}(\hat{x}_h(\tsnnhat))|\\
& \leq  L_{\psi} \|\tilde{x}(\tsnnhat)-\hat{x}_h(\tsnnhat)\|
 = O(h^p).
\end{align*}
We again use a first-order expansion:
\begin{align*}
	\psi_{i,j}({x}(t)) &= \psi_{i,j}({x}(\tsnnhat))+\frac{\dd}{\dd t}\psi_{i,j}({x}(\tsnnhat)) (t-\tsnnhat)  + o(|t-\tsnnhat|). 
\end{align*}

Once again, due to assumption \ref{ass:lcp_switching}, we have that $\frac{\dd }{\dd t}\psi_{i,j}({x}(\tsnn)) <0$.
Note that $ \frac{\dd }{\dd t}\psi_{i,j}({x}(t_{\mathrm{s},n+1}^-)) <0$. 
From $\tsnnhat \to \tsnn$ and continuity of $\psi_{i,j}(\cdot)$ it follows that for sufficiently small $h>0$: 
\begin{align*}
\frac{\dd}{\dd t}\psi_{i,j}({x}(\tsnnhat))<0. 
\end{align*}
\color{black}
Using similar reasoning as in case I (see right plot of Figure \ref{fig:switching_arg} for an illustration of the argument), there exists some $a_{II}>0$ such that from the last equation at $t = \tsnn$ it follows $0 \leq O(h^p)  - a_{II}(\tsnn-\tsnnhat)$, i.e., 
\begin{align}\label{eq:time_bound2}
		 \tsnn - \tsnnhat  &\leq O({h}^p).
\end{align}

Putting \eqref{eq:time_bound1} and \eqref{eq:time_bound2} together, we obtain the first part of the induction statement, i.e., 
\begin{align}\label{eq:time_error_estimate}
		| \tsnn - \tsnnhat| &= O({h}^p).
\end{align}

To complete the induction step we must prove that \eqref{eq:fesd_convergence_x} holds for $ t=\tsnnhat$. 
For $ \tsnnhat \leq \tsnn$ we have $\tilde{x}(\tsnnhat) = x(\tsnnhat)$ and the assertion follows directly from \eqref{eq:irk_solver_accuracry}. 
Note that for any other $t_n \in \mathcal{G}$ that is not a switching point \eqref{eq:fesd_convergence_x} follows immediately from Assumption \ref{ass:irk_scheme}. 

It is left to investigate the case of $ \tsnnhat > \tsnn$. 
Using Lipschitz continuity of $x(\cdot)$, $\tilde{x}(\cdot)$, the fact that $\tilde{x}(\tsnn) = x(\tsnn)$ and \eqref{eq:time_error_estimate} one obtains
\begin{align*}
	\begin{split}
		&\|\hat{x}_h(\tsnnhat)-{x}(\tsnnhat)\| \leq\|\hat{x}_h(\tsnnhat)-\tilde{x}(\tsnnhat)\|+\|{x}(\tsnn)-{x}(\tsnnhat)\| \\ 
		& 
		+ \|\tilde{x}(\tsnnhat)-\tilde{x}(\tsnn)\|
		\leq O(h^p) + (L_x+ L_{\tilde{x}})|\tsnn-\tsnnhat| = O(h^p).
	\end{split}
\end{align*}
Moreover, from Lipschitz continuity of $g_i(\cdot), i \in \mathcal{J}$ and the last inequality for sufficiently small $h>0$ we have that $\I(x(\tsnn)) = \I(\hat{x}(\tsnnhat))$, which completes the induction step for $n+1$. 
With the use of an interpolation scheme, from \eqref{eq:fesd_convergence} it follows that we can make a continuous-time approximation  $\hat{x}_h(t)$ for $t\in[0,T]$ with the accuracy $O(h^{\bar{p}})$, $1\leq\bar{p}\leq p$ for $t\notin \mathcal{G}$. 
Therefore, it follows that the sequence of approximations $\hat{x}_h(t)$ generated by the FESD method converges to a solution $x(t)$ in the sense of Definition \ref{def:piecewise_active} for $t\in [0, T]$. 
The proof is completed. \qed
In the next subsection, we illustrate the results of this theorem for several RK schemes. We compare the results obtained via FESD to the ones obtained with the standard RK discretization from Section \ref{sec:dcs_irk}.

\subsection{Illustrating the integration order }\label{sec:integrator_order}
Consider the nonsmooth IVP
\begin{align}\label{eq:ode_oscilator}
	\dot{x}(t) &= \begin{cases}
		A_1 x,  \ &c(x) < 0,\\
		A_2 x,  \ &c(x) >0,
	\end{cases} & 
\end{align}
with $ A_1 = \begin{bmatrix}
	&1 & \omega \\ &-\omega &1
\end{bmatrix}, 
\ A_2 = \begin{bmatrix}
	&1 & -\omega \\ &\omega &1
\end{bmatrix}, \
c(x)= x_1^2 + x_2^2 -1,\ \omega  =2 \pi \ \text{and } 
x(0) = (e^{-1},0)$ for $t\in [0,T]$.
The example trajectory is given in Figure \ref{fig:oscilator_solution}. 
Following Section \ref{sec:stewarts_dcs}, we can write \eqref{eq:ode_oscilator} in the form of the DCS \eqref{eq:dcs_1}, where:
\begin{align*}
	F(x) &= \begin{bmatrix}
		A_1x & A_2x
	\end{bmatrix},\; g(x) = [1 \; -\!1]^\top c(x),
\end{align*} 
where $g(x)$ was obtained and via Eq. \eqref{eq:indicator_func_formula}. 
It can be shown that the switch happens at $\ts = 1$ and that  $x(T) = (e^{(T-1)}\cos(2\pi(T-1)),-e^{(T-1)}\sin(2\pi(T-1)))$ for $T>\ts$. 
Hence, we can determine the global integration error $E(T) = \| x(T)-\hat{x}_h(T)\|$. 
\begin{figure}[t]
	\centering
	{ \input{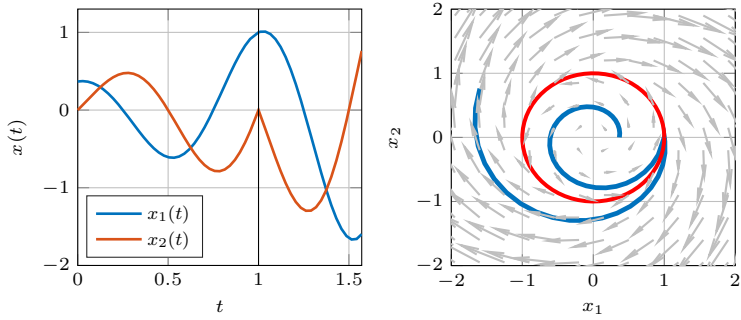}}
	\caption{Illustration of the solution to the nonsmooth IVP given by \eqref{eq:ode_oscilator}.}
	\label{fig:oscilator_solution}
\end{figure}
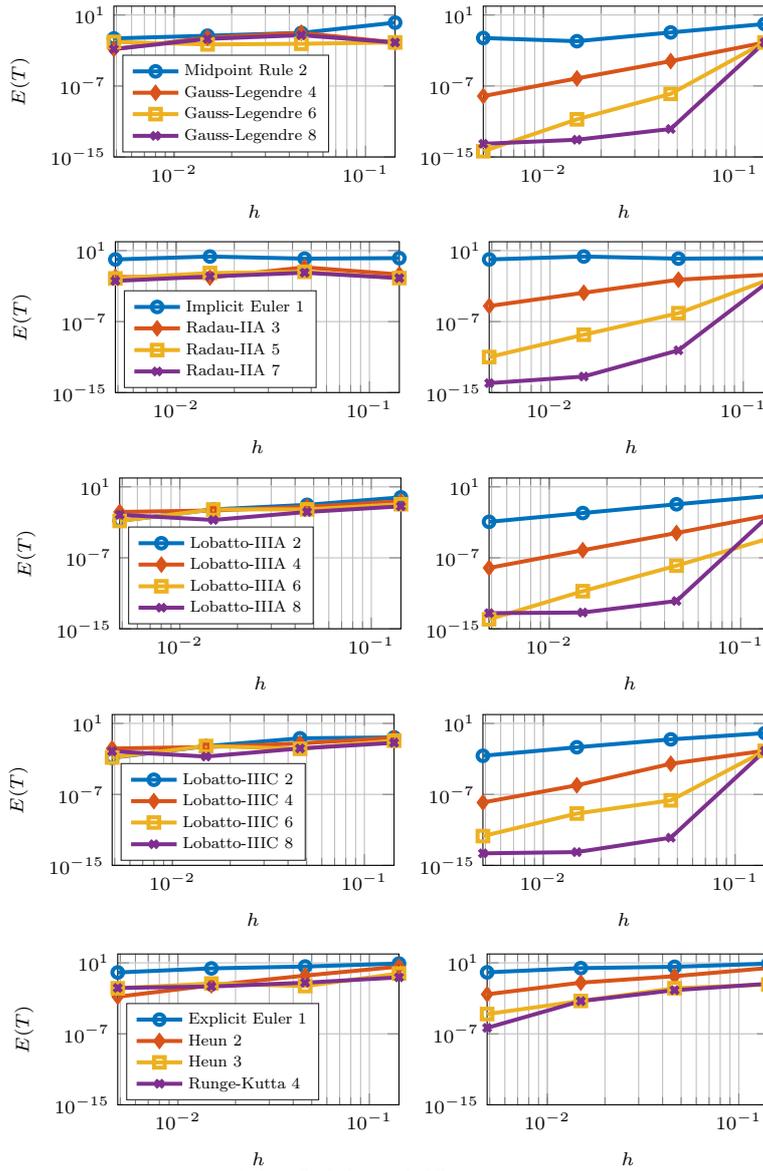
\begin{figure}[t]
	\centering
	{ \setlength{\fwidth}{9.0cm}
\setlength{\fheight}{2.0cm}
\definecolor{mycolor1}{rgb}{0.00000,0.44700,0.74100}%
\definecolor{mycolor2}{rgb}{0.85000,0.32500,0.09800}%
\definecolor{mycolor3}{rgb}{0.92900,0.69400,0.12500}%
\definecolor{mycolor4}{rgb}{0.49400,0.18400,0.55600}%
\begin{tikzpicture}

\begin{axis}[%
width=0.411\fwidth,
height=\fheight,
at={(0\fwidth,0\fheight)},
scale only axis,
xmode=log,
xmin=0.00484813681109536,
xmax=0.142799666072263,
xminorticks=true,
xlabel style={font=\color{white!15!black}},
xlabel={$h$},
ymode=log,
ymin=1e-15,
ymax=100,
yminorticks=true,
ylabel style={font=\color{white!15!black}},
ylabel={$E(T)$},
axis background/.style={fill=white},
xmajorgrids,
xminorgrids,
ymajorgrids,
yminorgrids,
legend style={at={(0.03,0.03)}, anchor=south west, legend cell align=left, align=left, draw=white!15!black,nodes={scale=0.75, transform shape}},
xlabel style={font={\scriptsize}},ylabel style={font=\scriptsize},  ylabel shift={-0cm},ticklabel style={font=\scriptsize}
]
\addplot [color=mycolor1, line width=1.5pt, mark=o, mark options={solid, mycolor1}]
  table[row sep=crcr]{%
0.142799666072263	1.39049559988208\\
0.0461998919645558	0.0966447269905166\\
0.0149599650170943	0.0467060468485782\\
0.00484813681109536	0.0229559479619011\\
};
\addlegendentry{Midpoint Rule 2}

\addplot [color=mycolor2, line width=1.5pt, mark=diamond, mark options={solid, mycolor2}]
  table[row sep=crcr]{%
0.142799666072263	0.00792854809898558\\
0.0461998919645558	0.0986931339922765\\
0.0149599650170943	0.028428094108133\\
0.00484813681109536	0.00144064024985568\\
};
\addlegendentry{Gauss-Legendre 4}

\addplot [color=mycolor3, line width=1.5pt, mark=square, mark options={solid, mycolor3}]
  table[row sep=crcr]{%
0.142799666072263	0.00806644399839152\\
0.0461998919645558	0.00573454518045635\\
0.0149599650170943	0.0048194123736941\\
0.00484813681109536	0.0093857333704811\\
};
\addlegendentry{Gauss-Legendre 6}

\addplot [color=mycolor4, line width=1.5pt, mark=x, mark options={solid, mycolor4}]
  table[row sep=crcr]{%
0.142799666072263	0.00806686750400564\\
0.0461998919645558	0.0534503146711053\\
0.0149599650170943	0.0204709023101535\\
0.00484813681109536	0.00144067335223297\\
};
\addlegendentry{Gauss-Legendre 8}

\end{axis}

\begin{axis}[%
width=0.411\fwidth,
height=\fheight,
at={(0.54\fwidth,0\fheight)},
scale only axis,
xmode=log,
xmin=0.00484813681109536,
xmax=0.142799666072263,
xminorticks=true,
xlabel style={font=\color{white!15!black}},
xlabel={$h$},
ymode=log,
ymin=1e-15,
ymax=100,
yminorticks=true,
axis background/.style={fill=white},
xmajorgrids,
xminorgrids,
ymajorgrids,
yminorgrids,
legend style={legend cell align=left, align=left, draw=white!15!black},
xlabel style={font={\scriptsize}},ylabel style={font=\scriptsize},  ylabel shift={-0cm},ticklabel style={font=\scriptsize}
]
\addplot [color=mycolor1, line width=1.5pt, mark=o, mark options={solid, mycolor1}]
  table[row sep=crcr]{%
0.142799666072263	0.967160314552169\\
0.0461998919645558	0.108621305175554\\
0.0149599650170943	0.0111860002772262\\
0.00484813681109536	0.0256762117929524\\
};

\addplot [color=mycolor2, line width=1.5pt, mark=diamond, mark options={solid, mycolor2}]
  table[row sep=crcr]{%
0.142799666072263	0.00792854809898489\\
0.0461998919645558	6.30285066358647e-05\\
0.0149599650170943	7.18099964935881e-07\\
0.00484813681109535	7.30849991636972e-09\\
};

\addplot [color=mycolor3, line width=1.5pt, mark=square, mark options={solid, mycolor3}]
  table[row sep=crcr]{%
0.142799666072263	0.00806644399838584\\
0.0461998919645558	1.30896526950863e-08\\
0.0149599650170943	1.73882019893767e-11\\
0.00484813681109535	4.44089209850063e-15\\
};

\addplot [color=mycolor4, line width=1.5pt, mark=x, mark options={solid, mycolor4}]
  table[row sep=crcr]{%
0.142799666072263	0.00806686750400697\\
0.0461998919645558	1.43973721833391e-12\\
0.0149599650170943	8.781864124785e-14\\
0.00484813681109535	3.03090885722668e-14\\
};

\end{axis}

\begin{axis}[%
width=1.227\fwidth,
height=1.227\fheight,
at={(-0.16\fwidth,-0.135\fheight)},
scale only axis,
xmin=0,
xmax=1,
ymin=0,
ymax=1,
axis line style={draw=none},
ticks=none,
axis x line*=bottom,
axis y line*=left,
legend style={legend cell align=left, align=left, draw=white!15!black},
xlabel style={font={\scriptsize}},ylabel style={font=\scriptsize},  ylabel shift={-0cm},ticklabel style={font=\scriptsize}
]
\end{axis}
\end{tikzpicture}%
}
	\vspace{0.01cm}
	\hspace{-0.25cm}
	\centering
	{ \setlength{\fwidth}{9.0cm}
\setlength{\fheight}{2.0cm}
\definecolor{mycolor1}{rgb}{0.00000,0.44700,0.74100}%
\definecolor{mycolor2}{rgb}{0.85000,0.32500,0.09800}%
\definecolor{mycolor3}{rgb}{0.92900,0.69400,0.12500}%
\definecolor{mycolor4}{rgb}{0.49400,0.18400,0.55600}%
\begin{tikzpicture}

\begin{axis}[%
width=0.415\fwidth,
height=\fheight,
at={(0\fwidth,0\fheight)},
scale only axis,
xmode=log,
xmin=0.00484813681109536,
xmax=0.142799666072263,
xminorticks=true,
xlabel style={font=\color{white!15!black}},
xlabel={$h$},
ymode=log,
ymin=1e-15,
ymax=100,
yminorticks=true,
ylabel style={font=\color{white!15!black}},
ylabel={$E(T)$},
axis background/.style={fill=white},
xmajorgrids,
xminorgrids,
ymajorgrids,
yminorgrids,
legend style={at={(0.03,0.03)}, anchor=south west, legend cell align=left, align=left, draw=white!15!black,nodes={scale=0.75, transform shape}},
xlabel style={font={\scriptsize}},ylabel style={font=\scriptsize},  ylabel shift={-0cm},ticklabel style={font=\scriptsize}
]
\addplot [color=mycolor1, line width=1.5pt, mark=o, mark options={solid, mycolor1}]
  table[row sep=crcr]{%
0.142799666072263	1.43146191975967\\
0.0461998919645558	1.2493863871042\\
0.0149599650170943	2.15946662357718\\
0.00484813681109536	1.01695406988224\\
};
\addlegendentry{Implicit Euler 1}

\addplot [color=mycolor2, line width=1.5pt, mark=diamond, mark options={solid, mycolor2}]
  table[row sep=crcr]{%
0.142799666072263	0.0206779440033487\\
0.0461998919645558	0.131258103035379\\
0.0149599650170943	0.00921452973781733\\
0.00484813681109536	0.0107635048397164\\
};
\addlegendentry{Radau-IIA 3}

\addplot [color=mycolor3, line width=1.5pt, mark=square, mark options={solid, mycolor3}]
  table[row sep=crcr]{%
0.142799666072263	0.00804452481618323\\
0.0461998919645558	0.0404277348513472\\
0.0149599650170943	0.0300130162948864\\
0.00484813681109536	0.00744475798356914\\
};
\addlegendentry{Radau-IIA 5}

\addplot [color=mycolor4, line width=1.5pt, mark=x, mark options={solid, mycolor4}]
  table[row sep=crcr]{%
0.142799666072263	0.00806684866774976\\
0.0461998919645558	0.0321175675981631\\
0.0149599650170943	0.0122001212865104\\
0.00484813681109536	0.00385003832678388\\
};
\addlegendentry{Radau-IIA 7}

\end{axis}

\begin{axis}[%
width=0.415\fwidth,
height=\fheight,
at={(0.546\fwidth,0\fheight)},
scale only axis,
xmode=log,
xmin=0.00484813681109536,
xmax=0.142799666072263,
xminorticks=true,
xlabel style={font=\color{white!15!black}},
xlabel={$h$},
ymode=log,
ymin=1e-15,
ymax=100,
yminorticks=true,
axis background/.style={fill=white},
xmajorgrids,
xminorgrids,
ymajorgrids,
yminorgrids,
legend style={legend cell align=left, align=left, draw=white!15!black},
xlabel style={font={\scriptsize}},ylabel style={font=\scriptsize},  ylabel shift={-0cm},ticklabel style={font=\scriptsize}
]
\addplot [color=mycolor1, line width=1.5pt, mark=o, mark options={solid, mycolor1}]
  table[row sep=crcr]{%
0.142799666072263	1.43146191975967\\
0.0461998919645558	1.2493863871042\\
0.0149599650170943	2.15946662357718\\
0.00484813681109536	1.01695406988224\\
};

\addplot [color=mycolor2, line width=1.5pt, mark=diamond, mark options={solid, mycolor2}]
  table[row sep=crcr]{%
0.142799666072263	0.0206779440033484\\
0.0461998919645558	0.00527284380871606\\
0.0149599650170943	0.000177486615084654\\
0.00484813681109535	5.70928806997717e-06\\
};

\addplot [color=mycolor3, line width=1.5pt, mark=square, mark options={solid, mycolor3}]
  table[row sep=crcr]{%
0.142799666072263	0.00804452481589434\\
0.0461998919645558	8.80576789374921e-07\\
0.0149599650170943	3.3116764885932e-09\\
0.00484813681109535	1.00424113469444e-11\\
};

\addplot [color=mycolor4, line width=1.5pt, mark=x, mark options={solid, mycolor4}]
  table[row sep=crcr]{%
0.142799666072263	0.00806684866775154\\
0.0461998919645558	5.8423599291757e-11\\
0.0149599650170943	6.45039577307217e-14\\
0.00484813681109535	1.15463194561017e-14\\
};

\end{axis}

\begin{axis}[%
width=1.24\fwidth,
height=1.24\fheight,
at={(-0.161\fwidth,-0.147\fheight)},
scale only axis,
xmin=0,
xmax=1,
ymin=0,
ymax=1,
axis line style={draw=none},
ticks=none,
axis x line*=bottom,
axis y line*=left,
legend style={legend cell align=left, align=left, draw=white!15!black},
xlabel style={font={\scriptsize}},ylabel style={font=\scriptsize},  ylabel shift={-0cm},ticklabel style={font=\scriptsize}
]
\end{axis}
\end{tikzpicture}%
}
	\vspace{0.01cm}
	\centering
	\hspace{-0.25cm}
	{ \setlength{\fwidth}{9.0cm}
\setlength{\fheight}{2.0cm}
\definecolor{mycolor1}{rgb}{0.00000,0.44700,0.74100}%
\definecolor{mycolor2}{rgb}{0.85000,0.32500,0.09800}%
\definecolor{mycolor3}{rgb}{0.92900,0.69400,0.12500}%
\definecolor{mycolor4}{rgb}{0.49400,0.18400,0.55600}%
\begin{tikzpicture}

\begin{axis}[%
width=0.411\fwidth,
height=\fheight,
at={(0\fwidth,0\fheight)},
scale only axis,
xmode=log,
xmin=0.00484813681109536,
xmax=0.142799666072263,
xminorticks=true,
xlabel style={font=\color{white!15!black}},
xlabel={$h$},
ymode=log,
ymin=1e-15,
ymax=100,
yminorticks=true,
ylabel style={font=\color{white!15!black}},
ylabel={$E(T)$},
axis background/.style={fill=white},
xmajorgrids,
xminorgrids,
ymajorgrids,
yminorgrids,
legend style={at={(0.03,0.03)}, anchor=south west, legend cell align=left, align=left, draw=white!15!black,nodes={scale=0.75, transform shape}},
xlabel style={font={\scriptsize}},ylabel style={font=\scriptsize},  ylabel shift={-0cm},ticklabel style={font=\scriptsize}
]
\addplot [color=mycolor1, line width=1.5pt, mark=o, mark options={solid, mycolor1}]
  table[row sep=crcr]{%
0.142799666072263	0.656413870649981\\
0.0461998919645558	0.0940894718846521\\
0.0149599650170943	0.0287489790115592\\
0.00484813681109536	0.00147812270333947\\
};
\addlegendentry{Lobatto-IIIA 2}

\addplot [color=mycolor2, line width=1.5pt, mark=diamond, mark options={solid, mycolor2}]
  table[row sep=crcr]{%
0.142799666072263	0.234601528567681\\
0.0461998919645558	0.057620853773693\\
0.0149599650170943	0.0216566913679285\\
0.00484813681109536	0.0148214312043352\\
};
\addlegendentry{Lobatto-IIIA 4}

\addplot [color=mycolor3, line width=1.5pt, mark=square, mark options={solid, mycolor3}]
  table[row sep=crcr]{%
0.142799666072263	0.111109946801133\\
0.0461998919645558	0.0289372267269207\\
0.0149599650170943	0.028428905846237\\
0.00484813681109536	0.00144066637826212\\
};
\addlegendentry{Lobatto-IIIA 6}

\addplot [color=mycolor4, line width=1.5pt, mark=x, mark options={solid, mycolor4}]
  table[row sep=crcr]{%
0.142799666072263	0.0634227270332261\\
0.0461998919645558	0.0148514129583631\\
0.0149599650170943	0.00185146018637372\\
0.00484813681109536	0.00721828688724957\\
};
\addlegendentry{Lobatto-IIIA 8}

\end{axis}

\begin{axis}[%
width=0.411\fwidth,
height=\fheight,
at={(0.54\fwidth,0\fheight)},
scale only axis,
xmode=log,
xmin=0.00484813681109536,
xmax=0.142799666072263,
xminorticks=true,
xlabel style={font=\color{white!15!black}},
xlabel={$h$},
ymode=log,
ymin=1e-15,
ymax=100,
yminorticks=true,
axis background/.style={fill=white},
xmajorgrids,
xminorgrids,
ymajorgrids,
yminorgrids,
legend style={legend cell align=left, align=left, draw=white!15!black},
xlabel style={font={\scriptsize}},ylabel style={font=\scriptsize},  ylabel shift={-0cm},ticklabel style={font=\scriptsize}
]
\addplot [color=mycolor1, line width=1.5pt, mark=o, mark options={solid, mycolor1}]
  table[row sep=crcr]{%
0.142799666072263	0.967160314551492\\
0.0461998919645558	0.108621305175979\\
0.0149599650170943	0.0111860002772366\\
0.00484813681109536	0.00114616752327967\\
};

\addplot [color=mycolor2, line width=1.5pt, mark=diamond, mark options={solid, mycolor2}]
  table[row sep=crcr]{%
0.142799666072263	0.00666357162479992\\
0.0461998919645558	6.30285069133096e-05\\
0.0149599650170943	7.18100024998946e-07\\
0.00484813681109535	7.30856886121954e-09\\
};

\addplot [color=mycolor3, line width=1.5pt, mark=square, mark options={solid, mycolor3}]
  table[row sep=crcr]{%
0.142799666072263	1.69531633873454e-05\\
0.0461998919645558	1.30896512517964e-08\\
0.0149599650170943	1.73764336253157e-11\\
0.00484813681109535	1.17683640610267e-14\\
};

\addplot [color=mycolor4, line width=1.5pt, mark=x, mark options={solid, mycolor4}]
  table[row sep=crcr]{%
0.142799666072263	0.00806735923671264\\
0.0461998919645558	1.3656853425914e-12\\
0.0149599650170943	6.87228052242971e-14\\
0.00484813681109535	5.79536418854332e-14\\
};

\end{axis}

\begin{axis}[%
width=1.227\fwidth,
height=1.227\fheight,
at={(-0.16\fwidth,-0.135\fheight)},
scale only axis,
xmin=0,
xmax=1,
ymin=0,
ymax=1,
axis line style={draw=none},
ticks=none,
axis x line*=bottom,
axis y line*=left,
legend style={legend cell align=left, align=left, draw=white!15!black},
xlabel style={font={\scriptsize}},ylabel style={font=\scriptsize},  ylabel shift={-0cm},ticklabel style={font=\scriptsize}
]
\end{axis}
\end{tikzpicture}%
}
	\vspace{0.05cm}
	\centering
	\hspace{-0.4cm}
	{ \setlength{\fwidth}{9.0cm}
\setlength{\fheight}{2.0cm}
\definecolor{mycolor1}{rgb}{0.00000,0.44700,0.74100}%
\definecolor{mycolor2}{rgb}{0.85000,0.32500,0.09800}%
\definecolor{mycolor3}{rgb}{0.92900,0.69400,0.12500}%
\definecolor{mycolor4}{rgb}{0.49400,0.18400,0.55600}%
\begin{tikzpicture}

\begin{axis}[%
width=0.412\fwidth,
height=\fheight,
at={(0\fwidth,0\fheight)},
scale only axis,
xmode=log,
xmin=0.00484813681109536,
xmax=0.142799666072263,
xminorticks=true,
xlabel style={font=\color{white!15!black}},
xlabel={$h$},
ymode=log,
ymin=1e-15,
ymax=100,
yminorticks=true,
ylabel style={font=\color{white!15!black}},
ylabel={$E(T)$},
axis background/.style={fill=white},
xmajorgrids,
xminorgrids,
ymajorgrids,
yminorgrids,
legend style={at={(0.03,0.03)}, anchor=south west, legend cell align=left, align=left, draw=white!15!black,nodes={scale=0.75, transform shape}},
xlabel style={font={\scriptsize}},ylabel style={font=\scriptsize},  ylabel shift={-0cm},ticklabel style={font=\scriptsize}
]
\addplot [color=mycolor1, line width=1.5pt, mark=o, mark options={solid, mycolor1}]
  table[row sep=crcr]{%
0.142799666072263	0.26475728429479\\
0.0461998919645558	0.210085043276037\\
0.0149599650170943	0.0273387095995115\\
0.00484813681109536	0.00133178822365898\\
};
\addlegendentry{Lobatto-IIIC 2}

\addplot [color=mycolor2, line width=1.5pt, mark=diamond, mark options={solid, mycolor2}]
  table[row sep=crcr]{%
0.142799666072263	0.241953154400408\\
0.0461998919645558	0.0558617792249278\\
0.0149599650170943	0.0216606321920892\\
0.00484813681109536	0.0148215624746615\\
};
\addlegendentry{Lobatto-IIIC 4}

\addplot [color=mycolor3, line width=1.5pt, mark=square, mark options={solid, mycolor3}]
  table[row sep=crcr]{%
0.142799666072263	0.114684743967353\\
0.0461998919645558	0.0136934470055323\\
0.0149599650170943	0.0283919841092547\\
0.00484813681109536	0.00143691437248961\\
};
\addlegendentry{Lobatto-IIIC 6}

\addplot [color=mycolor4, line width=1.5pt, mark=x, mark options={solid, mycolor4}]
  table[row sep=crcr]{%
0.142799666072263	0.0648266409770485\\
0.0461998919645558	0.0150462926388084\\
0.0149599650170943	0.00184242092479914\\
0.00484813681109536	0.00721922624783133\\
};
\addlegendentry{Lobatto-IIIC 8}

\end{axis}

\begin{axis}[%
width=0.412\fwidth,
height=\fheight,
at={(0.542\fwidth,0\fheight)},
scale only axis,
xmode=log,
xmin=0.00484813681109536,
xmax=0.142799666072263,
xminorticks=true,
xlabel style={font=\color{white!15!black}},
xlabel={$h$},
ymode=log,
ymin=1e-15,
ymax=100,
yminorticks=true,
axis background/.style={fill=white},
xmajorgrids,
xminorgrids,
ymajorgrids,
yminorgrids,
legend style={legend cell align=left, align=left, draw=white!15!black},
xlabel style={font={\scriptsize}},ylabel style={font=\scriptsize},  ylabel shift={-0cm},ticklabel style={font=\scriptsize}
]
\addplot [color=mycolor1, line width=1.5pt, mark=o, mark options={solid, mycolor1}]
  table[row sep=crcr]{%
0.142799666072263	0.806427436201109\\
0.0461998919645558	0.159023869384874\\
0.0149599650170943	0.0205149041750922\\
0.00484813681109536	0.00224079696586599\\
};

\addplot [color=mycolor2, line width=1.5pt, mark=diamond, mark options={solid, mycolor2}]
  table[row sep=crcr]{%
0.142799666072263	0.00838200685231515\\
0.0461998919645558	0.000281903622194269\\
0.0149599650170943	1.05668719585506e-06\\
0.00484813681109535	1.22776131483704e-08\\
};

\addplot [color=mycolor3, line width=1.5pt, mark=square, mark options={solid, mycolor3}]
  table[row sep=crcr]{%
0.142799666072263	0.00807756877787746\\
0.0461998919645558	2.12798996201968e-08\\
0.0149599650170943	6.95116408877537e-10\\
0.00484813681109535	1.98407956730762e-12\\
};

\addplot [color=mycolor4, line width=1.5pt, mark=x, mark options={solid, mycolor4}]
  table[row sep=crcr]{%
0.142799666072263	0.00809598846442605\\
0.0461998919645558	1.31172850359462e-12\\
0.0149599650170943	3.13082892944294e-14\\
0.00484813681109535	2.26485497023531e-14\\
};

\end{axis}

\begin{axis}[%
width=1.232\fwidth,
height=1.232\fheight,
at={(-0.16\fwidth,-0.139\fheight)},
scale only axis,
xmin=0,
xmax=1,
ymin=0,
ymax=1,
axis line style={draw=none},
ticks=none,
axis x line*=bottom,
axis y line*=left,
legend style={legend cell align=left, align=left, draw=white!15!black},
xlabel style={font={\scriptsize}},ylabel style={font=\scriptsize},  ylabel shift={-0cm},ticklabel style={font=\scriptsize}
]
\end{axis}
\end{tikzpicture}%
}
	\vspace{0.01cm}
	\centering
	{\setlength{\fwidth}{9.0cm}
\setlength{\fheight}{2.0cm}
\definecolor{mycolor1}{rgb}{0.00000,0.44700,0.74100}%
\definecolor{mycolor2}{rgb}{0.85000,0.32500,0.09800}%
\definecolor{mycolor3}{rgb}{0.92900,0.69400,0.12500}%
\definecolor{mycolor4}{rgb}{0.49400,0.18400,0.55600}%
\begin{tikzpicture}

\begin{axis}[%
width=0.411\fwidth,
height=\fheight,
at={(0\fwidth,0\fheight)},
scale only axis,
xmode=log,
xmin=0.00484813681109536,
xmax=0.142799666072263,
xminorticks=true,
xlabel style={font=\color{white!15!black}},
xlabel={$h$},
ymode=log,
ymin=1e-15,
ymax=100,
yminorticks=true,
ylabel style={font=\color{white!15!black}},
ylabel={$E(T)$},
axis background/.style={fill=white},
xmajorgrids,
xminorgrids,
ymajorgrids,
yminorgrids,
legend style={at={(0.03,0.03)}, anchor=south west, legend cell align=left, align=left, draw=white!15!black,nodes={scale=0.75, transform shape}},
xlabel style={font={\scriptsize}},ylabel style={font=\scriptsize},  ylabel shift={-0cm},ticklabel style={font=\scriptsize}
]
\addplot [color=mycolor1, line width=1.5pt, mark=o, mark options={solid, mycolor1}]
  table[row sep=crcr]{%
0.142799666072263	8.31506747724377\\
0.0461998919645558	3.86775333991472\\
0.0149599650170943	2.37365546309777\\
0.00484813681109536	0.823788236070209\\
};
\addlegendentry{Explicit Euler 1}

\addplot [color=mycolor2, line width=1.5pt, mark=diamond, mark options={solid, mycolor2}]
  table[row sep=crcr]{%
0.142799666072263	3.76855314346167\\
0.0461998919645558	0.36567116917783\\
0.0149599650170943	0.0293217113639036\\
0.00484813681109536	0.00151675541998186\\
};
\addlegendentry{Heun 2}

\addplot [color=mycolor3, line width=1.5pt, mark=square, mark options={solid, mycolor3}]
  table[row sep=crcr]{%
0.142799666072263	0.63070960722971\\
0.0461998919645558	0.024159448597602\\
0.0149599650170943	0.046879020725453\\
0.00484813681109536	0.0136574016600953\\
};
\addlegendentry{Heun 3}

\addplot [color=mycolor4, line width=1.5pt, mark=x, mark options={solid, mycolor4}]
  table[row sep=crcr]{%
0.142799666072263	0.230021416780682\\
0.0461998919645558	0.0577620464555856\\
0.0149599650170943	0.0214792886705394\\
0.00484813681109536	0.0148023190779655\\
};
\addlegendentry{Runge-Kutta 4}

\end{axis}

\begin{axis}[%
width=0.411\fwidth,
height=\fheight,
at={(0.54\fwidth,0\fheight)},
scale only axis,
xmode=log,
xmin=0.00484813681109536,
xmax=0.142799666072263,
xminorticks=true,
xlabel style={font=\color{white!15!black}},
xlabel={$h$},
ymode=log,
ymin=1e-15,
ymax=100,
yminorticks=true,
axis background/.style={fill=white},
xmajorgrids,
xminorgrids,
ymajorgrids,
yminorgrids,
legend style={legend cell align=left, align=left, draw=white!15!black},
xlabel style={font={\scriptsize}},ylabel style={font=\scriptsize},  ylabel shift={-0cm},ticklabel style={font=\scriptsize}
]
\addplot [color=mycolor1, line width=1.5pt, mark=o, mark options={solid, mycolor1}]
  table[row sep=crcr]{%
0.142799666072263	8.17395439190817\\
0.0461998919645558	3.58529097537404\\
0.0149599650170943	2.54523707657732\\
0.00484813681109536	0.847707017144985\\
};

\addplot [color=mycolor2, line width=1.5pt, mark=diamond, mark options={solid, mycolor2}]
  table[row sep=crcr]{%
0.142799666072263	2.62940871881158\\
0.0461998919645558	0.311660324065769\\
0.0149599650170943	0.0585255326757082\\
0.00484813681109536	0.00284953822375922\\
};

\addplot [color=mycolor3, line width=1.5pt, mark=square, mark options={solid, mycolor3}]
  table[row sep=crcr]{%
0.142799666072263	0.0364569764121518\\
0.0461998919645558	0.0137916440498528\\
0.0149599650170943	0.000509409552742146\\
0.00484813681109535	1.68933467172838e-05\\
};

\addplot [color=mycolor4, line width=1.5pt, mark=x, mark options={solid, mycolor4}]
  table[row sep=crcr]{%
0.142799666072263	0.0439031081394205\\
0.0461998919645558	0.00828320728711973\\
0.0149599650170943	0.000485019629422711\\
0.00484813681109535	4.69029565297597e-07\\
};

\end{axis}

\begin{axis}[%
width=1.227\fwidth,
height=1.227\fheight,
at={(-0.16\fwidth,-0.135\fheight)},
scale only axis,
xmin=0,
xmax=1,
ymin=0,
ymax=1,
axis line style={draw=none},
ticks=none,
axis x line*=bottom,
axis y line*=left,
legend style={legend cell align=left, align=left, draw=white!15!black},
xlabel style={font={\scriptsize}},ylabel style={font=\scriptsize},  ylabel shift={-0cm},ticklabel style={font=\scriptsize}
]
\end{axis}
\end{tikzpicture}%
}
	\vspace{-0.3cm}
	\caption{Integration error $E(T) = \| x(T)-\hat{x}_h(T)\|$ vs. the step size $h$ for different RK methods: standard (left plot) vs. FESD (right plot). The legend provides the name of the used RK method and its order of the global integration error.}
	\label{fig:integrator_order}
\end{figure}
We regard solution approximations to this IVP obtained by standard explicit and implicit RK methods \eqref{eq:dcs_irk} and FESD \eqref{eq:fesd_compact} with $\NFE=2$ and different step sizes. 
Both integration methods are available in \texttt{NOSNOC} via the function \texttt{integrator\_fesd()}. 
The regarded RK methods are listed in Table \ref{tab:irk_schemes} together with their global integration error estimate.  
For explicit methods, we consider in this paper $\Nstg\leq 4$. Several other IRK methods are available in \texttt{NOSNOC}, but we omit the full comparison for brevity.
The nominal step $h$ size is obtained by dividing $T$ by the number of simulation steps $N_{\mathrm{sim}}$. We take an irrational number for $T=\frac{\pi}{2}$ so that $\ts=1$ never coincides with a finite element boundary. Therefore, we avoid accidental switch detection via the discretization grid for the standard discretization.

\begin{table}
	\centering
	\begin{tabular}{|c|c|}	
		\hline
		\textbf{Method}  & \textbf{Global error estimate }   \\
		\hline
		Radau-IIA & $h^{2\Nstg-1}$\\
		\hline
		Gauss-Legendre &  $h^{2\Nstg}$\\
		\hline
		Lobatto-IIIA  & $h^{2\Nstg-2}$\\
		\hline
		Lobatto-IIIC & $h^{2\Nstg-2}$\\
		\hline
		Explicit-RK  & $h^{\Nstg}$\\
		\hline
	\end{tabular}
	\caption{List of analyzed RK methods and their accuracy order for ODE \cite{Hairer1991}. Note that for Explicit-RK methods the assertion in the table is true for $\Nstg\leq 4$, otherwise $p < \Nstg$.}
	\label{tab:irk_schemes}
\end{table}

The numerical error as a function of the step sizes for the RK and FESD methods from Table \ref{tab:irk_schemes} are depicted in Figure \ref{fig:integrator_order}. 
We can see that the standard discretization achieves in all cases only first-order accuracy (left plots). 
In contrast, the FESD method recovers in all cases the high accuracy order that the underlying RK method has for smooth ODE (right plots). 
This verifies the result of Theorem \ref{th:integration_order} in practice and demonstrates how FESD can be used as an event-based integrator without an external switch detection procedure.  
The \textit{saturation} in the right plots of some high-accuracy methods is due to the round-off errors in floating point arithmetic which limit the possible accuracy on a computer. 

\subsection{Convergence of discrete-time sensitivities} \label{sec:sensitivity_convergence}
One of the most important tools for computing stationary points and verifying their optimality in direct optimal control is the numerical sensitivities, i.e., the derivatives of the numerical solution approximation $\hat{x}_h(t;x_0)$ w.r.t. the initial value $x_0$ (and parameters) for some $t\in\left(0,T\right]$. 
Here we denote them by $\hat{X}_h(t;0,x_0)$ or sometimes more compactly by $\hat{X}_h(t,x_0)$.

A fundamental limitation of standard discretization methods for nonsmooth systems (e.g., with the RK discretization from Subsection \ref{sec:dcs_irk}) is that $\hat{X}_h(t,x_0) \not\to {X}(t,x_0)$ as $h\downarrow 0$ \cite{Stewart2010}. 
In direct optimal control, this can cause convergence to artificial stationary points arbitrarily close to the initialization point \cite{Nurkanovic2020,Stewart2010}. 
Fortunately, the sensitivities of the solutions generated by the FESD method converge to their true values (cf. Subsection \ref{sec:sensitivites}).
This is shown in the next theorem, but before we state it, we make one more assumption on the time derivatives of the solution approximation.

\begin{assumption}(RK derivatives)
	\label{ass:irk_scheme_derivative}
	Regard the RK methods from Assumption \ref{ass:irk_scheme} applied to the differential algebraic equations \eqref{eq:dcs_dae}. 
	Suppose that the derivatives of the numerical approximation for the same RK method converge with order $1 \leq q\leq p$, i.e., $\| \dot{\hat{x}}_h(t) - \dot{x}(t) \| = O(h^q),\; t\in \mathcal{G}$.
\end{assumption}
The aim of the assumption about the convergence of the derivatives of the numerical approximation is to cover a broad class of RK methods and the value of $q$ depends on the specific choice of an RK method. 
For example, for collocation-based implicit RK methods for ODE in general it holds that $q = p-1$ \cite[Theorem 7.10]{Hairer1993}. 
More specifically, methods that contain the boundary point of a finite element denoted by $\hat{t}$, that is $c_{\Nstg} =1$ satisfy $p = q$. This assertion follows directly from Lipschitz continuity and the fact that the numerical approximations satisfy the ODE at every stage point:  \\
$
\|  \dot{\hat{x}}_h(\bar{t}) - \dot{x}(\bar{t})\| = \|  f_i(\hat{x}_h(\bar{t})) - f_i({x}(\bar{t}))\|  \leq L_{f_i} \|  \hat{x}_h(\bar{t}) - {x}(\bar{t})\| = O(h^p).
$

\begin{theorem1}[Convergence to exact sensitivities]\label{th:exact_sensitivites}
Suppose the assumptions of Theorem \ref{th:integration_order} and Assumption \ref{ass:irk_scheme_derivative} hold. Assume that a single active-set change happens at time $t_{\mathrm{s},n}$, i.e., $||\I_n| - |\I_{n+1}|| \leq 1, n = 0,\dots \Nswitch$. Then for $h\downarrow 0$ it holds that $\hat{X}_h(t,x_0)  \to  X(t,x_0)$ with the convergence rate
\begin{align}
\| \hat{X}_h(t,x_0) - X(t,x_0)\| = O(h^{q}), \text{ for all } t\in \mathcal{G}.
\end{align}
\end{theorem1}
\textit{Proof.} Regard partition of $[0,T]$: $0 < \tilde{t}_1 <\dots < \tilde{t}_k < \dots <\tilde{t}_{N_k}  < T$ such that in the open interval between every two neighboring points there is a single switching point with $N_k>\Nswitch$. Then by the chain rule we have
\begin{align*}
	X(T;0,x_0) = X(T;\tilde{t}_{N_k},x(\tilde{t}_{N_k}))\cdots X(\tilde{t}_{k+1};\tilde{t}_{k},x(\tilde{t}_k))\cdots X(\tilde{t}_1;0,x_0).
\end{align*}
W.l.o.g.  assume that on $[0,\tilde{t}_1]$ a single switch occurs at $\ts \in (0,\tilde{t}_1)$. We show convergence of the sensitivities on this interval. Convergence on $[0,T]$ follows by inductively applying the same arguments on every sub-interval $[\tilde{t}_{k},\tilde{t}_{k+1}]$.

Regard the two smooth pieces of the approximation $\hat{x}_h(t)$: 
1) $\hat{x}_{h,1}(t,x_0)$ for $t\leq \hatts$ and, 
2) $\hat{x}_{h,2}(t,y_0(x_0))$ where $y_0(x_0) =\hat{x}_{h,2}(0,y_0(x_0)) = \hat{x}_{h,1}(\hatts,x_0)$. With this definition, we have for $t \geq \hatts$
\begin{align}\label{eq:conitnuity_x_hat}
	\hat{x}_{h,2}(t-\hatts,y_0(x_0)) &= \hat{x}_h(t,x_0).
\end{align}
From Theorem \ref{th:integration_order} we know that $| \hatts - \ts | =O(h^p)$. Obviously, the value of $\hatts$ depends on $x_0$ and we know from Theorem \ref{th:integration_order} that we obtain implicitly at a switching point the condition
\begin{align}\label{eq:senstivity_at_switch}
	\psi_{i,j}(\hat{x}_{h,1}(\hatts(x_0),x_0))&= g_i(\hat{x}_{h,1}(\hatts(x_0),x_0))  - g_j(\hat{x}_{h,1}(\hatts(x_0),x_0))=0,
\end{align} 
where $i\notin \I_0$, $i \in \I_1$ and $j \in \I_0$.

For computing $\hat{X}_h(\cdot)$ on $[\hat{t}_{\mathrm{s}},t_1]$ from Eq. \eqref{eq:conitnuity_x_hat} we have
\begin{align}\label{eq:sensitivity_connecting}
\begin{split}
\partialder{\hat{x}_h(t,x_0)}{x_0} &= \partialder{\hat{x}_{h,2}(t-\hatts(x_0),y_0(x_0))}{x_0} 
= - \dot{\hat{x}}_{h,2}(t-\hatts(x_0),y_0(x_0)) \nabla_{x_0} \hatts(x_0)^\top \\&+ 
\partialder{\hat{x}_{h,2}(t-\hatts(x_0),y_0(x_0))}{y_0} \nabla_{x_0} y_0(x_0).
\end{split}
\end{align}
Next, we compute the expressions for the two unknowns  $\nabla_{x_0} \hatts(x_0)^\top$ and $\nabla_{x_0} y_0(x_0)$.
Denote by $\hat{X}_{h,1}(t;0,x_0) = \partialder{\hat{x}_{h,1}(t,x_0)}{x_0}$. 
Using the implicit function theorem for \eqref{eq:senstivity_at_switch}, we can compute
\begin{align}\label{eq:sensitivity_t_s_x_0}
\nabla_{x_0} \hatts(x_0)^\top = - \frac{\nabla \psi_{i,j}(\hat{x}_{h,1}(\hatts(x_0),x_0))^\top \hat{X}_{h,1}(\hatts,x_0)}{\nabla \psi_{i,j}(\hat{x}_{h,1}(\hatts(x_0),x_0))^\top \dot{\hat{x}}_{h,1}(\hatts(x_0),x_0)}.
\end{align}
At $t = \hatts$, we exploit the fact that $\hat{x}_{h,2}(0,y_0(x_0))  = y_0(x_0) = \hat{x}_{h,1}(\hatts,x_0)$, thus
\begin{align*}
\nabla_{x_0} y_0(x_0)= \partialder{\hat{x}_{h,1}(\hatts(x_0),x_0)}{x_0} = \dot{\hat{x}}_{h,1}(\hatts(x_0),x_0)  \nabla_{x_0} \hatts(x_0)^\top
+ \hat{X}_{h,1}(\hatts;x_0).
\end{align*}
Combing the last line with \eqref{eq:sensitivity_t_s_x_0} we obtain
\begin{align}\label{eq:sensitivity_y0_x0}
	\nabla_{x_0} y_0(x_0)= \Bigg[ I - \frac{ \dot{\hat{x}}_{h,1}(\hatts(x_0),x_0)\nabla \psi_{i,j}(\hat{x}_{h,1}(\hatts(x_0),x_0))^\top} {\nabla \psi_{i,j}(\hat{x}_{h,1}(\hatts(x_0),x_0))^\top \dot{\hat{x}}_{h,1}(\hatts(x_0),x_0)}\Bigg] \hat{X}_{h,1}(\hatts;x_0).
\end{align}
We are interested in $\partialder{x_h(t;x_0)}{x_0}$ when $t \downarrow \hatts$. 
Note that $ \partialder{\hat{x}_{h,2}(t-\hatts(x_0),y_0(x_0))}{y_0} \to I$ as $t \downarrow \hatts$. 
Thus from \eqref{eq:sensitivity_connecting}, \eqref{eq:sensitivity_t_s_x_0} and \eqref{eq:sensitivity_y0_x0} one obtains 
\begin{align}\label{eq:discrete_time_sens_jump}
\begin{split}
\partialder{x_h(\hatts^+;x_0)}{x_0} & = \hat{J}_h(\hat{x}_h(\hat{\ts};x_0)) 
\partialder{x_h(\hatts^-;x_0)}{x_0},\ \text{with}\\
\hat{J}_h(\hat{x}_h(\hatts;x_0)) &\coloneqq \Bigg[I + \frac{ (\dot{\hat{x}}_{h}(\hatts^+,x_0)- \dot{\hat{x}}_{h}(\hatts^-,x_0))\nabla \psi_{i,j}(\hat{x}_{h}(\hatts^-,x_0))^\top} {\nabla \psi_{i,j}(\hat{x}_{h}(\hatts^-,x_0))^\top \dot{\hat{x}}_{h}(\hatts^-,x_0)}\Bigg].
\end{split}
\end{align}
By the chain rule, we have that for $ t> \hatts$
\begin{align}
	\hat{X}_{h}(t;x_0) = \hat{X}_{h,2}(t;\hatts^+,y_0)\hat{J}_h(\hat{x}_h(\hatts;x_0)) \hat{X}_{h,1}(\hat{t}_{\mathrm{s}}^-;0,x_0).
\end{align}

First, note that for a fixed active set the FESD equations for $\hat{x}_{h,1}(t,x_0)$ and $\hat{x}_{h,2}(t,y_0)$ boil down to RK equations for the DAE \eqref{eq:dcs_dae} with fixed step size $h_n$, cf. Theorem \ref{th:fesd_locally_unique_solution}. 
Differentiating the RK equations to obtain $\hat{X}_h(\cdot)$ results in the same RK method applied to the variational differential equations of the system at hand, thus the numerical sensitivities converge in this setting to the continuous-time sensitivities with the same accuracy $O(h^p)$ as for the solution of the system \cite{Albersmeyer2010b}. 

Second, as $h \downarrow 0$, due to assumption of this theorem, in $\hat{J}_h(\hat{x}_h(\hatts;x_0))$ the functions $\dot{x}_h(\cdot)$ converge to $f(x(\cdot))$ with order $q$. 
Due to Theorem \ref{th:integration_order} all other terms converge with order $p$. Thus, $\|\hat{J}_h(\hat{x}_h(\hatts;x_0)) - {J}({x}(\ts;x_0))\| = O(h^{q})$. 

Summarizing the last two arguments and applying them inductively for every active-set change we conclude that $\hat{X}_{h}(t;x_0)\to{X}(t;x_0)$ as $h\downarrow 0$ with the order $q=\min(p,q)$. 
This completes the proof.\qed

The only restrictive assumption we make is that a single active-set change happens at a time. 
For multiple simultaneous switches the derivation becomes quite involved even in continuous-time case \cite{Filippov1988}, hence we made this assumption for simplicity. 
The results of Theorem \ref{th:exact_sensitivites} are illustrated on a numerical example in the next subsection.
\subsection{Illustration of numerical sensitivity convergence}\label{sec:ecc_problem}
\begin{figure}[t]
	\centering
	{ \input{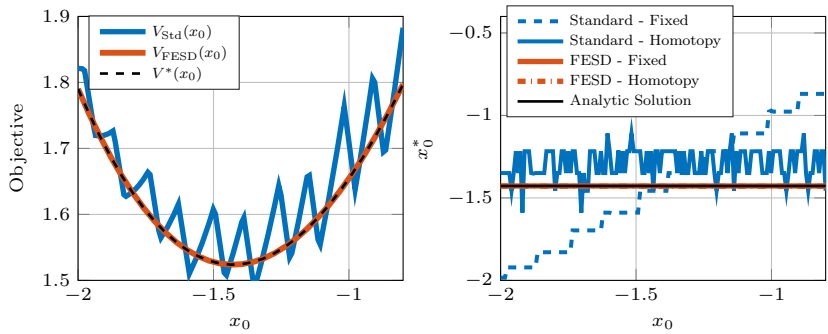}}
	\caption{Approximations of the objective $V(x_0)$ computed with FESD and a standard IRK Gauss-Legendre method compared to the true value are shown in the left plot. The right plot shows the value of the optimal solution $x^*_0$ against different initial values $x_0$ for which a initial feasible solution guess was computed. The optimal solution is computed with FESD and a standard method with and without a homotopy approach for the underlying optimization problem.}
	\label{fig:ivp_problem}
		\vspace{-0.45cm}					
\end{figure}

To demonstrate the improvements in FESD compared to standard methods we repeat the experiments from \cite{Nurkanovic2020}. 
For this purpose, we look at the optimal control problem from~\cite{Stewart2010}:
\begin{subequations}\label{eq:ivp_ocp}
	\begin{align}
		\min_{x_0\in \R, x(\cdot) \in \mathcal{C}^0} \quad & \int_0^2 \! x(t)^2 \, \dd t  +  (x(2) - 5/3)^2\\
		\textrm{s.t.} \quad &  x(0)= x_0,\\
		&\dot{x}(t) \in 2-\mathrm{sign}(x(t)),\ t \in [0,2] \label{eq:ivp_di}.
	\end{align}
\end{subequations}	
Note that the initial value $x_0$ is the only effective degree of freedom. Let $V_*(x_0)$ be the objective value for the unique feasible trajectory starting at $x(0)=x_0$. The equivalent reduced problem is given by 
\begin{align} \label{eq:ivp_ocp_reduced}
	\min_{x_0\in \R} \quad & V_*(x_0). 
\end{align}
In the first experiment, we evaluate $V_*(x_0)$ by simulating the trajectory of \eqref{eq:ivp_di} for various $x_0$ with a standard RK method \eqref{eq:dcs_irk} and FESD \eqref{eq:fesd_compact}. 
We pick the Gauss-Legendre method of fourth order with $\Nstg=2$ and $\NFE = 25$. 
The results are depicted in the left plot in Figure \ref{fig:ivp_problem}. 
It can be seen that for standard RK methods, the trajectories would converge to the analytic solution, but the derivatives w.r.t. $x_0$ do not. 
In fact, many artificial minima arise \cite{Nurkanovic2020,Stewart2010}. 
On the other hand, the FESD solution matches the analytic solution.

In the second experiment, we solve the OCP \eqref{eq:ivp_ocp} for different initial guesses $x_0$. 
The initial guess for the arising MPCC is obtained from a forward simulation of \eqref{eq:ivp_di} with the same method and same grid as in the discretized OCP. 
The arising MPCC are solved by a relaxation approach, where the complementarity constraints \eqref{eq:mpcc_GGENRAL_comp} are replaced by $w_1^\top w_2 \leq \sigma,\; w_1\geq0,\;w_2\geq0$. 
In the first experiment, the MPCC are solved for a fixed $\sigma = 10^{-15}$. 
In the second case, we solve a sequence of NLP in a homotopy procedure where $\sigma$ is decreased from $\sigma_0=1$ to $\sigma_N = 10^{-15}$, via the rule $\sigma_{k+1} = 0.1 \sigma_k$. 
It was demonstrated in \cite{Nurkanovic2020} that a homotopy approach with standard discretizations can lead to convergence to better local minima. The results are depicted in the right plot of Figure \ref{fig:ivp_problem}.

Remarkably, the FESD solution matches the analytic solution for all $x_0$ regardless of the MPCC strategy and initialization point (FESD - fixed and FESD - homotopy). 
This highlights that the numerical sensitivities converge to the true ones, cf. Theorem \ref{th:exact_sensitivites}.
On the other hand, as in \cite{Nurkanovic2020}, with the standard IRK method, the fixed parameter MPCC strategy leads to convergence to the artificial local minima close to the initialization. 
This results in a stair-like curve in right plot in Figure \ref{fig:ivp_problem} (Standard - fixed). 
This is a consequence of the fact that the numerical derivatives are always wrong, i.e., the numerical approximation is "trapped" in one of the local minima of the standard approach, which are visible in the left plot (blue).

The homotopy strategy yields better local minima since in the earlier homotopy iterations the derivatives are still correct \cite{Nurkanovic2021} (Standard - homotopy). 
However, the error of $O(h)$ to the analytic solution is still present since a standard method has at best accuracy of order one. 
This highlights the result of Theorems \ref{th:integration_order} and \ref{th:exact_sensitivites} and demonstrates the superiority of FESD to the na\"ive approach even on a very simple example.  
\section{FESD in direct optimal control}\label{sec:optimal_control} 
This section regards the use of FESD in direct optimal control, i.e., a first discretize then optimize approach.
We consider the following continuous-time optimal control problem on a control horizon $[0,\Tctrl]$:
\begin{subequations} \label{eq:OCP}
	\begin{align}
		\min_{x(\cdot),u(\cdot),z(\cdot)} \quad & \int_{0}^{\Tctrl} L_{\mathrm{int}}(x(t),u(t))\dd  t +  L_{\mathrm{end}}(x(\Tctrl)) \label{eq:OCP_objective}\\
		\textrm{s.t.} \quad  x_{0} &= \bar{x}_0 \label{eq:OCP_IV},\\
		\dot{x}(t) & = F(x(t),u(t))\theta(t), & t \in [0,\Tctrl], \label{eq:OCP_ode}\\
		0& = g(x(t)) - \lambda(t) - \mu(t)e, & t \in [0,\Tctrl], \label{eq:OCP_alg} \\
		1& = e^\top \theta(t),& t \in [0,\Tctrl], \label{eq:OCP_simplex}\\
		0& \leq \theta(t) \perp \lambda(t) \geq 0,& t \in [0,\Tctrl], \label{eq:OCP_cc}\\
		0&\leq G_{\mathrm{ineq}}(x(t),u(t)),& t \in [0,\Tctrl] \label{eq:OCP_path},\\
		0&\leq G_{\mathrm{end}}(x(\Tctrl)),&  \label{eq:OCP_term}
	\end{align}
\end{subequations}
where $\bar{x}_0$ is a given initial value, $u(t) \in \R^{n_u}$ is the control function, $z(t) = (\lambda(t),\theta(t),\mu(t)) \in \R^{2\Nsys+1}$ collects the algebraic variables. 
The function  $L_{\mathrm{int}}: \R^{n_x} \times \R^{n_u} \to \R$ is the Lagrange objective term to be integrated and $L_{\mathrm{end}}:\R^{n_x}\to \R$ is the terminal cost also called Mayer term. 
The path and terminal constraints are collected in the functions $G_{\mathrm{ineq}} : \R^{n_x}  \times \R^{n_u} \to \R^{n_{g1}}$ and $G_{\mathrm{end}} : \R^{n_x}  \to \R^{n_{g2}}$, respectively.
\subsection{A multiple shooting-type discretization}\label{sec:multiple_shooting_grids}
We discretize the OCP \eqref{eq:OCP} with $\Nctrl\geq 1$ control intervals indexed by $k$. 
The control function approximation is taken to be piecewise constant on an equidistant grid.
This is the usual choice in direct multiple shooting \cite{Bock1984} and is required by many practical applications of feedback control, but could be easily generalized to any other local control parameterization. 
The constant controls are collected in $\mathbf{q} = (q_0,\ldots,q_{\Nctrl-1})\in \R^{\Nctrl n_u}$. 
All internal variables of every control interval are additionally equipped with an index $k$. 
On every control interval $k$ with fixed duration $T = \frac{\Tctrl}{\Nctrl}$, we apply a discretization \eqref{eq:fesd_compact} with $N_{\mathrm{FE}}$ internal finite elements. The state values at the control interval boundaries are collected in  $\mathbf{s} = (s_0,\ldots,s_{\Nctrl})\in\R^{(\Nctrl+1)n_x}$. 
The following vectors collect all internal variables of all discretization steps: 
$\mathcal{H} = (\mathbf{h}_0,\ldots,\mathbf{h}_{\Nctrl-1})$ and
${\mathcal{Z}} = (\mathbf{{Z}}_0,\ldots,\mathbf{{Z}}_{\Nctrl-1})$.
The FESD discretization of the OCP \eqref{eq:OCP} together with the aforementioned control discretization reads as:
\begin{subequations}\label{eq:OCP_FESD}
	\begin{align}
		\min_{\mathbf{s},\mathbf{q},\mathcal{H},\mathcal{Z}} \quad & \sum_{k=1}^{
			\Nctrl-1} \hat{L}_{\mathrm{int}}(s_k,\mathbf{x}_k,q_k) +  \hat{L}_{\mathrm{end}}(s_{\Nctrl}) \label{eq:OCP_objective}\\
		\textrm{s.t.} \quad  &s_{0} = \bar{x}_0,\\
		&{s}_{k+1}  = F_{\fesd}(\mathbf{Z}_k),\; &k = 0,\ldots,\Nctrl-1,\\
		&0=G_{\fesd}(\mathbf{Z}_k,\mathbf{h}_k,s_k,q_k,T_k),\; &k = 0,\ldots,\Nctrl-1, \label{eq:OCP_FESD_internal}\\
		&0\leq G_{\mathrm{ineq}}(\mathbf{x}_k,q_k),\; &k = 0,\ldots,\Nctrl-1,\\
		& h_{\mathrm{min}}e \leq \mathbf{h}_k \leq h_{\mathrm{max}}e,\; &k = 0,\ldots,\Nctrl-1, \label{eq:OCP_FESD_box_h}\\
		&0\leq G_{\mathrm{end}}(s_{\Nctrl}).
	\end{align}
\end{subequations}
where $\hat{L}_{\mathrm{int}}:\R^{n_x}\times \R^{(\NFE+1)\Nstg n_x} \times \R^{n_u}\to \R$ and $\hat{L}_{\mathrm{end}}:\R^{n_x} \to \R$ are the discretized integral and terminal costs, respectively. The scalars $h_{\mathrm{min}}$ and $h_{\mathrm{max}}$ are the lower and upper bounds for the step sizes. The box constraint \eqref{eq:OCP_FESD_box_h}, prohibits negative step sizes and bounds the variability of the step size and thus the local discretization errors. We want to recall that with this formulation at every control interval the constraint $\sum_{n=1}^{\NFE} h_{n,k} =T$ is imposed as part of \eqref{eq:OCP_FESD_internal}, cf. Eq \eqref{eq:step_eq}.

\subsection{Solving the discretized OCP} \label{sec:solving_mpccs}
When solving the OCP \eqref{eq:OCP_FESD} in practice we do not use C-functions {(cf. Section \ref{sec:introduction} for a definition)}, but write the complementarity conditions explicitly. 
Therefore, the NLP \eqref{eq:OCP_FESD} is a Mathematical Program with Complementarity Constraints (MPCC) and it can be compactly written as 
\begin{subequations} \label{eq:mpcc_GGENRAL}
	\begin{align}
		\min_{w} \quad & \varphi(w)\\
		\textrm{s.t.} \quad  
		0 & \leq \zeta(w),\\
		0&\leq w_1 \perp w_2\geq 0, \label{eq:mpcc_GGENRAL_comp}
	\end{align}
\end{subequations} 
where $w=(w_0,w_1,w_2) \in \mathbb{R}^{n_w}$ is a decomposition of the problem variables.
MPCC are difficult nonsmooth optimization problems that violate the {Mangasarian-Fromovitz constraint qualification} at all feasible points \cite{Anitescu2007}, which makes standard NLP algorithms fail to converge.
Fortunately, MPCC can often be solved efficiently via reformulations and homotopy approaches \cite{Anitescu2007,Hall2022,Leyffer2006,Ralph2004}.  
In a homotopy procedure a sequence of more regular, smooth, and relaxed NLPs related to \eqref{eq:mpcc_GGENRAL} are solved. 
Several relaxation, smoothing, and penalty homotopy approaches are implemented in MATLAB package \texttt{NOSNOC} \cite{Nurkanovic2022b}. 
They differ in how the complementarity constraints \eqref{eq:mpcc_GGENRAL_comp} are treated. 
Under some regularity assumptions the complementarity constraints are satisfied exactly even after solving a finite sequence of NLP \cite{Anitescu2007,Ralph2004}.  
In \texttt{NOSNOC} the NLP is solved with \texttt{IPOPT} \cite{Waechter2006} via its \texttt{CasADi} \cite{Andersson2019} interface.
 \begin{figure}[b]
	\centering
	{ \pgfplotsset{compat=1.13}
\setlength{\fwidth}{9.5cm}
\setlength{\fheight}{2.45cm}
\definecolor{mycolor1}{rgb}{0.00000,0.44700,0.74100}%
\definecolor{mycolor2}{rgb}{0.85000,0.32500,0.09800}%
\begin{tikzpicture}

\begin{axis}[%
width=0.45\fwidth,
height=\fheight,
at={(0\fwidth,0\fheight)},
scale only axis,
xmin=-2.2,
xmax=2.2,
xlabel style={font=\color{white!15!black}},
xlabel={$x_1$},
ymin=-1.5,
ymax=1.5,
ylabel style={font=\color{white!15!black}},
ylabel={$x_2$},
axis background/.style={fill=white},
xmajorgrids,
ymajorgrids,
legend style={at={(0.515,0.042)}, anchor=south west, legend cell align=left, align=left, draw=white!15!black,nodes={scale=0.70, transform shape}},
xlabel style={font={\scriptsize}},ylabel style={font=\scriptsize},  ylabel shift={-0cm},ticklabel style={font=\scriptsize}
]
\addplot [color=mycolor1, line width=2.0pt]
  table[row sep=crcr]{%
2.0943951023932	1.0471975511966\\
1.42786881360382	0.379539275607538\\
1.11349927826026	0.0690303100381855\\
0.965744309115852	0.0450356543501567\\
0.812652335148182	0.026833935234793\\
0.760781013305388	0.0220165365546214\\
0.587015921682928	0.0101139230901244\\
0.367963388329131	0.00249105795914994\\
0.0956585032885018	4.37663922396858e-05\\
-0	-1.4210854715202e-14\\
-0.0027151213003318	-0.134539245834261\\
-0.00717649270981813	-0.218731078571998\\
-0.0198542965536457	-0.363815855923762\\
-0.0392777969560369	-0.508654894237407\\
-0.0568608969680882	-0.555655341077288\\
-0.0852830105966209	-0.596509827145444\\
-0.143226636405178	-0.647875249188915\\
-0.223048559094138	-0.693740348743471\\
-0.324748778663499	-0.73410512580911\\
-0.448327295113262	-0.768969580385833\\
-0.523598775606894	-0.785398163404992\\
};
\addlegendentry{$x(t)$}

\addplot [color=black]
  table[row sep=crcr]{%
-0.342015	-1.51\\
-0.273375	-1.35\\
-0.212415	-1.19\\
-0.159135	-1.03\\
-0.113535	-0.87\\
-0.075615	-0.71\\
-0.04704	-0.56\\
-0.025215	-0.41\\
-0.01014	-0.26\\
-0.0018149999999999	-0.11\\
-0.000240000000000018	0.04\\
-0.00541499999999995	0.19\\
-0.0173399999999999	0.34\\
-0.0360150000000001	0.49\\
-0.0614399999999999	0.64\\
-0.093615	0.79\\
-0.13254	0.94\\
-0.1815	1.1\\
-0.23814	1.26\\
-0.30246	1.42\\
-0.342015	1.51\\
};
\addlegendentry{$\varphi_1(x)=0$}

\addplot [color=black]
  table[row sep=crcr]{%
-2.21	-0.53969305\\
-2.09	-0.45646645\\
-1.97	-0.38226865\\
-1.85	-0.31658125\\
-1.73	-0.25888585\\
-1.61	-0.20866405\\
-1.48	-0.1620896\\
-1.35	-0.12301875\\
-1.21	-0.0885780500000002\\
-1.07	-0.0612521500000001\\
-0.92	-0.0389344\\
-0.75	-0.0210937499999999\\
-0.57	-0.00925965000000017\\
-0.35	-0.00214375000000011\\
-0.04	-3.20000000009202e-06\\
0.41	0.00344604999999998\\
0.63	0.0125023500000001\\
0.81	0.02657205\\
0.97	0.0456336500000001\\
1.12	0.0702463999999998\\
1.26	0.1000188\\
1.4	0.1372\\
1.53	0.17907885\\
1.66	0.2287148\\
1.78	0.2819876\\
1.9	0.34295\\
2.02	0.4121204\\
2.14	0.4900172\\
2.21	0.53969305\\
};
\addlegendentry{$\varphi_2(x)=0$}

\addplot [color=red, draw=none, mark=x, mark options={solid, red}, forget plot]
  table[row sep=crcr]{%
-0.523598775598299	-0.785398163397448\\
};
\end{axis}

\begin{axis}[%
width=0.45\fwidth,
height=\fheight,
at={(0.56\fwidth,0\fheight)},
scale only axis,
xmin=0,
xmax=4,
xlabel style={font=\color{white!15!black}},
xlabel={$t$},
ymin=-1,
ymax=2.5,
ylabel style={font=\color{white!15!black}},
ylabel={$x(t)$},
axis background/.style={fill=white},
xmajorgrids,
ymajorgrids,
legend style={at={(0.605,0.722)}, anchor=south west, legend cell align=left, align=left, draw=white!15!black,nodes={scale=0.70, transform shape}},
xlabel style={font={\scriptsize}},ylabel style={font=\scriptsize},  ylabel shift={-0cm},ticklabel style={font=\scriptsize}
]
\addplot [color=mycolor1, line width = 1.5pt]
  table[row sep=crcr]{%
0	2.0943951023932\\
0.666666666666667	1.42786881360382\\
1.33333333333333	0.760781013305388\\
2	0.0956585032885018\\
2.0996162275037	1.77635683940025e-15\\
2.66666666666667	-2.39808173319034e-14\\
3.25858562862696	-0.0027151213003318\\
3.39055987072614	-0.0125242100851137\\
3.54680037719366	-0.0317457347872949\\
3.60051908137676	-0.0400343690957277\\
3.64905333232829	-0.0568608969680877\\
3.69758758327982	-0.0852830105966209\\
3.76425424994648	-0.143226636405179\\
3.83092091661315	-0.223048559094138\\
3.89758758327982	-0.324748778663499\\
3.96425424994648	-0.448327295113262\\
4	-0.523598775606894\\
};
\addlegendentry{$x_1(t)$}

\addplot [color=mycolor2, line width = 1.5pt]
  table[row sep=crcr]{%
0	1.0471975511966\\
0.666666666666667	0.379539275607538\\
0.979551333326136	0.0693025788596549\\
1.11068114398771	0.0475989884210266\\
1.25453627124502	0.0296035657623843\\
1.33333333333333	0.0220165365546219\\
1.50698686768852	0.0101139230901248\\
1.6899839868844	0.0033041396473692\\
1.94740945071528	0.000161870994404012\\
2	4.37663922392417e-05\\
2.66666666666667	1.67865721323324e-13\\
3.02790863960859	-0.00307344828435063\\
3.07102298098143	-0.0129123340748132\\
3.08539442810571	-0.0176970799917795\\
3.12869222823602	-0.0366614344767857\\
3.17199002836634	-0.0624565802788686\\
3.21528782849665	-0.0950825173980263\\
3.25858562862696	-0.134539245834261\\
3.33500338298408	-0.220837831860786\\
3.4255893572404	-0.329940000322817\\
3.50152957421659	-0.413577557626418\\
3.60051908137676	-0.511885661863365\\
3.67332045780405	-0.576446983254587\\
3.76425424994648	-0.647875249188915\\
3.83092091661315	-0.693740348743471\\
3.89758758327982	-0.73410512580911\\
4	-0.785398163404992\\
};
\addlegendentry{$x_2(t)$}

\end{axis}

\begin{axis}[%
width=1.227\fwidth,
height=1.227\fheight,
at={(-0.16\fwidth,-0.135\fheight)},
scale only axis,
xmin=0,
xmax=1,
ymin=0,
ymax=1,
axis line style={draw=none},
ticks=none,
axis x line*=bottom,
axis y line*=left,
legend style={legend cell align=left, align=left, draw=white!15!black},
xlabel style={font={\scriptsize}},ylabel style={font=\scriptsize},  ylabel shift={-0cm},ticklabel style={font=\scriptsize}
]
\end{axis}
\end{tikzpicture}%
}
	\caption{A solution $x(t)$ to the OCP \eqref{eq:ocp_step_eq}.}
	\label{fig:step_eq_solution}
	\centering
	{ \pgfplotsset{compat=1.13}
\setlength{\fwidth}{9.5cm}
\setlength{\fheight}{2.45cm}
\definecolor{mycolor1}{rgb}{0.00000,0.44700,0.74100}%
\definecolor{mycolor2}{rgb}{0.85000,0.32500,0.09800}%
\begin{tikzpicture}

\begin{axis}[%
width=0.45\fwidth,
height=\fheight,
at={(0\fwidth,0\fheight)},
scale only axis,
xmin=0,
xmax=4,
xlabel style={font=\color{white!15!black}},
xlabel={$t$},
ymin=-3.51310036332557,
ymax=1,
ylabel style={font=\color{white!15!black}},
ylabel={$v(t)$},
axis background/.style={fill=white},
		legend style={at={(0.055,0.122)}, anchor=south west, legend cell align=left, align=left, draw=white!15!black,nodes={scale=0.70, transform shape}},
xlabel style={font={\scriptsize}},ylabel style={font=\scriptsize},  ylabel shift={-0cm},ticklabel style={font=\scriptsize}
]
\addplot [color=mycolor1, line width=1.5pt]
  table[row sep=crcr]{%
0	0\\
2.23871045233656	-0.00438933399497099\\
2.66666666666667	-0.0235834833783413\\
3.33333333333333	0.0880169290370247\\
4	-3.19372760302325\\
};
\addlegendentry{$v_1(t)$}

\addplot [color=mycolor2, line width=1.5pt]
  table[row sep=crcr]{%
0	0\\
0.87320713084938	0.00163333899250873\\
1.33333333333333	0.0118993070651454\\
2	-0.0446224014943\\
2.66666666666667	0.166590298912054\\
3.33333333333333	-2.26252481391827\\
4	-1.43747644058087\\
};
\addlegendentry{$v_2(t)$}

\addplot [color=black, dashed, forget plot]
  table[row sep=crcr]{%
0.666666666666667	-4\\
0.666666666666667	4\\
};
\addplot [color=black, dashed, forget plot]
  table[row sep=crcr]{%
1.33333333333333	-4\\
1.33333333333333	4\\
};
\addplot [color=black, dashed, forget plot]
  table[row sep=crcr]{%
2	-4\\
2	4\\
};
\addplot [color=black, dashed, forget plot]
  table[row sep=crcr]{%
2.66666666666667	-4\\
2.66666666666667	4\\
};
\addplot [color=black, dashed, forget plot]
  table[row sep=crcr]{%
3.33333333333333	-4\\
3.33333333333333	4\\
};
\addplot [color=black, dashed, forget plot]
  table[row sep=crcr]{%
4	-4\\
4	4\\
};
\end{axis}

\begin{axis}[%
width=0.45\fwidth,
height=\fheight,
at={(0.55\fwidth,0\fheight)},
scale only axis,
xmin=0,
xmax=4,
xlabel style={font=\color{white!15!black}},
xlabel={$t$},
ymin=-2.70743923894972,
ymax=1,
ylabel style={font=\color{white!15!black}},
ylabel={${u}(t)$},
axis background/.style={fill=white},
legend style={at={(0.055,0.122)}, anchor=south west, legend cell align=left, align=left, draw=white!15!black,nodes={scale=0.70, transform shape}},
xlabel style={font={\scriptsize}},ylabel style={font=\scriptsize},  ylabel shift={-0cm},ticklabel style={font=\scriptsize}
]
\addplot[const plot, color=mycolor1, line width=1.5pt] table[row sep=crcr] {%
0	0.000315850223816749\\
0.666666666666667	-0.00157925111908375\\
1.33333333333333	0.00600115425252401\\
2	-0.0224253658910127\\
2.66666666666667	0.0837003093115247\\
3.33333333333333	-2.4613083990452\\
4	-2.4613083990452\\
};
\addlegendentry{${u}_1(t)$}

\addplot[const plot, color=mycolor2, line width=1.5pt] table[row sep=crcr] {%
0	-0.00223112007470894\\
0.666666666666667	0.0111556003735678\\
1.33333333333333	-0.0423912814195839\\
2	0.158409525304765\\
2.66666666666667	-1.82183633462274\\
3.33333333333333	0.618786280003046\\
4	0.618786280003046\\
};
\addlegendentry{${u}_2(t)$}

\addplot [color=black, dashed, forget plot]
  table[row sep=crcr]{%
0.666666666666667	-3.24892708673967\\
0.666666666666667	3.24892708673967\\
};
\addplot [color=black, dashed, forget plot]
  table[row sep=crcr]{%
1.33333333333333	-3.24892708673967\\
1.33333333333333	3.24892708673967\\
};
\addplot [color=black, dashed, forget plot]
  table[row sep=crcr]{%
2	-3.24892708673967\\
2	3.24892708673967\\
};
\addplot [color=black, dashed, forget plot]
  table[row sep=crcr]{%
2.66666666666667	-3.24892708673967\\
2.66666666666667	3.24892708673967\\
};
\addplot [color=black, dashed, forget plot]
  table[row sep=crcr]{%
3.33333333333333	-3.24892708673967\\
3.33333333333333	3.24892708673967\\
};
\addplot [color=black, dashed, forget plot]
  table[row sep=crcr]{%
4	-3.24892708673967\\
4	3.24892708673967\\
};
\end{axis}

\begin{axis}[%
width=1.227\fwidth,
height=1.227\fheight,
at={(-0.16\fwidth,-0.135\fheight)},
scale only axis,
xmin=0,
xmax=1,
ymin=0,
ymax=1,
axis line style={draw=none},
ticks=none,
axis x line*=bottom,
axis y line*=left,
legend style={legend cell align=left, align=left, draw=white!15!black},
xlabel style={font={\scriptsize}},ylabel style={font=\scriptsize},  ylabel shift={-0cm},ticklabel style={font=\scriptsize}
]
\end{axis}
\end{tikzpicture}%
}
	\caption{The left plot shows the  solution trajectories for $v(t)$. 
		The right plot shows optimal controls obtained via the FESD discretization of the OCP \eqref{eq:ocp_step_eq}. The vertical dashed lines highlight the control discretization grid.}
	\label{fig:step_eq_controls} 
\end{figure}
{
In our numerical experiments, we do not provide an elaborate initial guess. 
For the states, we simply use the initial value for all stages and set the controls to zero.
The homotopy approaches implemented in \texttt{NOSNOC} are quite robust and usually succeed in finding a good solution, as shown in the benchmark in~\cite{Nurkanovic2023e}.
}
\subsection{A numerical optimal control example}\label{sec:step_equilb_example}
The following optimal control problem demonstrates several features developed in this paper: FESD for Cartesian Products of Filippov systems (Sec. \ref{sec:cartesian_filippov}), handling multiple sliding modes, step equilibration (Sec. \ref{sec:step_equilibration}) and equidistant control discretization grids (Sec. \ref{sec:multiple_shooting_grids}). 
We also include a benchmark where we compare the accuracy and computational time of the standard and FESD methods.
Regard the following OCP with $q\in \R^2, v\in \R^2, u\in \R^2$ and $x = (q,v)$:
\begin{subequations}\label{eq:ocp_step_eq}
	\begin{align}
		\min_{x(\cdot),u(\cdot)} \quad & \int_{0}^{4} u(t)^\top u(t) + v(t)^\top v(t) \, \mathrm{d} t +\rho\; 
		\|q(4) - q_{\mathrm{final}}\|_1 \\
		\textrm{s.t.} \quad & x(0) = (\frac{2\pi}{3},\frac{\pi}{3},0,0),\\
		&\dot{x}(t) = \begin{bmatrix}
			-\mathrm{sign}(c(x(t))) + v(t)\\
			u(t)
		\end{bmatrix}
		 ,\quad t \in [0,4],\\
		& -2 e \leq v(t) \leq  2e,\quad t \in [0,4],\\
		& -10 e \leq u(t) \leq  10e,\quad t \in [0,4].
	\end{align}
\end{subequations}	
where $\rho = 10^3$, $\varphi_1(x)  = q_1+0.15q_2^2$, $\varphi_2(x)  = -0.05q_1^3+q_2$ and the function
$c(x)=(\varphi_1(x),\varphi_2(x))$ defines the region boundaries. 
It can be seen that for $v(t) = 0$ the vector fields of $q$ point in all regions towards the origin in the $(q_1,q_2)$ plane. 
Settings the control functions to zero, this results in trajectories going to the origin with sliding modes on the surfaces of discontinuity defined by $c(x)=0$. 
On the other hand, by increasing the value of $v(t)$ via the control functions $u(t)$, the vector fields can change their direction and sliding modes can be left or not achieved at all.

The goal in the OCP is to reach  $q_{\mathrm{final}} = (-\frac{\pi}{6},-\frac{\pi}{4})$ with a minimum control effort. 
The trajectory has to:
(1) first reach $\varphi_1(x)  = 0$;
(2) slide towards $\mathcal{M} = \{q \in \R^2 \mid c(x) = 0\}$;
(3) stay there for some time;
(4) exit $\mathcal{M}$ and slide on $\varphi_2(x) = 0$;
(5) and then leave the sliding mode as late as possible to reach $q_{\mathrm{final}}$.
The seemingly simple example comprises several difficult switching cases in its solution.
The described solution is illustrated in Figure \ref{fig:step_eq_solution} and the optimal controls $u(t)$ and state $v(t)$ are depicted in Figure \ref{fig:step_eq_controls}. 
The OCP is discretized with the FESD Radau-IIA method of order 3 with $\Nstg =2$ and $\Nctrl = 6$ control intervals with $\NFE = 6$ finite elements on every control interval. 
This system is transformed as described in Section \ref{sec:cartesian_filippov} with $\Nsubsys  = 2$, where $c_{i,1}(x) = \varphi_i(x), i = 1,2$.  The functions $g_{i,j}(x)\in \R^2, i,j = 1,2$ are computed via Eq. \eqref{eq:indicator_func_formula}. 

This solution comprises four switches and different sliding modes on nonlinear manifolds including: twice sliding on co-dimension 1 manifolds and once on the co-dimension 2 manifold $\mathcal{M}$, and additionally leaving the sliding mode twice. 
One can see in Figure \ref{fig:step_eq_controls} that the control is discretized on an equidistant grid despite the variable length of the finite elements. This illustrates the multiple shooting-type discretizations described in Section \ref{sec:multiple_shooting_grids}.

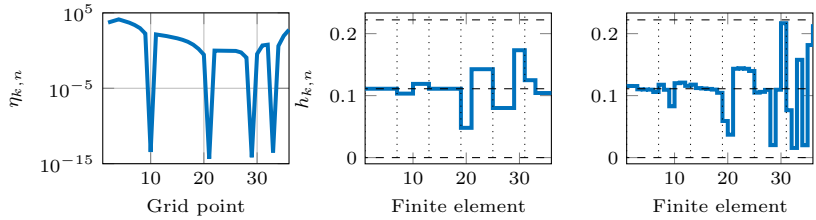
\begin{figure}[t]
	\centering
	{ \pgfplotsset{compat=1.13}
\setlength{\fwidth}{10cm}
\setlength{\fheight}{2.0cm}
\definecolor{mycolor1}{rgb}{0.00000,0.44700,0.74100}%
\begin{tikzpicture}

\begin{axis}[%
width=0.245\fwidth,
height=\fheight,
at={(0\fwidth,0\fheight)},
scale only axis,
xmin=1,
xmax=36,
xlabel style={font=\color{white!15!black}},
xlabel={Grid point},
ymode=log,
ymin=1e-15,
ymax=100000,
yminorticks=true,
ylabel style={font=\color{white!15!black}},
ylabel={$\eta_{k,n}$},
axis background/.style={fill=white},
xmajorgrids,
ymajorgrids,
yminorgrids,
legend style={legend cell align=left, align=left, draw=white!15!black},
xlabel style={font={\scriptsize}},ylabel style={font=\scriptsize},  ylabel shift={-0cm},ticklabel style={font=\scriptsize}
]
\addplot [color=mycolor1, line width=1.5pt]
  table[row sep=crcr]{%
2	5151.66590952098\\
4	13747.6706277476\\
6	4966.7281540256\\
7	2444.51397837466\\
8	941.152763230192\\
9	179.463761762105\\
10	3.54045128751418e-14\\
11	138.55172529678\\
13	83.261683924396\\
15	42.1748260888314\\
16	26.6389396370465\\
17	14.5699473560956\\
18	6.03876984843044\\
19	1.67855563186402\\
20	0.284703508473581\\
21	4.4922526470985e-15\\
22	0.995929613567003\\
26	0.897036771938966\\
27	0.562478104180656\\
28	0.115799987583859\\
29	7.03069454698289e-15\\
30	0.353376296315169\\
31	4.76866335274363\\
32	6.00030317743819\\
33	3.08852200418892e-14\\
34	4.08164782210334\\
35	88.7864984244316\\
36	577.669614437855\\
};

\end{axis}

\begin{axis}[%
width=0.245\fwidth,
height=\fheight,
at={(0.345\fwidth,0\fheight)},
scale only axis,
xmin=1,
xmax=36,
xlabel style={font=\color{white!15!black}},
xlabel={Finite element},
ymin=-0.01,
ymax=0.233333333333333,
ylabel style={font=\color{white!15!black}},
ylabel={$h_{k,n}$},
axis background/.style={fill=white},
legend style={legend cell align=left, align=left, draw=white!15!black},
xlabel style={font={\scriptsize}},ylabel style={font=\scriptsize},  ylabel shift={-0cm},ticklabel style={font=\scriptsize}
]
\addplot[const plot, color=mycolor1, line width=1.5pt] table[row sep=crcr] {%
1	0.111111111111114\\
7	0.103270232091354\\
10	0.118951990130867\\
13	0.111111111111114\\
19	0.0480291904465986\\
21	0.142652071443365\\
25	0.0800421990462539\\
29	0.173248935240828\\
31	0.124906863068965\\
33	0.104213235132185\\
36	0.104213235132185\\
};

\addplot [color=black, dotted]
  table[row sep=crcr]{%
1	0\\
1	0.257666666666667\\
};

\addplot [color=black, dotted]
  table[row sep=crcr]{%
7	0\\
7	0.257666666666667\\
};

\addplot [color=black, dotted]
  table[row sep=crcr]{%
13	0\\
13	0.257666666666667\\
};

\addplot [color=black, dotted]
  table[row sep=crcr]{%
19	0\\
19	0.257666666666665\\
};

\addplot [color=black, dotted]
  table[row sep=crcr]{%
25	0\\
25	0.257666666666665\\
};

\addplot [color=black, dotted]
  table[row sep=crcr]{%
31	0\\
31	0.257666666666665\\
};

\addplot [color=black, dashed]
  table[row sep=crcr]{%
0	0\\
36	0\\
};

\addplot [color=black, dashed]
  table[row sep=crcr]{%
0	0.222222222222221\\
36	0.222222222222221\\
};

\addplot [color=black, dashed]
  table[row sep=crcr]{%
0	0.111111111111114\\
36	0.111111111111114\\
};

\end{axis}

\begin{axis}[%
width=0.245\fwidth,
height=\fheight,
at={(0.689\fwidth,0\fheight)},
scale only axis,
xmin=1,
xmax=36,
xlabel style={font=\color{white!15!black}},
xlabel={Finite element},
ymin=-0.01,
ymax=0.233333333333333,
axis background/.style={fill=white},
legend style={legend cell align=left, align=left, draw=white!15!black},
xlabel style={font={\scriptsize}},ylabel style={font=\scriptsize},  ylabel shift={-0cm},ticklabel style={font=\scriptsize}
]
\addplot[const plot, color=mycolor1, line width=1.5pt] table[row sep=crcr] {%
	1	0.115592542844105\\
	2	0.115514097718524\\
	3	0.110187895431544\\
	4	0.109903134316554\\
	5	0.109637179271083\\
	6	0.105831817084855\\
	7	0.117581039634452\\
	8	0.109436555649431\\
	9	0.0827931009901803\\
	10	0.120102909032944\\
	11	0.121047250062027\\
	12	0.115705811297637\\
	13	0.117777491651616\\
	14	0.112194677019374\\
	15	0.111522832744519\\
	16	0.110846651725502\\
	17	0.109768763255772\\
	18	0.104556250269887\\
	19	0.0590776986582568\\
	20	0.0369806822349403\\
	21	0.143368519573535\\
	22	0.144189560621378\\
	23	0.143292895347471\\
	24	0.139757310231083\\
	25	0.104805496224401\\
	26	0.106081403059662\\
	27	0.109281896297958\\
	28	0.0201194638519695\\
	29	0.109528551168196\\
	30	0.216849856064478\\
	31	0.0764712791051778\\
	32	0.015555261798248\\
	33	0.157787185234511\\
	34	0.0199141498535056\\
	35	0.181606390879551\\
	36	0.215332399795678\\
};

\addplot [color=black, dotted]
  table[row sep=crcr]{%
1	0\\
1	0.257666666666667\\
};

\addplot [color=black, dotted]
  table[row sep=crcr]{%
7	0\\
7	0.257666666666667\\
};

\addplot [color=black, dotted]
  table[row sep=crcr]{%
13	0\\
13	0.257666666666667\\
};

\addplot [color=black, dotted]
  table[row sep=crcr]{%
19	0\\
19	0.257666666666665\\
};

\addplot [color=black, dotted]
  table[row sep=crcr]{%
25	0\\
25	0.257666666666665\\
};

\addplot [color=black, dotted]
  table[row sep=crcr]{%
31	0\\
31	0.257666666666665\\
};

\addplot [color=black, dashed]
  table[row sep=crcr]{%
0	0\\
36	0\\
};

\addplot [color=black, dashed]
  table[row sep=crcr]{%
0	0.222222222222221\\
36	0.222222222222221\\
};

\addplot [color=black, dashed]
  table[row sep=crcr]{%
0	0.111111111111114\\
36	0.111111111111114\\
};

\end{axis}

\begin{axis}[%
width=1.227\fwidth,
height=1.227\fheight,
at={(-0.185\fwidth,-0.135\fheight)},
scale only axis,
xmin=0,
xmax=1,
ymin=0,
ymax=1,
axis line style={draw=none},
ticks=none,
axis x line*=bottom,
axis y line*=left,
legend style={legend cell align=left, align=left, draw=white!15!black},
xlabel style={font={\scriptsize}},ylabel style={font=\scriptsize},  ylabel shift={-0cm},ticklabel style={font=\scriptsize}
]
\end{axis}
\end{tikzpicture}%
}
	\caption{The left plot depicts the switching indicator function $\eta(\cdot)$ at the solution of the OCP \eqref{eq:ocp_step_eq}. 
		The middle plot show the step size $h_{n,k}$ with step equilibration. 
		The right plot shows the step sizes without step equilibration. 
		The horizontal dashed lines correspond to the minimum, maximum, and nominal step size, the vertical dotted lines correspond to control interval boundaries.}
	\label{fig:step_eq_step_size}
	\vspace{-0.4cm}
\end{figure}

In the first experiment, we demonstrate the effects of step equilibration for a solution of the OCP \eqref{eq:ocp_step_eq}. The left plot in Figure \ref{fig:step_eq_step_size} depicts the indicator function $\eta(\cdot)$. Clearly, the function is only zero if a switch occurs, cf. the right plot of Figure \ref{fig:step_eq_solution}. The resulting step sizes $h_{k,n}$ with and without step equilibration are depicted in the middle and right plots, respectively. For the right plot we discard the step equilibration conditions $G_{\mathrm{eq}}(\mathbf{h},\mathbf{\Theta},\mathbf{\Lambda},T) = 0$.
Obviously, without them, the optimizer varies the step size in a somewhat random way. 
On the other hand, with step equilibration, we obtain a piecewise equidistant grid, where the step size changes only when a switch occurs.

In the second experiment, we compare the accuracy of an OCP solution obtained with the standard and FESD method as a function of the CPU time. 
We take the optimal controls and perform a high-accuracy simulation of the system dynamics in \eqref{eq:ocp_step_eq}, which we denote by $x_{\mathrm{int}}(t)$. 
As a metric, we take the terminal constraint satisfaction of the high accuracy solution, i.e.,  $E(T) = \| x_{\mathrm{int}}(T) - x_{\mathrm{final}}\|$. 
We set $\Nctrl=6$, and vary the number of finite elements per stage $\NFE$ from 1 to 7 and the number of stage points $\Nstg$ from 1 to 4. 
The experiment is performed for the following RK methods: Radau-IIA, Gauss-Legendre, Lobatto-IIIC, and Explicit-RK. 
 
\begin{figure}[t]
	\centering
	{ \input{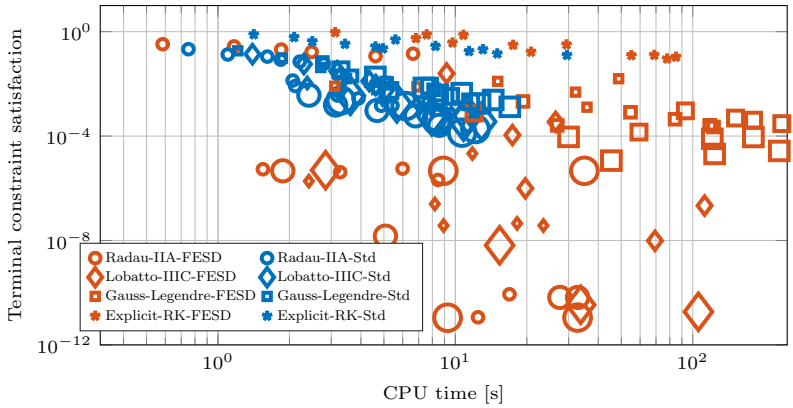}}
	\caption{Terminal constraint satisfaction vs. CPU time for the Standard and FESD method. The size of the marker indicates the number of stage points, the smallest corresponds to $\Nstg=1$ and the largest to $\Nstg = 4$.}
	\label{fig:pareto_lots}
	\vspace{-0.45cm}					
\end{figure}

For Radau-IIA-FESD and Lobatto-IIIC-FESD we solved the arising MPCC with an elastic mode homotopy approach \cite{Anitescu2007}. In the other scenarios, we were not able to solve all problems to convergence with the elastic mode approach. 
Therefore, the MPCC was solved with a relaxation homotopy approach, cf. \cite{Nurkanovic2022b} for implementation details. 
The second approach is slightly slower than elastic mode, but more robust, and all problems were solved successfully. 
The terminal error as a function of the total CPU time is given in Figure \ref{fig:pareto_lots}.

We can draw several conclusions from the experiments. Clearly, the FESD method outperforms the standard approach in all experiments. 
For example, for a CPU time of $\approx1$ second FESD achieves five orders of magnitude more accurate solutions than the standard {time-stepping} approach. 
A better solution than the most accurate one of the standard approaches can be achieved with FESD by an order of magnitude faster CPU time.
The Radau-IIA and Lobatto-IIIC methods are the most efficient ones in this benchmark, whereas the Gauss-Legendre and Explicit-RK methods perform poorly. 
This is no surprise, since the solution trajectories contain sliding mode arcs which require solving nonlinear DAE of index 2. 
Radau-IIA and Lobatto-IIIC usually perform well and have good theoretical properties for higher index DAE, whereas Gauss-Legendre and Explicit RK even lose the high accuracy orders that they have for ODE \cite{Hairer1991}. 
The source code of all examples is available in the repository of the open source tool \texttt{NOSNOC} \cite{Nurkanovic2022c}.
The same repository contains a few additional examples where FESD is used, including systems with state jumps that are transformed via time-freezing into PSS \cite{Nurkanovic2022a,Nurkanovic2022b} and a nonsmooth mechanics simulation problem from \cite{Stewart1996a}, among others.

\section{Summary}\label{sec:summary}
This paper presents a method that enables direct optimal control of a broad class of switched systems with high simulation accuracy, called Finite Elements with Switch Detection (FESD). We build a solid theoretical foundation and prove that FESD has the same accuracy as the underlying RK method has for smooth differential equations and that it delivers exact numerical sensitivities. An implementation of the FESD method described in this paper is available in the recently introduced open-source package \texttt{NOSNOC} \cite{Nurkanovic2022c}. Compared to standard discretizations, FESD achieves for a similar CPU time usually several orders of magnitude more accurate solutions.

A key advantage of the new approach for direct optimal control is that no guessing of the number or order of switches needs to be done. FESD can treat multiple or simultaneous switches and sliding modes. 
With time-freezing \cite{Nurkanovic2021,Nurkanovic2022a,Nurkanovic2023a} many nonsmooth systems with state jumps can be recast into a piecewise smooth system. This allows us to treat many classes of nonsmooth systems in direct optimal control in a unified way. 

In future work, one should relax some of the possibly restrictive assumptions in our theoretical analysis. We aim to extend FESD to other transformations of piecewise smooth systems into dynamic complementarity systems, e.g., the via the step function approach \cite{Acary2014,Dieci2011}.
Some further open questions to be answered are: Are all limit points of the solution approximations candidates $\hat{x}_h(t)$  indeed solutions to the Filippov DI \eqref{eq:FilippovDI}? Do unique Filippov solutions to a given problem also imply a unique solution to the corresponding FESD problem \eqref{eq:fesd_compact}?


\appendix

\section{Auxiliary results needed for Theorem \ref{th:integration_order}}\label{app:switch_detection_lcp}
\normalsize
This lemma is about the active sets of perturbed LCP, where strict complementarity holds at a solution.
\begin{lemma1}[{\cite{Stewart1990b}[Lemma A.2]}]\label{lem:lcp_active_sets}
	Suppose that all entries of $M$ are positive in \eqref{eq:generic_lcp} and all solutions of $\mathrm{LCP}(M,q)$ are strongly stable. 
	If $\hat{M}_n \to M$, $\hat{q}_n \to q$, then $\mathrm{SOL}(\hat{M}_n,\hat{q}_n) \to \mathrm{SOL}(M,q)$, as $n\to\infty$, in the Hausdorff metric. Moreover, if $(w,\theta)\in \mathrm{SOL}(M,q)$, such that $w+\theta>0$, then there is a $(\hat{w}_n,\hat{\theta}_n) \in \mathrm{SOL}(\hat{M}_n,\hat{q}_n)$ for sufficiently large $n$ such that $\{i \mid \hat{\theta}_{n,i}>0\}  = \{i \mid {\theta}_{i}>0\}$.
\end{lemma1}

\bibliographystyle{abbrv}


\begin{thebibliography}{10}

\bibitem{Nurkanovic2022c}
{NOSNOC}.
\newblock https://github.com/nurkanovic/nosnoc, 2022.

\bibitem{Acary2010}
Vincent Acary, Olivier Bonnefon, and Bernard Brogliato.
\newblock {\em Nonsmooth modeling and simulation for switched circuits},
  volume~69.
\newblock Springer Science \& Business Media, 2010.

\bibitem{Acary2008}
Vincent Acary and Bernard Brogliato.
\newblock {\em Numerical methods for nonsmooth dynamical systems: applications
  in mechanics and electronics}.
\newblock Springer Science \& Business Media, 2008.

\bibitem{Acary2010b}
Vincent Acary and Bernard Brogliato.
\newblock Implicit {E}uler numerical scheme and chattering-free implementation
  of sliding mode systems.
\newblock {\em Systems \& Control Letters}, 59(5):284--293, 2010.

\bibitem{Acary2014}
Vincent Acary, Hidde De~Jong, and Bernard Brogliato.
\newblock Numerical simulation of piecewise-linear models of gene regulatory
  networks using complementarity systems.
\newblock {\em Physica D: Nonlinear Phenomena}, 269:103--119, 2014.

\bibitem{Albersmeyer2010b}
Jan Albersmeyer.
\newblock {\em Adjoint-based algorithms and numerical methods for sensitivity
  generation and optimization of large scale dynamic systems}.
\newblock PhD thesis, University of Heidelberg, 2010.

\bibitem{Andersson2019}
Joel A.~E. Andersson, Joris Gillis, Greg Horn, James~B Rawlings, and Moritz
  Diehl.
\newblock {CasADi} -- a software framework for nonlinear optimization and
  optimal control.
\newblock {\em Mathematical Programming Computation}, 11(1):1--36, 2019.

\bibitem{Anitescu2007}
Mihai Anitescu, Paul Tseng, and Stephen~J. Wright.
\newblock Elastic-mode algorithms for mathematical programs with equilibrium
  constraints: global convergence and stationarity properties.
\newblock {\em Mathematical {P}rogramming}, 110(2):337--371, 2007.

\bibitem{Ban2012}
Xuegang~Jeff Ban, Jong-Shi Pang, Henry~X Liu, and Rui Ma.
\newblock Continuous-time point-queue models in dynamic network loading.
\newblock {\em Transportation Research Part B: Methodological}, 46(3):360--380,
  2012.

\bibitem{Baumrucker2009}
Brian~T. Baumrucker and Lorenz~T. Biegler.
\newblock {M}{P}{E}{C} strategies for optimization of a class of hybrid dynamic
  systems.
\newblock {\em Journal of Process Control}, 19(8):1248--1256, 2009.

\bibitem{Bemporad1999b}
Alberto Bemporad and Manfred Morari.
\newblock {C}ontrol of systems integrating logic, dynamics, and constraints.
\newblock {\em Automatica}, 35(3):407--427, 1999.

\bibitem{Bock1984}
Hans~G. Bock and K.~J. Plitt.
\newblock A multiple shooting algorithm for direct solution of optimal control
  problems.
\newblock In {\em Proceedings of the IFAC World Congress}, pages 242--247.
  Pergamon Press, 1984.

\bibitem{Bock2018}
Hans~Georg Bock, Christian Kirches, Andreas Meyer, and Andreas Potschka.
\newblock Numerical solution of optimal control problems with explicit and
  implicit switches.
\newblock {\em Optimization Methods and Software}, 33(3):450--474, 2018.

\bibitem{Brogliato2020}
Bernard Brogliato and Aneel Tanwani.
\newblock Dynamical systems coupled with monotone set-valued operators:
  {F}ormalisms, applications, well-posedness, and stability.
\newblock {\em SIAM Review}, 62(1):3--129, 2020.

\bibitem{Cojocaru2008}
Monica-Gabriela Cojocaru.
\newblock Dynamic equilibria of group vaccination strategies in a heterogeneous
  population.
\newblock {\em Journal of Global Optimization}, 40(1):51--63, 2008.

\bibitem{Colombo2020}
Giovanni Colombo, Boris~S Mordukhovich, and Dao Nguyen.
\newblock Optimization of a perturbed sweeping process by constrained
  discontinuous controls.
\newblock {\em SIAM Journal on Control and Optimization}, 58(4):2678--2709,
  2020.

\bibitem{Dieci2011}
Luca Dieci and Luciano Lopez.
\newblock Sliding motion on discontinuity surfaces of high co-dimension. a
  construction for selecting a filippov vector field.
\newblock {\em Numerische Mathematik}, 117(4):779--811, 2011.

\bibitem{Dontchev1992}
Asen Dontchev and Frank Lempio.
\newblock Difference methods for differential inclusions: {A} survey.
\newblock {\em SIAM review}, 34(2):263--294, 1992.

\bibitem{Dontchev2014}
Asen~L. Dontchev and Tyrrell~R. Rockafellar.
\newblock {\em Implicit Functions and Solution Mappings: A View from
  Variational Analysis}.
\newblock Springer, 2014.

\bibitem{Facchinei2003}
Francisco Facchinei and Jong-Shi Pang.
\newblock {\em Finite-dimensional variational inequalities and complementarity
  problems}, volume 1-2.
\newblock Springer-Verlag, 2003.

\bibitem{Filippov1962}
Aleksei~Fedorovich Filippov.
\newblock On certain questions in the theory of optimal control.
\newblock {\em Journal of the Society for Industrial and Applied Mathematics,
  Series A: Control}, 1(1):76--84, 1962.

\bibitem{Filippov1964}
Aleksei~Fedorovich Filippov.
\newblock {D}ifferential {E}quations with discontinuous right hand side.
\newblock {\em AMS Transl.}, 42:199--231, 1964.

\bibitem{Filippov1988}
Alexei~F. Filippov.
\newblock {\em Differential Equations with Discontinuous Righthand Sides:
  Control Systems}, volume~18.
\newblock Springer Science \& Business Media, 1988.

\bibitem{Guo2016}
Lei Guo and Jane~J. Ye.
\newblock Necessary optimality conditions for optimal control problems with
  equilibrium constraints.
\newblock {\em SIAM Journal on Control and Optimization}, 54(5):2710--2733,
  2016.

\bibitem{Hairer1991}
E.~Hairer and G.~Wanner.
\newblock {\em {S}olving {O}rdinary {D}ifferential {E}quations {II} -- {S}tiff
  and {D}ifferential-{A}lgebraic {P}roblems}.
\newblock Springer, Berlin Heidelberg, 2nd edition, 1991.

\bibitem{Hairer1993}
Ernst Hairer, Syvert~.P. N{\o}rsett, and Gerhard Wanner.
\newblock {\em {S}olving {O}rdinary {D}ifferential {E}quations {I}}.
\newblock Springer Series in Computational Mathematics. Springer, Berlin, 2nd
  edition, 1993.

\bibitem{Hall2022}
Jonas Hall, Armin Nurkanovi\'c, Florian Messerer, and Moritz Diehl.
\newblock A sequential convex programming approach to solving quadratic
  programs and optimal control problems with linear complementarity
  constraints.
\newblock {\em IEEE Control Systems Letters}, 6:536--541, 2022.

\bibitem{Hauswirth2021}
Adrian Hauswirth, Saverio Bolognani, Gabriela Hug, and Florian D{\"o}rfler.
\newblock Optimization {A}lgorithms as {R}obust {F}eedback {C}ontrollers.
\newblock {\em arXiv preprint arXiv:2103.11329}, 2021.

\bibitem{Kastner1990}
Alois Kastner-Maresch.
\newblock {I}mplicit {R}unge-{K}utta methods for differential inclusions.
\newblock {\em Numerical functional analysis and optimization},
  11(9-10):937--958, 1990.

\bibitem{Katayama2020}
Sotaro Katayama, Masahiro Doi, and Toshiyuki Ohtsuka.
\newblock A moving switching sequence approach for nonlinear model predictive
  control of switched systems with state-dependent switches and state jumps.
\newblock {\em International Journal of Robust and Nonlinear Control},
  30(2):719--740, 2020.

\bibitem{Kirches2006}
Christian Kirches.
\newblock {A} {N}umerical {M}ethod for {N}onlinear {R}obust {O}ptimal {C}ontrol
  with {I}mplicit {D}iscontinuities and an {A}pplication to {P}owertrain
  {O}scillations.
\newblock Diploma thesis, University of Heidelberg, October 2006.

\bibitem{Kirches2022}
Christian Kirches, Jeffrey Larson, Sven Leyffer, and Paul Manns.
\newblock Sequential linearization method for bound-constrained mathematical
  programs with complementarity constraints.
\newblock {\em SIAM Journal on Optimization}, 32(1):75--99, 2022.

\bibitem{Leyffer2006}
Sven Leyffer, Gabriel L{\'o}pez-Calva, and Jorge Nocedal.
\newblock Interior methods for mathematical programs with complementarity
  constraints.
\newblock {\em SIAM Journal on Optimization}, 17(1):52--77, 2006.

\bibitem{Matsaglia1974}
George Matsaglia and George PH~Styan.
\newblock Equalities and inequalities for ranks of matrices.
\newblock {\em Linear and multilinear Algebra}, 2(3):269--292, 1974.

\bibitem{Nurkanovic2023a}
Armin Nurkanovi\'c, Sebastian Albrecht, Bernard Brogliato, and Moritz Diehl.
\newblock The time-freezing reformulation for numerical optimal control of
  complementarity lagrangian systems with state jumps.
\newblock {\em Automatica}, 158:111295, 2023.

\bibitem{Nurkanovic2020}
Armin Nurkanovi\'c, Sebastian Albrecht, and Moritz Diehl.
\newblock Limits of {MPCC} {F}ormulations in {D}irect {O}ptimal {C}ontrol with
  {N}onsmooth {D}ifferential {E}quations.
\newblock In {\em 2020 European Control Conference (ECC)}, pages 2015--2020,
  2020.

\bibitem{Nurkanovic2022a}
Armin Nurkanovi{\'c} and Moritz Diehl.
\newblock Continuous optimization for control of hybrid systems with hysteresis
  via time-freezing.
\newblock {\em IEEE Control Systems Letters}, 6:3182--3187, 2022.

\bibitem{Nurkanovic2022b}
Armin Nurkanovi{\'c} and Moritz Diehl.
\newblock Nosnoc: A software package for numerical optimal control of nonsmooth
  systems.
\newblock {\em IEEE Control Systems Letters}, 6:3110--3115, 2022.

\bibitem{Nurkanovic2023e}
Armin Nurkanovi{\'c}, Anton Pozharskiy, and Moritz Diehl.
\newblock Solving mathematical programs with complementarity constraints
  arising in nonsmooth optimal control.
\newblock {\em arXiv preprint arXiv:2312.11022}, 2023.

\bibitem{Nurkanovic2021}
Armin Nurkanovi\'c, Tommaso Sartor, Sebastian Albrecht, and Moritz Diehl.
\newblock A {T}ime-{F}reezing {A}pproach for {N}umerical {O}ptimal {C}ontrol of
  {N}onsmooth {D}ifferential {E}quations with {S}tate {J}umps.
\newblock {\em IEEE Control Systems Letters}, 5(2):439--444, 2021.

\bibitem{Pontryagin1962a}
Lev~Semenovich Pontryagin.
\newblock {\em The mathematical theory of optimal processes}.
\newblock Wiley, 1962.

\bibitem{Pytlak2021}
Radoslaw Pytlak and Damian Suski.
\newblock Algorithms for optimal control of hybrid systems with sliding motion.
\newblock {\em arXiv preprint arXiv:2101.04754}, 2021.

\bibitem{Ralph2004}
Daniel Ralph and Stephen~J. Wright.
\newblock Some properties of regularization and penalization schemes for mpecs.
\newblock {\em Optimization Methods and Software}, 19(5):527--556, 2004.

\bibitem{Rawlings2017}
James~B. Rawlings, David~Q. Mayne, and Moritz Diehl.
\newblock {\em Model Predictive Control: Theory, Computation, and Design}.
\newblock Nob Hill, 2nd edition, 2017.

\bibitem{Shaikh2007}
M.~Shahid Shaikh and Peter~E. Caines.
\newblock On the hybrid optimal control problem: Theory and algorithms.
\newblock {\em IEEE Transactions on Automatic Control}, 52(9):1587--1603, 2007.

\bibitem{Stewart1990b}
David~E. Stewart.
\newblock A high accuracy method for solving {ODE}s with discontinuous
  right-hand side.
\newblock {\em Numerische Mathematik}, 58(1):299--328, 1990.

\bibitem{Stewart1990c}
David~E. Stewart.
\newblock {\em High accuracy numerical methods for ordinary differential
  equations with discontinuous right-hand side}.
\newblock PhD thesis, University of Queensland, St. Lucia, Queensland 4072,
  Australia, 1990.

\bibitem{Stewart1996a}
David~E. Stewart.
\newblock A numerical method for friction problems with multiple contacts.
\newblock {\em The ANZIAM Journal}, 37(3):288--308, 1996.

\bibitem{Stewart2000}
David~E. Stewart.
\newblock {R}igid-body dynamics with friction and impact.
\newblock {\em SIAM {R}eview}, 42(1):3--39, 2000.

\bibitem{Stewart2011}
David~E. Stewart.
\newblock {\em Dynamics with Inequalities: Impacts and Hard Constraints},
  volume~59.
\newblock SIAM, 2011.

\bibitem{Stewart2010}
David~E. Stewart and Mihai Anitescu.
\newblock Optimal control of systems with discontinuous differential equations.
\newblock {\em Numerische Mathematik}, 114(4):653--695, 2010.

\bibitem{Taubert1981}
Klaus Taubert.
\newblock Converging multistep methods for initial value problems involving
  multivalued maps.
\newblock {\em Computing}, 27(2):123--136, 1981.

\bibitem{Vieira2019}
Alexandre Vieira, Bernard Brogliato, and Christophe Prieur.
\newblock Quadratic optimal control of linear complementarity systems:
  First-order necessary conditions and numerical analysis.
\newblock {\em IEEE Transactions on Automatic Control}, 65(6):2743--2750, 2020.

\bibitem{Waechter2006}
Andreas W\"achter and Lorenz~T. Biegler.
\newblock On the implementation of an interior-point filter line-search
  algorithm for large-scale nonlinear programming.
\newblock {\em Mathematical Programming}, 106(1):25--57, 2006.

\end{thebibliography}
\end{document}